\documentclass{elsart}               

\usepackage{amssymb}
\usepackage{amsmath}
\usepackage{amsfonts}
\usepackage{graphicx}
\usepackage{graphicx}          
\usepackage{color}                               

\newtheorem{theorem}{Theorem}[section]

\newtheorem{definition}{Definition}[section]

\newtheorem{lemma}{Lemma}[section]

\newtheorem{proposition}{Proposition}[section]

\def\3n{\negthinspace \negthinspace \negthinspace }
\def\2n{\negthinspace \negthinspace }
\def\1n{\negthinspace }

\def\ds{\displaystyle}
\def\qq{\qquad}

\def\({\Big (}
\def\){\Big )}
\def\[{\Big[}
\def\]{\Big]}

\begin{document}
\footnotesize
\begin{frontmatter}

\title{A Class of Discrete-time Mean-field Stochastic Linear-quadratic Optimal Control Problems with Financial Application} 
\thanks{Xun Li is with the Department of
Applied Mathematics, The Hong Kong Polytechnic University, Hunghom,
Kowloon, Hong Kong. E-mail: {\tt malixun@polyu.edu.hk}. This author
is supported by the Hong Kong RGC grants 520412 and 15209614.} 

\author[Li]{Xun Li}~~
\ead{malixun@polyu.edu.hk}
%
%
\author[Li]{Allen H. Tai}~~
\ead{allen.tai@polyu.edu.hk}
%
%
\author[Li]{Fei Tian}~~
\ead{tianfei0618@foxmail.com}

\address[Li]{Department of Applied Mathematics, The Hong Kong Polytechnic University, Hong Kong, China}
%
%
%
%

\begin{abstract}                          
This paper is concerned with a discrete-time mean-field stochastic linear-quadratic optimal control problem arose from financial application. Through matrix dynamical optimization method, a group of linear feedback controls is investigated. The problem is then reformulated as an operator stochastic linear-quadratic optimal control problem by a sequence of bounded linear operators over Hilbert space, the optimal control with six algebraic Riccati difference equations is obtained by backward induction. The two above approaches are proved to be coincided by the classical method of completing the square. Finally, after discussing the solution of the problem under multidimensional noises, a financial application example is given.
\end{abstract}

\end{frontmatter}

\section{Introduction\label{1}}

In this paper, we consider a discrete time stochastic linear-quadratic optimal control problem of mean-field type. Here the terms `mean-field' and `linear-quadratic' refer to a dynamic model exhibiting macroscopic behaviour of an attractive mean-field interaction and linear stochastic systems with quadratic performance criterion, respectively.
In the combination of these two issue, the investigation of classical mean-field SDE problem can be traced back to 1960s, when McKean \cite{McKean1966} first discusses a similar connection between a series of Markov processes and certain non-linear parabolic equations. Then, many scientific results are emerging.
Dawson \cite{Dawson1983} investigates the dynamics and fluctuations of mean-field system in the critical condition by adopting approach based on the Papanicolaou et al.'s \cite{Papanicolaou1977} theory for Markov processes.
Dawson and G\"{a}rtner \cite{Dawsont1987} examine the conversion from the large deviations from the McKean-Vlasov limit to a generalization of the theory of Freidlin and Wentzell \cite{Freidlin1984}.
G\"{a}rtner \cite{Gartner1988} systematically give research results for a system of diffusions in a domain range of $\mathbb{R}^d$ with long-range weak interaction.
Similar issues can refer to \cite{Chan1994,Dai1996,Bossy1997}.
Buckdahn et al. \cite{Buckdahn2009} study a special approximation by solution of some decoupled forward-backward equation and find out the convergence speed, and later, they investigate such problem under a more general framework \cite{Buckdahn2009-2}.

In theoretical research, the field of stochastic optimal control has made great progress. Rockafellar and Wets \cite{Rockafellar1990} consider some generalized stochastic linear-quadratic optimal control problems in discrete time.
In a discrete-time system, {\it Riccati difference equation} plays an important role in the synthesis of the optimal control. Beghi and D'alessandro \cite{Beghi1998} derive the optimal control for a discrete-time linear-quadratic problem with control-dependent noise. Moore et al. \cite{Moore1999} consider some partially observed stochastic models where the stochastic disturbances depend on both the states and the controls. Ait Rami et al. \cite{AitRami2002} extend Beghi and D'alessandro's result through allowing the weighting matrices in the cost functional are indefinite. Huang et al. \cite{Huang2008} discuss the problem with an infinite horizon, in which the concepts of stochastic stabilizability and exact observability are introduced.
Attention is also focused on the solution of optimal control with mean-field terms. By introducing the mean terms to the cost functional, the variations of the state process and the control process can be minimized so that they are not too sensitive to random events \cite{Yong2013}. Andersson et al. \cite{Andersson2011} study this problem with functional depend coefficients under the assumption of convex action space, which assumption is consensus in further research. 
Du et al. \cite{Du2011} obtain an existence of the solution and a theorem of comparison for one dimensional mean-field backward stochastic differential equations.
Yong \cite{Yong2013} considers the deterministic coefficient continuous-time mean-field stochastic differential equations with linear-quadratic optimal control problem. 
Elliott et al. \cite{Elliott2013} discuss a discrete-time model for obtaining the optimal control with different methodology for solving the problem that necessary and sufficient conditions for the solvability of the problem are presented.
For recent research, this problem is extended in two aspects: with indefinite weight matrices in the cost functional and with an infinite horizon, respectively \cite{Ni2014,Ni2015,Ni2015-2,Ni2016}.
Sun and Yong \cite{Sun2016} demonstrate that the non-emptiness of the admissible control set for all initial state is equivalent to the $\mathcal{L}^2$-stabilizability of the control system by concerning continuous-time model in an infinite horizon with constant coefficient. Li et al. systematically summarize their recent research results for a linear mean-field stochastic differential equation with a quadratic cost functional, included in \cite{Li2016}. On other mean-field type control problems, readers may refer to literature such as \cite{Buckdahn2011,Andersson2011,Bensoussan2013,Meyer-Brandis2012}. 
On the other hand, mean-field game is also an area closely related to mean-field theory. 
Huang et al. \cite{Huang2006} decompose a class of stochastic games into optimal controls problems and designate the {\it Nash certainty equivalence principle} as property of solvability.
%
%
Bensoussan et al. \cite{Bensoussan2016} study the unique existence of an equilibrium strategies of linear-quadratic MFGs by adjoint equation method. 
For relevant literatures, readers can refer to \cite{Bauso2012,Bensoussan2013,Carmona2013,Carmona2015,Gueant2011}.

Stochastic optimal control theory is widely applied in varies practical problems since Wonham's work in 1968 \cite{Wonham1968}. The development of mathematical mean-field stochastic linear-quadratic optimal control theory has greatly promoted the research of applications in recent works.
Zhou and Yin \cite{Zhou2003} study a continuous-time regime-switching model for portfolio selection, where a Markov chain modulated diffusion formulation is used to model the problem.
%
%
Xie and Wang \cite{Xie2008} investigate mean-variance portfolio selection problems by using general stochastic control technique, where an incomplete market is studied with correlative multiple risky assets and a liability according to a Brownian motion with drift.
By adopting the techniques of \cite{Zhou2003}, Chen et al. \cite{Chen2008} investigate the feasibility, obtain the optimal strategy, delineate the efficient frontier, and establish the associated mutual fund theorem over a continuous-time Markov regime-switching model.
Cui et al. \cite{Cui2014} propose a new mean-field framework that provides a more efficient modelling tool and accurate solution to solve sustainability problems.
Dang et al. \cite{Dang2016} consider Markowitz's problem through method of transforming the problem into an equivalent one with bankruptcy prohibition but without portfolio constraints and then treat by martingale theory.
Literatures can be referred to \cite{Cui2015,Hou2016,Zhang2016,Ziemba2003}.

The problem studied in this paper arises from a practical problem in finance. 
Existing research works in financial applications only consider the investment entity's equity assets, and there is no debt. In order to adapt practical application, the system state is adjusted to two linear stochastic difference equations with several cost functional affected variables. Under a series of necessary and sufficient conditions, the obtained Ricatti equations of the adjusted model are provided and a more general frame from mean-field linear-quadratic controls theory to financial applications.

The remainder of the paper is organized as follows. In Section 2, we give the formulation for the problem. Preliminaries for the analyses are presented in Section 3. In Section 4, the closed-loop optimal control is obtained and then it is represented via Riccati equations. Some financial applications of the problem are presented in Section 5. The paper is then concluded in Section 6.

\section{Problem Formulation\label{2}}

Let $N$ be a positive integer. The system equation is the following set of linear stochastic difference equations with $k\in\{0,1,2,\cdots,N-1\}\equiv \mathbb{N}$,
\begin{eqnarray}\label{1-system-1}
\left\{\begin{array}{l}
x_{k+1}=(A_kx_k+\bar{A}_k\mathbb{E}x_k+B_ku_k+\bar{B}_k\mathbb{E}u_k)\\
\quad +(C_kx_k+\bar{C}_k\mathbb{E}x_k+D_ku_k+\bar{D}_k\mathbb{E}u_k)w_k, \\
y_{k+1}=(F_ky_k+\bar{F}_k\mathbb{E}y_k)+ (G_ky_k+\bar{G}_k\mathbb{E}y_k)v_k, \\
x_0=\zeta^x, y_0=\zeta^y,
\end{array}\right.
\end{eqnarray}
where $x_k, y_k \in \mathbb{R}^n$. $A_k,\bar{A}_k,C_k,\bar{C}_k,F_k,\bar{F}_k,G_k,\bar{G}_k\in \mathbb{R}^{n\times n}$, and $B_k,\bar{B}_k,D_k,\bar{D}_k\in \mathbb{R}^{n\times m}$ are given deterministic matrices. $\mathbb{E}$ is the expectation operator. 
Denote the set  $\{0,1,2,\cdots,$ $ N\}$ by $\bar{\mathbb{N}}$. 
In (\ref{1-system-1}), $\{x_k, k\in \bar{\mathbb{N}}\}$ and $\{y_k, k\in \bar{\mathbb{N}}\}$ are the state processes  
and $\{u_k \in \mathbb{R}^m, k\in \mathbb{N}\}$ is a control process. 
$\{w_k,v_k, k\in \mathbb{N}\}$ are defined on probability space $(\Omega, \mathcal{F}, P)$, 
represent the stochastic distribution for the two state processes, and are assumed to be martingale difference sequences
\begin{eqnarray}\label{1-moment-1}
\begin{array}{l}
\mathbb{E}[w_{k+1}|\mathcal{F}_k]=0,~
\mathbb{E}[(w_{k+1})^2|\mathcal{F}_k]=1,~
\mathbb{E}[v_{k+1}|\mathcal{F}_k]=0,\\
\mathbb{E}[(v_{k+1})^2|\mathcal{F}_k]=1,~
\mathbb{E}[w_{k+1}v_{k+1}|\mathcal{F}_k]=\rho,
\end{array}
\end{eqnarray}
where $\mathcal{F}_k$ is the $\sigma$-algebra generated by $\{\zeta^x,
w_l, l=0, 1,\cdots,k\}$ and $\{\zeta^y,v_l, l=0, 1,\cdots,k\}$. 
The cost functional associated with (\ref{1-system-1}) is
\begin{eqnarray}\label{1-cost-1}
\begin{array}{l}
J(\zeta^x,\zeta^y,u) = \ds\sum_{k=0}^{N-1}\mathbb{E}\Big((x_k-y_k)^T Q_k (x_k-y_k)+u_k^TR_ku_k\\
\quad +\mathbb{E}(x_k-y_k)^T\bar{Q}_k\mathbb{E}(x_k-y_k)+(\mathbb{E}u_k)^T\bar{R}_k\mathbb{E}u_k\Big)\\
\quad +\mathbb{E}\left((x_N-y_N)^T{Q}_N (x_N-y_N)\right)\\
\quad +\mathbb{E}(x_N-y_N)^T\bar{Q}_N\mathbb{E}(x_N-y_N),
\end{array}
\end{eqnarray}
where $Q_k, \bar{Q}_k, k\in\mathbb{\bar{N}}$ and $R_k, \bar{R}_k, k\in\mathbb{N}$ are deterministic symmetric matrices with appropriate dimensions. We introduce the following admissible control set of $u=(u_0,u_1,\cdots,u_{N-1})$
\begin{eqnarray*}
\begin{array}{l}
\mathcal{U}_{ad}\equiv\left\{u~
\left\vert ~u_k  \in \mathbb{R}^m, \mbox{ is }\mathcal{F}_k \mbox{-measurable},~\mathbb{E}|u_k|^2<\infty\right.\right\}.
\end{array}
\end{eqnarray*}
The optimal control problem considered in this paper is then stated as follows:

\noindent
\textbf{Problem (MF-LQ)}.
\emph{For any given square-integrable initial values $\zeta^x$ and $\zeta^y$, find  $u^o\in \mathcal{U}_{ad}$ such that}
\begin{eqnarray} \label{1-problem-1}
J(\zeta^x,\zeta^y,u^o)=\inf_{u\in \mathcal{U}_{ad}}J(\zeta^x,\zeta^y,u).
\end{eqnarray}
\emph{We then call $u^o$ an optimal control for Problem (MF-LQ).}

\section{Preliminaries\label{3}}

In this section, we convert Problem (MF-LQ) to a quadratic optimization problem in Hilbert space. After a series of statements of standard notation and definition, we give necessary and sufficient conditions for the solvability of the Problem (MF-LQ). Firstly, some spaces are introduced as follows: for $k \in \bar{\mathbb{N}}$,
\begin{eqnarray*}
\begin{array}{l}
\mathcal{Z}_k=L_{\mathcal{F}_{k}}^{2}(\mathbb{R}^n)=\left\{\xi:\Omega \mapsto \mathbb{R}^n 
\left\vert\begin{array}{l}
\xi\mbox{ is }\mathcal{F}_{k}\mbox{-measurable,}\\[-2mm]
\mathbb{E}|\xi|^2<\infty
\end{array}\right.\right\},\\[6mm]
\mathcal{Z}[0,k]=\left\{(z_0,...,z_k)
\left\vert\begin{array}{l}
~z_k \in \mathcal{Z}_k, \mbox{ is }\mathcal{F}_{k}\mbox{-measurable,}\\[-2mm]
\sum^k_{l=0}\mathbb{E}|z_l|^2<\infty
\end{array}\right.\right\},
\end{array}
\end{eqnarray*}

and for $l \in \mathbb{N}$,
\begin{eqnarray*}
\begin{array}{l}
\mathcal{U}_l=L_{\mathcal{F}_{l}}^{2}(\mathbb{R}^m)=\left\{\eta:\Omega \mapsto \mathbb{R}^m
\left\vert\begin{array}{l}
\eta\mbox{ is }\mathcal{F}_{l}\mbox{-measurable,}\\[-2mm]
\mathbb{E}|\eta|^2<\infty
\end{array}\right.\right\}.

\end{array}
\end{eqnarray*}
Here, $\mathcal{Z}_k$, $\mathcal{U}_l$ and $\mathcal{Z}[0,k]$ are Hilbert spaces. There are two cases of the domain and range of expectation operator: $\mathbb{E}$ maps $\mathcal{Z}_k$ to $\mathbb{R}^{n}$ or $\mathcal{U}_l$ to $\mathbb{R}^{m}$ \cite{Elliott2013}. Therefore, the notation $\mathbb{E}$ and adjoint operator $\mathbb{E}^*$ may differ from place to place. Let $\mathcal{H}=\mathcal{Z}_k,\mathcal{U}_l$. For illustration, we now use $\mathbb{E}_{\mathcal{H}}$ and $\mathbb{E}^*_{\mathcal{H}}$ to emphasize $\mathcal{H}$. $\mathbb{E}$ and $\mathbb{E}^*$ may appear in the form of $M\mathbb{E}_{\mathcal{H}}$ and $\mathbb{E}^*_{\mathcal{H}}N\mathbb{E}_{\mathcal{H}'}$ where $M$, $N$ are matrices with appropriate dimensions and $\mathcal{H}$, $\mathcal{H}'$ can be different. 
To simplify the expressions in this paper, $\bar{A}\mathbb{E}z$, $\bar{B}\mathbb{E}u$, $\mathbb{E}^*\bar{Q}\mathbb{E}z$, $\mathbb{E}^*\bar{L}\mathbb{E}u$, $\mathbb{E}^*\bar{R}\mathbb{E}u$ are used to denote $\bar{A}\mathbb{E}_{\mathcal{Z}_k}(z)$, $\bar{B}\mathbb{E}_{\mathcal{U}_k}(u)$, $\mathbb{E}^*_{\mathcal{Z}_k}\bar{Q}\mathbb{E}_{\mathcal{Z}_k}(z)$, $\mathbb{E}^*_{\mathcal{Z}_k}\bar{L}\mathbb{E}_{\mathcal{U}_k}(u)$, $\mathbb{E}^*_{\mathcal{U}_k}\bar{R}\mathbb{E}_{\mathcal{U}_k}(u)$, respectively. Here, $z\in\mathcal{Z}_k$, $u\in\mathcal{U}_k$, $\bar{A},\bar{Q}\in\mathbb{R}^{n \times n}$, $\bar{B},\bar{L}\in\mathbb{R}^{n \times m}$, $\bar{R}\in\mathbb{R}^{m \times m}$. We then give the following notion.
\begin{definition}\label{1-definition-1}
(i). Problem (MF-LQ) is said to be finite for $\zeta^x$ and $\zeta^y$ if
\begin{eqnarray*}
\inf_{u\in \mathcal{U}_{ad}}J(\zeta^x,\zeta^y,u)>-\infty.
\end{eqnarray*}
Problem (MF-LQ) is said to be finite if it is finite for any $\zeta^x$ and $\zeta^y$.\\
(ii). Problem is said to be uniquely solvable for $\zeta^x$ and $\zeta^y$ if there exists a unique $u^o \in \mathcal{U}_{ad}$ such that (\ref{1-problem-1}) holds for $\zeta^x$ and $\zeta^y$. Problem (MF-LQ) is said to be uniquely solvable for any $\zeta^x$ and $\zeta^y$.
\end{definition}

We express the system states explicitly in terms of $k$.
Let 
\begin{eqnarray*}
\left\{\begin{array}{l}
\bar{\Phi}(k,l)=\prod\limits^k_{i=l}(A_{i}+\bar{A}_{i}),~k\geq l, \\
\bar{\Phi}(k,l)=I,~k<l,\\
{\Phi}(k,l)=\prod\limits^k_{i=l}(A_{i}+w_{i}C_{i}),~ k\geq l, \\
{\Phi}(k,l)=I, k<l.
\end{array}\right.
\end{eqnarray*}
and
\begin{eqnarray*}
\left\{\begin{array}{l}
\bar{\Xi}(k,l)=\prod\limits^k_{i=l}(F_{i}+\bar{F}_{i}),~k\geq l, \\
\bar{\Xi}(k,l)=I,~k<l,\\
{\Xi}(k,l)=\prod\limits^k_{i=l}(F_{i}+v_{i}G_{i}),~k\geq l, \\
{\Xi}(k,l)=I,~k<l.
\end{array}\right.
\end{eqnarray*}
Define the following operators on $\zeta^x,\zeta^y \in \mathcal{Z}_0$, $u\in \mathcal{U}_{ad}$ for $k \in \bar{\mathbb{N}}$:
\begin{eqnarray*}
\left\{\begin{array}{rl}
\Gamma_k (\zeta^x)=&{\Phi}(k-1,0)\zeta^x, \\ 
\bar{\Gamma}_k(\zeta^x)=&\sum_{l=1}^{k-1}\big(\Phi(k-1,l)(\bar{A}_{l-1}+w_{l-1}\bar{C}_{l-1})\\
&\cdot\bar{\Phi}(l-2,0)\mathbb{E}\zeta^x\big), \\
\Psi_k (\zeta^y)=&{\Xi}(k-1,0)\zeta^y, \\
\bar{\Psi}_k(\zeta^y)=&\sum_{l=1}^{k-1}\big(\Xi(k-1,l)(\bar{F}_{l-1}+v_{l-1}\bar{G}_{l-1})\\
&\cdot\bar{\Xi}(l-2,0)\mathbb{E}\zeta^y\big), \\
L_k(u)=&\sum_{l=1}^{k-1}\Phi(k-1,l)(B_{l-1}+w_{l-1}D_{l-1})u_{l-1}\\
&+(B_{k-1}+w_{k-1}D_{k-1}u_{k-1}), \\
\bar{L}_k(u)=&\sum_{l=1}^{k-1}\Phi(k-1,l+1)(\bar{A}_{l-1}+w_{l-1}\bar{C}_{l-1})\\
&\cdot\sum_{i=1}^{l-1}\bar{\Phi}(l-2, i)(B_{i-1}+\bar{B}_{i-1})\mathbb{E}u_{i-1} \\
&+\sum_{l=1}^{k-1}\Phi(k-1,l)(\bar{B}_{l-1}+w_{l-1}\bar{D}_{l-1})\mathbb{E}u_{l-1}\\
&+(\bar{B}_{k-1}+w_{k-1}\bar{D}_{k-1})\mathbb{E}u_{k-1},
\end{array}\right.
\end{eqnarray*}
%
where
\begin{eqnarray*}
\left\{\begin{array}{l}
\Gamma_k,\bar{\Gamma}_k,\Psi_k,\bar{\Psi}_k:\mathcal{Z}_0\mapsto \mathcal{Z}[0,k], ~k\in\bar{\mathbb{N}},\\
L_k,\bar{L}_k:\mathcal{U}_{ad}\mapsto \mathcal{Z}[0,k],~k\in\bar{\mathbb{N}},
\end{array}\right.
\end{eqnarray*}
are linear and bounded. Then the system states can be expressed as
\begin{eqnarray*}
\left\{\begin{array}{l}
x_k=\Gamma_k (\zeta^x)+\bar{\Gamma}_k(\zeta^x)+L_k(u)+\bar{L}_k(u)\\
y_k=\Psi_k (\zeta^y)+\bar{\Psi}_k (\zeta^y).
\end{array}\right.
\end{eqnarray*}
The cost functional $J(\zeta^x,\zeta^y,u)$ has the following form of usual inner products
\begin{eqnarray*}
\begin{array}{l}
J=\ds\sum_{k=0}^{N-1}\Big(\langle Q_kx_k, x_k\rangle -2\langle Q_kx_k,y_k \rangle + \langle Q_ky_k,y_k \rangle\\
\quad +\langle \bar{Q}\mathbb{E}x_k, \mathbb{E}x_k\rangle-2\langle \bar{Q}_k\mathbb{E}x_k, \mathbb{E}y_k \rangle + \langle \bar{Q}_k\mathbb{E}y_k, \mathbb{E}y_k\rangle\\
\quad +\langle R_ku_k, u_k\rangle +\langle \bar{R}_k\mathbb{E}u_k, \mathbb{E}u_k \rangle\Big) +\langle Q_Nx_N, x_N\rangle\\
\quad -2\langle Q_Nx_N,y_N \rangle + \langle Q_Ny_N,y_N \rangle +\langle \bar{Q}\mathbb{E}x_N, \mathbb{E}x_N\rangle \\
\quad -2\langle \bar{Q}_N\mathbb{E}x_N, \mathbb{E}y_N \rangle + \langle \bar{Q}_N\mathbb{E}y_N, \mathbb{E}y_N\rangle,
\end{array}
\end{eqnarray*}
that is
\begin{eqnarray}\label{1-cost-4}
\begin{array}{l}
J(\zeta^x,\zeta^y,u)= \ds\sum_{k=0}^{N-1}\Big(\big\langle Q_k\big(\Gamma_k (\zeta^x)+\bar{\Gamma}_k(\zeta^x)+L_k(u)\\
\quad +\bar{L}_k(u)\big), \big(\Gamma_k( \zeta^x)+\bar{\Gamma}_k(\zeta^x)+L_k(u)+\bar{L}_k(u)\big)\big\rangle\\
\quad -2\big\langle Q_k\big(\Gamma_k (\zeta^x)+\bar{\Gamma}_k(\zeta^x)+L_k(u)+\bar{L}_k(u)\big),\\
\quad \big(\Psi_k(\zeta^y)+\bar{\Psi}_k(\zeta^y)\big) \big\rangle + \big\langle Q_k\big(\Psi_k(\zeta^y)+\bar{\Psi}_k(\zeta^y)\big),\\
\quad \big(\Psi_k(\zeta^y)+\bar{\Psi}_k (\zeta^y)\big) \big\rangle +\big\langle \bar{Q}\mathbb{E}\big(\Gamma_k \zeta^x+\bar{\Gamma}\zeta^x\\
\quad +L_k(u)+\bar{L}_k(u)\big), \mathbb{E}\big(\Gamma (\zeta^x)+\bar{\Gamma}_N(\zeta^x)+L_k(u)\\
\quad +\bar{L}_k(u)\big)\big\rangle -2\big\langle \bar{Q}_k\mathbb{E}\big(\Gamma_k( \zeta^x)+\bar{\Gamma}_k(\zeta^x)+L_k(u)\\
\quad +\bar{L}_k(u)\big), \mathbb{E}\big(\Psi_k (\zeta^y)+\bar{\Psi}_k (\zeta^y)) \big\rangle +\big\langle \bar{Q}_k\mathbb{E}(\Psi_k(\zeta^y)\\
\quad +\bar{\Psi}_k(\zeta^y)\big), \mathbb{E}\big(\Psi_k (\zeta^y)+\bar{\Psi}_k(\zeta^y )\big)\big\rangle  +\langle R_ku_k, u_k\big\rangle\\
\quad +\big\langle \bar{R}_k\mathbb{E}u_k, \mathbb{E}u_k \big\rangle\Big)+\big\langle Q_N\big(\Gamma_N (\zeta^x)+\bar{\Gamma}_N(\zeta^x)\\
\quad +L_N(u)+\bar{L}_N(u)\big), \big(\Gamma_N( \zeta^x)+\bar{\Gamma}_N(\zeta^x)+L_N(u)\\
\quad +\bar{L}_N(u)\big)\big\rangle -2\big\langle Q_N\big(\Gamma_N (\zeta^x)+\bar{\Gamma}_N(\zeta^x)+L_N(u)\\
\quad +\bar{L}_N(u)\big), \big(\Psi_N(\zeta^y)+\bar{\Psi}_N(\zeta^y)) \big\rangle + \big\langle Q_N(\Psi_N(\zeta^y)\\ 
\quad +\bar{\Psi}_N(\zeta^y)\big), \big(\Psi_N(\zeta^y)+\bar{\Psi}_N (\zeta^y)\big) \big\rangle +\big\langle \bar{Q}_N\mathbb{E}\big(\Gamma_N \zeta^x\\
\quad +\bar{\Gamma}_N\zeta^x+L_N(u)+\bar{L}_N(u)\big), \mathbb{E}\big(\Gamma_N (\zeta^x)+\bar{\Gamma}_N(\zeta^x)\\
\quad +L_N(u)+\bar{L}_N(u)\big)\big\rangle -2\big\langle \bar{Q}_N\mathbb{E}\big(\Gamma_N( \zeta^x)+\bar{\Gamma}_N(\zeta^x)\\
\quad +L_N(u)+\bar{L}_N(u)\big), \mathbb{E}\big(\Psi_N (\zeta^y)+\bar{\Psi}_N (\zeta^y)\big) \big\rangle\\
\quad +\big\langle \bar{Q}_N\mathbb{E}\big(\Psi_N(\zeta^y)+\bar{\Psi}_N(\zeta^y)\big), \mathbb{E}\big(\Psi_N (\zeta^y)+\bar{\Psi}_N(\zeta^y )\big)\big\rangle\\
=  2\langle \Theta_1u,\zeta^x\rangle +\langle \Theta_2u, u\rangle + 2\langle \Theta_3u, \zeta^y\rangle + \langle \Lambda_1 \zeta^x, \zeta^x\rangle\\
\quad + 2\langle \Lambda_2\zeta^x, \zeta^y\rangle + \langle \Lambda_3\zeta^y, \zeta^y \rangle.
\end{array}
\end{eqnarray}
Recall that $\langle Q_kx_k, y_k\rangle$ denotes $\mathbb{E}(y_k^TQ_kx_k)$ with similar meanings for related notation. Here,
\begin{eqnarray*}
\begin{array}{rl}
\Theta_1 =& \ds\sum^{N}_{k=0}\Big((\Gamma_k+\bar{\Gamma}_k)^*Q_k(L_k+\bar{L}_k)\\
&+(\Gamma_k+\bar{\Gamma}_k)^*\mathbb{E}^*\bar{Q}_k\mathbb{E}(L_k+\bar{L}_k)\Big),\\
\Theta_2 =& \ds\sum_{k=0}^{N-1}\Big(R_k+\mathbb{E}^*\bar{R}_k\mathbb{E}+(L_k+\bar{L}_k)^*Q_k(L_k+\bar{L}_k)\\
&+(L_k+\bar{L}_k)^*\mathbb{E}^*\bar{Q}_k\mathbb{E}(L_k+\bar{L}_k) \Big)\\
&+({L}_N+{\bar{L}_N})^*Q_N({L}_N+{\bar{L}_N})\\
&+({L}_N+{\bar{L}_N})^*\mathbb{E}^*\bar{Q}_N\mathbb{E}({L}_N+{\bar{L}_N}),\\
\end{array}
\end{eqnarray*}
\begin{eqnarray*}
\begin{array}{rl}
\Theta_3 =& \ds\sum_{k=0}^{N}\Big((\Psi_k+\bar{\Psi}_k)^*Q_k(L_k+\bar{L}_k)\\
&+(\Psi_k+\bar{\Psi}_k)^*\mathbb{E}^*\bar{Q}_k\mathbb{E}(L_k+\bar{L}_k)\Big),\\
\Lambda_1 =& \ds\sum_{k=0}^{N}\Big((\Gamma_k+\bar{\Gamma}_k)^*Q_k(\Gamma_k+\bar{\Gamma}_k)\\
&+(\Gamma_k+\bar{\Gamma}_k)^*\mathbb{E}^*\bar{Q}_k\mathbb{E}(\Gamma_k+\bar{\Gamma}_k)\Big),\\
\Lambda_2 =& -\ds\sum_{k=0}^{N}\Big((\Psi_k+\bar{\Psi}_k)^*Q_k(\Gamma_k+\bar{\Gamma}_k)\\
&+(\Psi_k+\bar{\Psi}_k)^*\mathbb{E}^*\bar{Q}_k\mathbb{E}(\Gamma_k+\bar{\Gamma}_k)\Big),\\ 
\Lambda_3 =& \ds\sum_{k=0}^{N}\Big((\Psi_k+\bar{\Psi}_k)^*Q_k(\Psi_k+\bar{\Psi}_k)\\
&+(\Psi_k+\bar{\Psi}_k)^*\mathbb{E}^*\bar{Q}_k\mathbb{E}(\Psi_k+\bar{\Psi}_k)\Big).\\
\end{array}
\end{eqnarray*}
Note that in this paper we use numerator layout for matrices calculus, i.e., for any matrix $Y$, $\frac{\partial}{\partial Y}Tr(AYB)=BA$ if $AYB$ is meaningful. We then have the following result.
\begin{proposition}\label{1-proposition-1}
(i). If $J(\zeta^x, \zeta^y, u)$ has a minimum, then
$$\Theta_2\geq 0.$$ 
(ii). Problem (MF-LQ) is (uniquely) solvable if and only if $\Theta_2\geq 0$ and there  exists a (unique) ${u}$ such that
\begin{eqnarray*}
u^T\Theta_2+x^T\Theta_1+ y^T\Theta_3=0.
\end{eqnarray*}
(iii). If $\Theta_2>0$, then for any $\zeta^x$ and $\zeta^y$, $J(\zeta^x, \zeta^y, u)$ admits a pathwise unique minimizer $u^o$ given by
\begin{eqnarray}\label{1-control-1}
u^o_k=-\big(\Theta_2^{-1}(\Theta_1^*\zeta^x + \Theta_3^*\zeta^y) \big)(k), ~~k\in\mathbb{N}.
\end{eqnarray}
In addition, if
\begin{eqnarray}\label{1-condition-1}
Q_k, Q_k+\bar{Q}_k\geq 0,~k\in\bar{\mathbb{N}},~R_k, R_k+\bar{R}_k>0,~k\in\mathbb{N},
\end{eqnarray}
then $\Theta_2>0$.
\end{proposition}
\emph{Proof}. The proofs of (i), (ii) and the first part of (iii) are well known and therefore omitted here \cite{Moore1999,Yong2013}. We now prove the second part of (iii). From (\ref{1-condition-1}), for $k\in\bar{\mathbb{N}}$, we have
\begin{eqnarray*}
\begin{array}{l}
(L_k+\bar{L}_k)^*Q_k(L_k+\bar{L}_k)\\
+(L_k+\bar{L}_k)^*\mathbb{E}^*\bar{Q}_k\mathbb{E}(L_k+\bar{L}_k) \geq 0.
\end{array}
\end{eqnarray*}
Also,
\begin{eqnarray*}
\begin{array}{l}
\langle R_k u_k,u_k \rangle +\langle \bar{R}_k\mathbb{E}u_k,\mathbb{E}u_k \rangle\\
= \mathbb{E} \big[u_k^TR_ku_k+(\mathbb{E}u_k)^T\bar{R}_k\mathbb{E}u_k \big]\\
= \mathbb{E}\big[(u_k-\mathbb{E}u_k)^TR_k(u_k-\mathbb{E}u_k)\big]\\
\quad +(\mathbb{E}u_k )^T(R_k+\bar{R}_k)\mathbb{E}u_k > 0,~k\in\mathbb{N}
\end{array}
\end{eqnarray*}
for any non-zero $u\in\mathcal{U}_{ad}$, which implies $\Theta_2>0$. This completes the proof.
\hfill $\square$

\section{Closed-loop Optimal Control via Riccati Equations\label{4}}

In this section, we first find the optimal control within the class of linear state feedback controls by using matrix minimum principle. Secondly, several sequences of bounded linear operators are presented and problem (MF-LQ) is reformulated as an operator stochastic linear-quadratic optimal control problem. We then find the optimal control via Riccati equations. Finally, we show that the two above approaches employed for solving Problem (MF-LQ) coincide by completing the square.
\subsection{Linear Feedback Control\label{4.1}}
The linear feedback controls a linear functional of the system states, which gives the control based on the current system states. 
Suppose that a control having the following form is in used:
\begin{eqnarray}\label{1-control-2}
u_k=L^x_kx_k+\bar{L}^x_k\mathbb{E}x_k+L^y_ky_k+\bar{L}^y_k\mathbb{E}y_k,~~k\in\mathbb{N},
\end{eqnarray} 
where $L^x_k, \bar{L}^x_k, L^y_k, \bar{L}^y_k \in \mathbb{R}^{m\times n}$. Under (\ref{1-control-2}), the closed-loop system (\ref{1-system-1}) becomes
\begin{eqnarray}\label{1-system-2}
\left\{\begin{array}{l}
x_{k+1}=  A_kx_k+\bar{A}_k \mathbb{E}x_k+B_k(L^x_kx_k+\bar{L}^x_k\mathbb{E}x_k+L^y_ky_k\\
\quad +\bar{L}^y_k\mathbb{E}y_k)+ \bar{B}_k\big[(L^x_k +\bar{L}^x_k)\mathbb{E}x_k+(L^y_k+\bar{L}^y_k)\mathbb{E}y_k\big]\\
\quad +\big\{C_kx_k+\bar{C}_k \mathbb{E}x_k+D_k(L^x_kx_k+\bar{L}^x_k\mathbb{E}x_k+L^y_ky_k\\
\quad +\bar{L}^y_k\mathbb{E}y_k)+ \bar{D}_k\big[(L^x_k +\bar{L}^x_k)\mathbb{E}x_k+(L^y_k+\bar{L}^y_k)\mathbb{E}y_k\big]\big\}w_k,\\
y_{k+1}= (F_ky_k+\bar{F}_k\mathbb{E}y_k)+ (G_ky_k+\bar{G}_k\mathbb{E}y_k)v_k, \\
x_0=  \zeta^x, y_0=\zeta^y,
\end{array}\right.
\end{eqnarray}
and the cost functional (\ref{1-cost-1}) may be represented as
\begin{eqnarray}\label{1-cost-2}
\begin{array}{l}
J(\zeta^x,\zeta^y,u)\\
= \ds\sum_{k=0}^{N-1}\mathbb{E}\Big((x_k-y_k)^TQ_k(x_k-y_k)+\mathbb{E}(x_k-y_k)^T\bar{Q}_k\mathbb{E}\\
\quad \cdot(x_k-y_k)+(L^x_kx_k+\bar{L}^x_k\mathbb{E}x_k+L^y_ky_k+\bar{L}^y_k\mathbb{E}y_k)^T R_k\\
\quad \cdot(L^x_kx_k+\bar{L}^x_k\mathbb{E}x_k+L^y_ky_k+\bar{L}^y_k\mathbb{E}y_k)+\big((L^x_k+\bar{L}^x_k)\mathbb{E}x_k\\
\quad +(L^y_k+\bar{L}^y_k)\mathbb{E}y_k\big)^T\bar{R}_k\big((L^x_k+\bar{L}^x_k)\mathbb{E}x_k+(L^y_k+\bar{L}^y_k)\mathbb{E}y_k\big)\Big)\\
\quad +\mathbb{E}\big((x_{N}-y_N)^TQ_N(x_{N}-y_N)\big)\\
\quad +\mathbb{E}(x_{N}-y_N)^T\bar{Q}_N\mathbb{E}(x_{N}-y_N)\\
= \ds\sum_{k=0}^{N-1}\big\{
Tr\big[Q_k\big(\mathbb{E}(x_kx_k^T)-\mathbb{E}(x_ky_k^T)-\mathbb{E}(y_kx_k^T)+\mathbb{E}(y_ky_k^T)\big)\big]\\
\quad +Tr\big[ \bar{Q}_k \big(\mathbb{E}x_k\mathbb{E}x_k^T-\mathbb{E}x_k\mathbb{E}y_k^T-\mathbb{E}y_k\mathbb{E}x_k^T+\mathbb{E}y_k\mathbb{E}y_k^T\big)\big]\\
\quad  +Tr\big[(L^x_k)^TR_kL^x_k \mathbb{E}(x_kx_k^T)+\big((L^x_k)^TR_k\bar{L}^x_k+(\bar{L}^x_k)^TR_kL^x_k\\
\quad +(\bar{L}^x_k)^TR_k\bar{L}^x_k\big)\mathbb{E}x_k\mathbb{E}x_k^T \big]+Tr\big[(L^x_k)^TR_kL^y_k\mathbb{E}(y_kx_k^T) \\
\quad +\big((L^x_k)^TR_k\bar{L}^y_k+(\bar{L}^x_k)^TR_kL^y_k+(\bar{L}^x_k)^TR_k\bar{L}^y_k\big)\mathbb{E}y_k\mathbb{E}x_k^T \big]\\
\quad +Tr\big[(L^y_k)^TR_kL^x_k \mathbb{E}(x_ky_k^T) +\big((L^y_k)^TR_k\bar{L}^x_k+(\bar{L}^y_k)^TR_kL^x_k\\
\quad +(\bar{L}^y_k)^TR_k\bar{L}^x_k\big)\mathbb{E}x_k\mathbb{E}y_k^T \big]+Tr\big[(L^y_k)^TR_kL^y_k\mathbb{E}(y_ky_k^T)\\
\quad +\big((L^y_k)^TR_k\bar{L}^y_k+(\bar{L}^y_k)^TR_kL^y_k+(\bar{L}^y_k)^TR_k\bar{L}^y_k\big)\mathbb{E}y_k\mathbb{E}y_k^T \big]\\
\quad +Tr\big[(L^x_k+\bar{L}^x_k)^T\bar{R}_k(L^x_k+\bar{L}^x_k\big)\mathbb{E}x_k\mathbb{E}x_k^T+ (L^x_k+\bar{L}^x_k)^T\\
\quad \cdot\bar{R}_k(L^y_k+\bar{L}^y_k\big)\mathbb{E}y_k\mathbb{E}x_k^T \big]+Tr\big[(L^y_k+\bar{L}^y_k)^T\bar{R}_k(L^x_k+\bar{L}^x_k\big)\\
\quad \cdot\mathbb{E}x_k\mathbb{E}y_k^T+ (L^y_k+\bar{L}^y_k)^T\bar{R}_k(L^y_k+\bar{L}^y_k\big)\mathbb{E}y_k\mathbb{E}y_k^T \big]\big\}\\
\quad +Tr\big[Q_N\big(\mathbb{E}(x_Nx_N^T)-\mathbb{E}(x_Ny_N^T)-\mathbb{E}(y_Nx_N^T)+\mathbb{E}(y_Ny_N^T)\big)\big]\\
\quad +Tr\big[ \bar{Q}_N \big(\mathbb{E}x_N\mathbb{E}x_N^T-\mathbb{E}x_N\mathbb{E}y_N^T-\mathbb{E}y_N\mathbb{E}x_N^T+\mathbb{E}y_N\mathbb{E}y_N^T\big)\big].
\end{array}
\end{eqnarray}
From the form (\ref{1-control-2}) of the control, we may view $\{(L^x_k, \bar{L}^x_k, L^y_k, \bar{L}^y_k), k\in \mathbb{N}\}$ as the new control input. Also (\ref{1-cost-2}) reminds us that $\mathbb{E}\left(x_kx_k^T\right)$, $\mathbb{E}x_k(\mathbb{E}x_k)^T$, $\mathbb{E}\left(x_ky_k^T\right)$, $\mathbb{E}x_k(\mathbb{E}y_k)^T$, $\mathbb{E}\left(y_ky_k^T\right)$ and $\mathbb{E}y_k(\mathbb{E}y_k)^T$ may be considered as the new system states. Write $X_k=\mathbb{E}\left(x_kx_k^T\right)$, $\bar{X}_k=\mathbb{E}x_k(\mathbb{E}x_k)^T$, $XY_k=\mathbb{E}\left(x_ky_k^T\right)$, $\bar{XY}_k=\mathbb{E}x_k(\mathbb{E}y_k)^T$, $Y_k=\mathbb{E}\left(y_ky_k^T\right)$ and $\bar{Y}_k=\mathbb{E}y_k(\mathbb{E}y_k)^T$. Then by (\ref{1-system-2}), we have
\begin{eqnarray}\label{1-X-1}
\left\{
\begin{array}{l}
X_{k+1}=(A_k+B_kL^x_k)X_k(A_k+B_kL^x_k)^T+\big[\bar{A}_k+B_k\bar{L}^x_k\\
\quad +\bar{B}_k(L^x_k+\bar{L}^x_k)\big]\bar{X}_k\big[\bar{A}_k+B_k\bar{L}^x_k+\bar{B}_k(L^x_k+\bar{L}^x_k)\big]^T\\
\quad +(A_k+B_kL^x_k)\bar{X}_k\big[\bar{A}_k+B_k\bar{L}^x_k+\bar{B}_k(L^x_k+\bar{L}^x_k)\big]^T\\
\quad +\big[\bar{A}_k+B_k\bar{L}^x_k+\bar{B}_k(L^x_k+\bar{L}^x_k)\big]\bar{X}_k(A_k+B_kL^x_k)^T\\
\quad +(A_k+B_kL^x_k)XY_k(B_kL^y_k)^T+\big[\bar{A}_k+B_k\bar{L}^x_k\\
\quad +\bar{B}_k(L^x_k+\bar{L}^x_k)\big]\bar{XY}_k\big[B_k\bar{L}^y_k+\bar{B}_k(L^y_k+\bar{L}^y_k)\big]^T\\
\quad +(A_k+B_kL^x_k)\bar{XY}_k\big[B_k\bar{L}^y_k+\bar{B}_k(L^y_k+\bar{L}^y_k)\big]^T\\
\quad +\big[\bar{A}_k+B_k\bar{L}^x_k+\bar{B}_k(L^x_k+\bar{L}^x_k)\big]\bar{XY}_k(B_kL^y_k)^T\\
\end{array}\right.
\end{eqnarray}
\begin{eqnarray}\label{1-barX-1}
\left\{\begin{array}{l}
\bar{X}_{k+1}=\big[A_k+\bar{A}_k+(B_k+\bar{B}_k)(L^x_k+\bar{L}^x_k)\big]\bar{X}_k\\
\quad \cdot\big[A_k+\bar{A}_k+(B_k+\bar{B}_k)(L^x_k+\bar{L}^x_k)\big]^T+\big[A_k+\bar{A}_k\\
\quad +(B_k+\bar{B}_k)(L^x_k+\bar{L}^x_k)\big] \bar{XY}_k \big[(B_k+\bar{B}_k)(L^y_k+\bar{L}^y_k)\big]^T\\
\quad +\big[(B_k+\bar{B}_k)(L^y_k+\bar{L}^y_k)\big] (\bar{XY}_k)^T \big[A_k+\bar{A}_k\\
\quad +(B_k+\bar{B}_k)(L^x_k+\bar{L}^x_k)\big]^T+\big[(B_k+\bar{B}_k)(L^y_k+\bar{L}^y_k)\big]\\
\quad \cdot\bar{Y}_k\big[(B_k+\bar{B}_k)(L^y_k+\bar{L}^y_k)\big]^T\\
\quad\equiv\bar{\mathcal{X}}_k(L^x_k,\bar{L}^x_k,L^y_k,\bar{L}^y_k),\\
\bar{X}_0=\mathbb{E}\zeta^x \big(\mathbb{E}\zeta^x\big)^T,
\end{array}\right.
\end{eqnarray}
\begin{eqnarray}\label{1-XY-1}
\left\{\begin{array}{l}
XY_{k+1}=
 (A_k+B_kL^x_k)XY_k F_k^T+\big[\bar{A}_k+B_k\bar{L}^x_k+\bar{B}_k\\
\quad \cdot(L^x_k+\bar{L}^x_k)\big]\bar{XY}_k\bar{F}_k^T+(A_k+B_kL^x_k)\bar{XY}_k\bar{F}_k^T\\
\quad +\big[\bar{A}_k+B_k\bar{L}^x_k+\bar{B}_k(L^x_k+\bar{L}^x_k)\big]\bar{XY}_kF_k^T+(B_kL^y_k)\\
\quad \cdot Y_kF_k^T+(B_kL^y_k)\bar{Y}_k\bar{F}_k^T +\big[B_k\bar{L}^y_k+\bar{B}_k(L^y_k+\bar{L}^y_k)\big]\\
\quad \cdot\bar{Y}_kF_k^T+\big[B_k\bar{L}^y_k+\bar{B}_k(L^y_k+\bar{L}^y_k)\big]\bar{Y}_k\bar{F}_k^T+\rho\big\{(C_k\\
\quad +D_kL^x_k)XY_k G_k^T+\big[\bar{C}_k+D_k\bar{L}^x_k+\bar{D}_k(L^x_k+\bar{L}^x_k)\big]\\
\quad \cdot\bar{XY}_k\bar{G}_k^T+\big[\bar{C}_k+D_k\bar{L}^x_k+\bar{D}_k(L^x_k+\bar{L}^x_k)\big]\bar{XY}_kG_k^T\\
\quad +(C_k+D_kL^x_k)\bar{XY}_k\bar{G}_k^T+(D_kL^y_k)Y_kG_k^T\\
\quad +(D_kL^y_k)\bar{Y}_k\bar{G}_k^T +\big[D_k\bar{L}^y_k+\bar{D}_k(L^y_k+\bar{L}^y_k)\big]\bar{Y}_kG_k^T\\
\quad +\big[D_k\bar{L}^y_k+\bar{D}_k(L^y_k+\bar{L}^y_k)\big]\bar{Y}_k\bar{G}_k^T\big\}\\
\quad\equiv\mathcal{XY}_k(L^x_k,\bar{L}^x_k,L^y_k,\bar{L}^y_k),\\
XY_0=\mathbb{E}\big[\zeta^x(\zeta^y)^T\big],
\end{array}\right.
\end{eqnarray}
\begin{eqnarray}\label{1-barXY-1}
\left\{\begin{array}{l}
\bar{XY}_{k+1}=\big[(A_k+\bar{A}_k)+(B_k+\bar{B}_k)(L^x_k+\bar{L}^x_k)\big]\bar{XY}_k\\
\quad \cdot(F_k+\bar{F}_k)^T+\big[(B_k+\bar{B}_k)(L^y_k+\bar{L}^y_k)\big]\bar{Y}_k(F_k+\bar{F}_k)^T\\
\quad \equiv\bar{\mathcal{XY}}_k(L^x_k,\bar{L}^x_k,L^y_k,\bar{L}^y_k), \\
\bar{XY}_0=\mathbb{E}\zeta^x \big(\mathbb{E}\zeta^y\big)^T,
\end{array}\right.
\end{eqnarray}
\begin{eqnarray}\label{1-Y-1}
\left\{\begin{array}{l}
Y_{k+1}=F_kY_kF_k^T+F_k\bar{Y}_k\bar{F}_k^T+\bar{F}_k\bar{Y}_kF_k^T+\bar{F}_k\bar{Y}_k\bar{F}_k^T\\
\quad +G_kY_kG_k^T+G_k\bar{Y}_k\bar{G}_k^T+\bar{G}_k\bar{Y}_kG_k^T+\bar{G}_k\bar{Y}_k\bar{G}_k^T\\
\quad \equiv\bar{\mathcal{Y}}_k, \\
Y_0=\mathbb{E}[\zeta^y(\zeta^y)^T],
\end{array}\right.
\end{eqnarray}
\begin{eqnarray}\label{1-barY-1}
\left\{\begin{array}{l}
\bar{Y}_{k+1}=(F_k+\bar{F}_k)\bar{Y}_k(F_k+\bar{F}_k)^T\\
\equiv\bar{\mathcal{Y}}_k,\\
\bar{Y}_0=\mathbb{E}\zeta^y \big(\mathbb{E}\zeta^y\big)^T.
\end{array}\right.
\end{eqnarray}
In terms of ${X}$, $\bar{X}$, $XY$, $\bar{XY}$, $Y$ and $\bar{Y}$, $J(\zeta^x,\zeta^y,u)$ with $u$ defined in (\ref{1-control-2}) may be represented as
\begin{eqnarray}\label{1-cost-3}
\begin{array}{l}
J(\zeta^x,\zeta^y,u)\\
= \ds\sum_{k=0}^{N-1}\big\{Tr\big[Q_k\big(X_k-XY_k-(XY_k)^T+Y_k\big)\big]\\
\quad +Tr\big[ \bar{Q}_k \big(\bar{X}_k-\bar{XY}_k-(\bar{XY}_k)^T+\bar{Y}_k\big)\big]\\
\quad +Tr\big[(L^x_k)^TR_kL^x_k X_k +\big((L^x_k)^TR_k\bar{L}^x_k+(\bar{L}^x_k)^TR_kL^x_k\\
\quad +(\bar{L}^x_k)^TR_k\bar{L}^x_k\big)\bar{X}_k \big]+Tr\big[(L^x_k)^TR_kL^y_k (XY_k)^T\\
\quad +\big((L^x_k)^TR_k\bar{L}^y_k+(\bar{L}^x_k)^TR_kL^y_k+(\bar{L}^x_k)^TR_k\bar{L}^y_k\big)(\bar{XY}_k)^T \big]\\
\quad +Tr\big[(L^y_k)^TR_kL^x_k XY_k +\big((L^y_k)^TR_k\bar{L}^x_k+(\bar{L}^y_k)^TR_kL^x_k\\
\quad +(\bar{L}^y_k)^TR_k\bar{L}^x_k\big)\bar{XY}_k\big]+Tr\big[(L^y_k)^TR_kL^y_k Y_k\\
\quad +\big((L^y_k)^TR_k\bar{L}^y_k+(\bar{L}^y_k)^TR_kL^y_k+(\bar{L}^y_k)^TR_k\bar{L}^y_k\big)\bar{Y}_k \big]\\
\quad +Tr\big[(L^x_k+\bar{L}^x_k)^T\bar{R}_k(L^x_k+\bar{L}^x_k\big)\bar{X}_k+ (L^x_k+\bar{L}^x_k)^T\\
\quad \cdot\bar{R}_k(L^y_k+\bar{L}^y_k\big)(\bar{XY}_k)^T \big]+Tr\big[(L^y_k+\bar{L}^y_k)^T\bar{R}_k\\
\quad \cdot(L^x_k+\bar{L}^x_k\big)\bar{XY}_k+ (L^y_k+\bar{L}^y_k)^T\bar{R}_k(L^y_k+\bar{L}^y_k\big)\bar{Y}_k \big]\big\}\\
\quad +Tr\big[Q_N\big(X_N-XY_N-(XY_N)^T+Y_N\big)\big]\\
\quad +Tr\big[ \bar{Q}_N \big(\bar{X}_N-\bar{XY}_N-(\bar{XY}_N)^T+\bar{Y}_N\big)\big]\\
\equiv\mathcal{J}(X_0,\bar{X}_0,XY_0,\bar{XY}_0,Y_0,\bar{Y}_0, \mathcal{L}),
\end{array}
\end{eqnarray}
where $\mathcal{L}\equiv \{(L^x_k, \bar{L}^x_k, L^y_k, \bar{L}^y_k), k\in \mathbb{N}\}$. Therefore, Problem (MF-LQ) is equivalent to the following problem:
\begin{eqnarray}\label{1-problem-2}
\left\{\begin{array}{l}
\ds\min_{L^x_k, \bar{L}^x_k, L^y_k, \bar{L}^y_k\in \mathbb{R}^{m\times n}, k\in \mathbb{N}}\mathcal{J}(X_0,\bar{X}_0,XY_0,\bar{XY}_0,Y_0,\bar{Y}_0, \mathcal{L}),\\
\mbox{subject to } (\ref{1-X-1})\mbox{--}(\ref{1-barY-1}).
\end{array}\right.
\end{eqnarray}
Clearly, this is a matrix dynamical optimization problem. A natural
way to deal with this class of problems is by the matrix minimum
principle \cite{Athans1967}. 
Following the framework above, we can obtain the optimal control of form (\ref{1-control-2}).
Define the optimal feedback gains $L^o_k=( L_k^{xo}, L_k^{yo}, \bar{L}_k^{xo}, \bar{L}_k^{yo} ),~k\in\mathbb{N}$. We first introduce some notation and then present the result.

Let
\begin{eqnarray}\label{1-WH-1}
\left\{\begin{array}{l}
\bar{W}_k^{(1)}=R_k+B_k^TP^x_{k+1}B_k+D_k^TP^x_{k+1}D_k,\\
\bar{W}_k^{(2)}=R_k+\bar{R}_k+(B_k+\bar{B}_k)^T(P^x_{k+1}+\bar{P}^x_{k+1})\\
\quad \cdot(B_k+\bar{B}_k)+(D_k+\bar{D}_k)^TP^x_{k+1}(D_k+\bar{D}_k),\\
\bar{H}_k^{(1)}=A_k^TP^x_{k+1}B_k+C_k^TP^x_{k+1}D_k, \\
\bar{H}_k^{(2)}=(A_k+\bar{A}_k)^T(P^x_{k+1}+\bar{P}^x_{k+1})(B_k+\bar{B}_k)\\
\quad +(C_k+\bar{C}_k)^TP^x_{k+1}(D_k+\bar{D}_k),\\
\bar{H}_k^{(3)}=F_k^T P^{xy}_{k+1}B_k + \rho G_k^T P^{xy}_{k+1}D_k,\\
\bar{H}_k^{(4)}=(F_k+\bar{F}_k)^T (P^{xy}_{k+1}+\bar{P}^{xy}_{k+1})(B_k+\bar{B}_k)\\
\quad +\rho(G_k+\bar{G}_k)^T P^{xy}_{k+1}(D_k+\bar{D}_k),\\
\end{array}\right.
\end{eqnarray}
with
\begin{eqnarray}\label{1-Px-1}
\left\{\begin{array}{l}
P^x_k =Q_k+(L_k^{xo})^TR_kL^{xo}_k+(A_k+B_kL^{xo}_k)^TP^x_{k+1}(A_k\\
\quad +B_kL^{xo}_k)+(C_k+D_kL^{xo}_k)^TP^x_{k+1}(C_k+D_kL^{xo}_k),\\
P^x_N = Q_N,
\end{array}
\right.
\end{eqnarray}
\begin{eqnarray}\label{1-barPx-1}
\left\{\begin{array}{l}
\bar{P}^x_k = \bar{Q}_k+(L_k^{xo})^TR_k\bar{L}^{xo}_k+(\bar{L}_k^{xo})^TR_kL^{xo}_k\\
\quad +(\bar{L}_k^{xo})^TR_k\bar{L}^{xo}_k+(L^{xo}_k+\bar{L}^{xo}_k)^T\bar{R}_k(L^{xo}_k+\bar{L}^{xo}_k)\\
\quad +(A_k+B_kL^{xo}_k)^TP^x_{k+1}\big[\bar{A}_k+B_k\bar{L}^{xo}_k+\bar{B}_k(L^{xo}_k+\bar{L}^{xo}_k)\big]\\
\quad +\big[\bar{A}_k+B_k\bar{L}^{xo}_k+\bar{B}_k(L^{xo}_k+\bar{L}^{xo}_k)\big]^TP^x_{k+1}(A_k+B_kL^{xo}_k)\\
\quad +\big[\bar{A}_k+B_k\bar{L}^{xo}_k+\bar{B}_k(L^{xo}_k+\bar{L}^{xo}_k)\big]^TP^x_{k+1}\big[\bar{A}_k+B_k\bar{L}^{xo}_k\\
\quad +\bar{B}_k(L^{xo}_k+\bar{L}^{xo}_k)\big]+\big[A_k+\bar{A}_k+(B_k+\bar{B}_k)(L^{xo}_k\\
\quad +\bar{L}^{xo}_k)\big]^T\bar{P}^x_{k+1} \big[A_k+\bar{A}_k+(B_k+\bar{B}_k)(L^{xo}_k+\bar{L}^{xo}_k)\big]\\
\quad +(C_k+D_kL^{xo}_k)^TP^x_{k+1}\big[\bar{C}_k+D_k\bar{L}^{xo}_k+\bar{D}_k(L^{xo}_k+\bar{L}^{xo}_k)\big]\\
\quad +\big[\bar{C}_k+D_k\bar{L}^{xo}_k+\bar{D}_k(L^{xo}_k+\bar{L}^{xo}_k)\big]^TP^x_{k+1}(C_k+D_kL^{xo}_k)\\
\quad +\big[\bar{C}_k+D_k\bar{L}^{xo}_k+\bar{D}_k(L^{xo}_k+\bar{L}^{xo}_k)\big]^TP^x_{k+1}\big[\bar{C}_k+D_k\bar{L}^{xo}_k\\
\quad +\bar{D}_k(L^{xo}_k+\bar{L}^{xo}_k)\big],\\
\bar{P}^x_N=\bar{Q}_N,
\end{array}\right.
\end{eqnarray}
\begin{eqnarray}\label{1-Pxy-1}
\left\{\begin{array}{l}
P^{xy}_k= -Q_k+(L^{yo}_k)^TR_kL^{xo}_k+(B_kL^{yo}_k)^T P^x_{k+1}\\
\quad \cdot(A_k+B_kL^{xo}_k)+(D_kL^{yo}_k)^T P^x_{k+1}(C_k+D_kL^{xo}_k)\\
\quad +F_k^T P^{xy}_{k+1}(A_k+B_kL^{xo}_k)+\rho G_k^T P^{xy}_{k+1} (C_k+D_kL^{xo}_k),\\
P^{xy}_N=-Q_N,
\end{array}\right.
\end{eqnarray}
\begin{eqnarray}\label{1-barPxy-1}
\left\{\begin{array}{l}
\bar{P}^{xy}_k
= -\bar{Q}_k+(L^{yo}_k)^TR_k \bar{L}^{xo}_k+(\bar{L}^{yo}_k)^T R_k L^{xo}_k\\
\quad +(\bar{L}^{yo}_k)^TR_k\bar{L}^{xo}_k+(L^{yo}_k+\bar{L}^{yo}_k)^T\bar{R}_k(L^{xo}_k+\bar{L}^{xo}_k)\\
\quad +\big[B_k\bar{L}^{yo}_k+\bar{B}_k(L^{yo}_k+\bar{L}^{yo}_k)\big]^TP^x_{k+1}(A_k+B_kL^{xo}_k)\\
\quad +(B_kL^{yo}_k)^T P^x_{k+1}\big[\bar{A}_k+B_k\bar{L}^{xo}_k+\bar{B}_k(L^{xo}_k+\bar{L}^{xo}_k)\big]\\
\quad +\big[B_k\bar{L}^{yo}_k+\bar{B}_k(L^{yo}_k+\bar{L}^{yo}_k)\big]^T P^x_{k+1}\big[\bar{A}_k+B_k\bar{L}^{xo}_k\\
\quad +\bar{B}_k(L^{xo}_k+\bar{L}^{xo}_k)\big]+\big[(B_k+\bar{B}_k)(L^{yo}_k+\bar{L}^{yo}_k)\big]^T\\
\quad \cdot\bar{P}^x_{k+1}\big[(A_k+\bar{A}_k)+(B_k+\bar{B}_k)(L^{xo}_k+\bar{L}^{xo}_k)\big]\\
\quad +\big[D_k\bar{L}^{yo}_k+\bar{D}_k(L^{yo}_k+\bar{L}^{yo}_k)\big]^T P^x_{k+1}(C_k+D_kL^{xo}_k)\\
\quad +(D_kL^{yo}_k)^T P^x_{k+1}\big[\bar{C}_k+D_k\bar{L}^{xo}_k+\bar{D}_k(L^{xo}_k+\bar{L}^{xo}_k)\big]\\
\quad +\big[D_k\bar{L}^{yo}_k+\bar{D}_k(L^{yo}_k+\bar{L}^{yo}_k)\big]^T P^x_{k+1}\big[\bar{C}_k+D_k\bar{L}^{xo}_k\\
\quad +\bar{D}_k(L^{xo}_k+\bar{L}^{xo}_k)\big]+F_k^T P^{xy}_{k+1} \big[\bar{A}_k+B_k\bar{L}^{xo}_k+\bar{B}_k(L^{xo}_k\\
\quad +\bar{L}^{xo}_k)\big]+\bar{F}_k^TP^{xy}_{k+1}\big[A_k+\bar{A}_k+(B_k+\bar{B}_k)(L^{xo}_k+\bar{L}^{xo}_k)\big]\\
\quad +(F_k+\bar{F}_k)^T\bar{P}^{xy}_{k+1}\big[A_k+\bar{A}+(B_k+\bar{B}_k)(L^{xo}_k+\bar{L}^{xo}_k)\big]\\
\quad +\rho G_k^T P^{xy}_{k+1} \big[\bar{C}_k+D_k\bar{L}^{xo}_k+\bar{D}_k(L^{xo}_k+\bar{L}^{xo}_k)\big]\\
\quad +\rho\bar{G}_k^TP^{xy}_{k+1}\big[C_k+\bar{C}_k+(D_k+\bar{D}_k)(L^{xo}_k+\bar{L}^{xo}_k)\big],\\
\bar{P}^{xy}_N=-Q_N,
\end{array}
\right.
\end{eqnarray}
\begin{eqnarray}\label{1-Py-1}
\left\{\begin{array}{l}
P^y_k
= Q_k+(L_k^{yo})^TR_kL^{yo}_k+(B_kL^{yo}_k)^TP^x_{k+1}B_kL^{yo}_k\\
\quad +(D_kL^{yo}_k)^TP^x_{k+1}D_kL^{yo}_k+F_k^TP^{xy}_{k+1}B_kL_k^{yo}\\
\quad +\rho G_k^TP^{xy}_{k+1}D_kL_k^{yo}+(B_kL_k^{yo})^T(P^{xy}_{k+1})^T F_k\\
\quad +\rho(D_kL_k^{yo})^T(P^{xy}_{k+1})^T G_k+F_k^TP^y_{k+1}F_k+G_k^TP^y_{k+1}G_k,\\
P^y_N=Q_N,
\end{array}\right.
\end{eqnarray}
and
\begin{eqnarray}\label{1-barPy-1}
\left\{\begin{array}{l}
\bar{P}^y_k
= \bar{Q}_k+(L_k^{yo})^TR_k\bar{L}^{yo}_k+(\bar{L}_k^{yo})^TR_kL^{yo}_k\\
\quad  +(\bar{L}_k^{yo})^TR_k\bar{L}^{yo}_k+(L^{yo}_k+\bar{L}^{yo}_k)^T\bar{R}_k(L^{yo}_k+\bar{L}^{yo}_k)\\
\quad +(B_kL^{yo}_k)^TP^x_{k+1}\big[B_k\bar{L}^{yo}_k+\bar{B}_k(L^{yo}_k+\bar{L}^{yo}_k)\big]\\
\quad +\big[B_k\bar{L}^{yo}_k+\bar{B}_k(L^{yo}_k+\bar{L}^{yo}_k)\big]^TP^x_{k+1}B_kL^{yo}_k+\big[B_k\bar{L}^{yo}_k\\
\quad +\bar{B}_k(L^{yo}_k+\bar{L}^{yo}_k)\big]^TP^x_{k+1}\big[B_k\bar{L}^{yo}_k+\bar{B}_k(L^{yo}_k+\bar{L}^{yo}_k)\big]\\
\quad +\big[(B_k+\bar{B}_k)(L^{yo}_k+\bar{L}^{yo}_k)\big]^T\bar{P}^x_{k+1} \big[(B_k+\bar{B}_k)(L^{yo}_k\\
\quad +\bar{L}^{yo}_k)\big]+(D_kL^{yo}_k)^TP^x_{k+1}\big[D_k\bar{L}^{yo}_k+\bar{D}_k(L^{yo}_k+\bar{L}^{yo}_k)\big]\\
\quad +\big[D_k\bar{L}^{yo}_k+\bar{D}_k(L^{yo}_k+\bar{L}^{yo}_k)\big]^TP^x_{k+1}D_kL^{xo}_k+\big[D_k\bar{L}^{yo}_k\\
\end{array}\right.
\end{eqnarray}
The optimal control can be obtained by the following theorem.
\begin{theorem}\label{1-theorem-1}
For {Problem (MF-LQ)}, under the  condition
\begin{eqnarray}\label{1-condition-2}
Q_k, Q_k+\bar{Q}_k\geq 0,~k\in\bar{\mathbb{N}},~R_k, R_k+\bar{R}_k>0,~k\in\mathbb{N},
\end{eqnarray}
the unique optimal control within the class of controls of form (\ref{1-control-2}) is
\begin{eqnarray}\label{1-control-3}
\begin{array}{l}
u^o_k = \big[-(\bar{W}_k^{(1)})^{-1}(\bar{H}_k^{(1)})^T\big]x_k+\big[-(\bar{W}_k^{(2)})^{-1}(\bar{H}_k^{(2)})^T\\
\quad +(\bar{W}_k^{(1)})^{-1}(\bar{H}_k^{(1)})^T\big]\mathbb{E}x_k+\big[-(\bar{W}_k^{(1)})^{-1}(\bar{H}_k^{(3)})^T\big]y_k\\
\quad +\big[-(\bar{W}_k^{(2)})^{-1}(\bar{H}_k^{(4)})^T+(\bar{W}_k^{(1)})^{-1}(\bar{H}_k^{(3)})^T\big]\mathbb{E}y_k,\\
\equiv L^{xo}_kx_k+\bar{L}^{xo}_k\mathbb{E}x_k+L^{yo}_ky_k+\bar{L}^{yo}_k\mathbb{E}y_k,~k\in\mathbb{N},
\end{array}
\end{eqnarray}
with the property
\begin{eqnarray}\label{1-nonnegative-1}
\begin{array}{l}
P^x_k,~P^x_k+\bar{P}^x_k\geq 0,~k \in \bar{\mathbb{N}}.
\end{array}
\end{eqnarray}
\end{theorem}
\emph{Proof}. See the Appendix.

\textbf{Remark 4.1}  Taking expectation on the optimal control (\ref{1-control-3}) gives
\begin{eqnarray*}
\begin{array}{l}
\mathbb{E}u^o_k = -(\bar{W}_k^{(2)})^{-1}(\bar{H}_k^{(2)})^T\mathbb{E}x_k-(\bar{W}_k^{(2)})^{-1}(\bar{H}_k^{(4)})^T\mathbb{E}y_k.
\end{array}
\end{eqnarray*}
We combine (\ref{1-control-3}) with the expected system state equation, that is
\begin{eqnarray*}
\left\{\begin{array}{l}
\mathbb{E}x_{k+1}= (A_k+\bar{A}_k)\mathbb{E}x_k+(B_k+\bar{B}_k)\mathbb{E}u_k\\
\quad = \big[A_k+\bar{A}_k-(B_k+\bar{B}_k)(\bar{W}_k^{(2)})^{-1}(\bar{H}_k^{(2)})^T\big]\mathbb{E}x_k\\
\qquad +\big[-(B_k+\bar{B}_k)(\bar{W}_k^{(2)})^{-1}(\bar{H}_k^{(4)})^T\big]\mathbb{E}y_k\\
\quad \equiv N_k\mathbb{E}x_k+M_k\mathbb{E}y_k,\\
\mathbb{E}y_{k+1}= (F_k+\bar{F}_k)\mathbb{E}y_k\\
\quad \equiv O_k\mathbb{E}y_k,\\
\mathbb{E}x_{0}=\mathbb{E}\zeta^x,~\mathbb{E}y_{0}=\mathbb{E}\zeta^y.
\end{array}\right.
\end{eqnarray*}
Therefore, the relative open-loop expected system is
\begin{eqnarray}\label{2-remark-system-1}
\left\{\begin{array}{l}
\mathbb{E}x_{k+1}= \prod\limits_{i=1}^k N_i\mathbb{E}\zeta^x+\sum^{k}_{i=0}\prod\limits_{j=i+1}^kN_j M_i\prod\limits_{j=1}^{i-1}O_j\mathbb{E}\zeta^y,\\
\mathbb{E}y_{k+1}= \prod\limits_{i=1}^k O_i\mathbb{E}\zeta^y,
\end{array}\right.
\end{eqnarray}
and hence the optimal control is given by
\begin{eqnarray*}
\begin{array}{l}
u_k^o=\big[-(\bar{W}_k^{(1)})^{-1}(\bar{H}_k^{(1)})^T\big]x_k+\big[-(\bar{W}_k^{(2)})^{-1}(\bar{H}_k^{(2)})^T\\
\quad \ds +(\bar{W}_k^{(1)})^{-1}(\bar{H}_k^{(1)})^T\big]\cdot\big[ \prod\limits_{i=1}^k N_i\mathbb{E}\zeta^x+ \sum^{k}_{i=0}\prod\limits_{j=i+1}^k N_j M_i\\
\quad \ds \prod\limits_{j=1}^{i-1}O_j\mathbb{E}\zeta^y\big]+\big[-(\bar{W}_k^{(1)})^{-1}(\bar{H}_k^{(3)})^T\big]y_k\\
\quad \ds +\big[-(\bar{W}_k^{(2)})^{-1}(\bar{H}_k^{(4)})^T+(\bar{W}_k^{(1)})^{-1}(\bar{H}_k^{(3)})^T\big]\cdot\prod\limits_{i=1}^k O_i\mathbb{E}\zeta^y.
\end{array}
\end{eqnarray*}
It is easy to realize the control rules if $x_k$, $y_k$ and the distributions of $\zeta^x$, $\zeta^y$ are available. We remark that the distributions of $x_k$ and $y_k$ are not needed.
\subsection{Operator Linear-quadratic Problem\label{4.2}}
In this subsection, the optimal control via operator linear-quadratic theory shall be derived. Firstly, we reformulate the discrete-time operator LQ problem. A linear controlled system in abstract form rewritten from (\ref{1-system-1}) is
\begin{eqnarray}\label{2-system-1}
\left\{\begin{array}{l}
x_{k+1}=(\mathcal{A}_kx_k+\mathcal{B}_ku_k)+(\mathcal{C}_kx_k+\mathcal{D}_ku_k)w_k, \\
y_{k+1}=\mathcal{F}_ky_k+ \mathcal{G}_ky_kv_k, \\
x_0=\zeta^x, y_0=\zeta^y,
\end{array}\right.
\end{eqnarray}
with several sequences of operators
\begin{eqnarray}\label{2-Operator-1}
\left\{\begin{array}{rl}
\mathcal{A}_kz=&A_kz+\bar{A}_k\mathbb{E}z ,\\
\mathcal{C}_kz=&C_kz+\bar{C}_k\mathbb{E}z ,\\
\mathcal{F}_kz=&F_kz+\bar{F}_k\mathbb{E}z ,\\
\mathcal{G}_kz=&G_kz+\bar{G}_k\mathbb{E}z ,\\
\mathcal{B}_ku=&B_ku+\bar{B}_k\mathbb{E}u ,\\
\mathcal{D}_ku=&D_ku+\bar{D}_k\mathbb{E}u ,\\
\end{array}\right.
\end{eqnarray}
where $z\in \mathcal{Z}_k$, $u\in \mathcal{U}_k$, $\mathcal{A}_k$,$\mathcal{C}_k$,$\mathcal{F}_k$,$\mathcal{G}_k$ are from $\mathcal{Z}_k$ to $\mathcal{Z}_k$ and $\mathcal{B}_k$, $\mathcal{D}_k$ are from $\mathcal{U}_k$ to $\mathcal{U}_k$, $k\in \mathbb{N}$.  Furthermore, the performance functional in the form of inner products (\ref{1-cost-4}) has been presented in previous section. We now define operators
\begin{eqnarray}\label{2-Operator-2}
\left\{\begin{array}{l}
\mathcal{Q}_kz=(Q_k+\mathbb{E}^*\bar{Q}_k\mathbb{E})z,~~z\in \mathcal{Z}_k,~k\in \bar{\mathbb{N}},\\
\mathcal{R}_ku=(R_k+\mathbb{E}^*\bar{R}_k\mathbb{E})u,~~u\in \mathcal{U}_k,~k\in \mathbb{N}.\\
\end{array}
\right.
\end{eqnarray}
Recall that $\mathbb{E}\big[y_k^T(Q_k + \mathbb{E}^* \bar{Q}_k\mathbb{E})x_k\big]=\langle \mathcal{Q}_kx_k, y_k \rangle$ with similar meanings for related notation. Hence,
\begin{eqnarray}\label{2-cost-1}
\begin{array}{l}
J(\zeta^x,\zeta^y,u)= \ds \sum_{k=0}^{N-1}\Big(\langle \mathcal{Q}_kx_k, x_k\rangle -2\langle \mathcal{Q}_kx_k,y_k \rangle\\
\quad  + \langle \mathcal{Q}_ky_k,y_k \rangle +\langle \mathcal{R}_ku_k, u_k\rangle\Big) +\langle \mathcal{Q}_Nx_N, x_N\rangle\\
\quad -2\langle \mathcal{Q}_Nx_N,y_N \rangle + \langle \mathcal{Q}_Ny_N,y_N \rangle.
\end{array}
\end{eqnarray}
Problem (MF-LQ) in abstract form can be represented as
\begin{eqnarray}\label{2-problem-1}
\left\{\begin{array}{l}
 \mbox{Minimize}~(\ref{2-cost-1}),\\[-1mm]
\mbox{subject~to}~u\in \mathcal{U}_{ad},~\mbox{with}~(x.,y.,u.),~\mbox{satisfying}~({\ref{2-system-1}}).
\end{array}\right.
\end{eqnarray}
We then use $u^*$ to represent the optimal control for Problem (MF-LQ) in abstract form.

Next, we shall solve the problem above. Suppose we have a sequence of self-adjoint linear operators $\left\{\mathcal{P}^x_k:\mathcal{Z}_k\mapsto \mathcal{Z}_k; ~k\in \mathbb{\bar{N}}\right\}$, $\left\{\mathcal{P}^y_k:\mathcal{Z}_k\mapsto \mathcal{Z}_k; ~k\in \mathbb{\bar{N}}\right\}$ and linear operator $\left\{\mathcal{P}^{xy}_k:\mathcal{Z}_k\mapsto \mathcal{Z}_k; ~k\in \mathbb{\bar{N}}\right\}$ determined by
\begin{eqnarray*}
\left\{\begin{array}{rl}
\mathcal{P}^x_k\mathcal{Z}_l\equiv& \left\{\mathcal{P}^x_kz\mid z\in \mathcal{Z}_l\right\}\subseteq \mathcal{Z}_l,~~l\leq k,\\
\mathcal{P}^{xy}_k\mathcal{Z}_l\equiv& \left\{\mathcal{P}^{xy}_kz\mid z\in \mathcal{Z}_l\right\}\subseteq \mathcal{Z}_l,~~l\leq k,\\
\mathcal{P}^y_k\mathcal{Z}_l\equiv& \left\{\mathcal{P}^y_kz\mid z\in \mathcal{Z}_l\right\}\subseteq \mathcal{Z}_l,~~l\leq k.\\
\end{array}\right.
\end{eqnarray*}
where $\mathcal{Z}_l\subseteq\mathcal{Z}_k$, $l\leq k$. Clearly,
\begin{eqnarray}\label{2-Px-1}
\begin{array}{l}
\langle \mathcal{P}^x_{k+1}x_{k+1},x_{k+1}\rangle\\
=\big\langle \mathcal{P}^x_{k+1}\big[(\mathcal{A}_kx_k+\mathcal{B}_ku_k)+(\mathcal{C}_kx_k+\mathcal{D}_ku_k)w_k\big],\\
\quad \big[(\mathcal{A}_kx_k+\mathcal{B}_ku_k)+(\mathcal{C}_kx_k+\mathcal{D}_ku_k)w_k\big] \big\rangle\\
=\mathbb{E}\big[(\mathcal{A}_kx_k+\mathcal{B}_ku_k)^T\mathcal{P}^x_{k+1}(\mathcal{A}_kx_k+\mathcal{B}_ku_k)\\
\quad +((\mathcal{C}_kx_k+\mathcal{D}_ku_k)w_k)^T\mathcal{P}^x_{k+1}(\mathcal{A}_kx_k+\mathcal{B}_ku_k)\big]\\
\quad +\mathbb{E}\big[(\mathcal{A}_kx_k+\mathcal{B}_ku_k)^T\mathcal{P}^x_{k+1}(\mathcal{C}_kx_k+\mathcal{D}_ku_k)w_k\\
\quad +((\mathcal{C}_kx_k+\mathcal{D}_ku_k)w_k)^T\mathcal{P}^x_{k+1}(\mathcal{C}_kx_k+\mathcal{D}_ku_k)w_k\big]\\
=\mathbb{E}\big[(\mathcal{A}_kx_k+\mathcal{B}_ku_k)^T\mathcal{P}^x_{k+1}(\mathcal{A}_kx_k+\mathcal{B}_ku_k)\big]\\
\quad +\mathbb{E}\big[(\mathcal{C}_kx_k+\mathcal{D}_ku_k)^T\mathcal{P}^x_{k+1}(\mathcal{C}_kx_k+\mathcal{D}_ku_k)\big]\\
=\big\langle \mathcal{P}^x_{k+1}(\mathcal{A}_kx_k+\mathcal{B}_ku_k), (\mathcal{A}_kx_k+\mathcal{B}_ku_k) \rangle\\
\quad +\langle \mathcal{P}^x_{k+1}(\mathcal{C}_kx_k+\mathcal{D}_ku_k), (\mathcal{C}_kx_k+\mathcal{D}_ku_k) \big\rangle,
\end{array}
\end{eqnarray}
\begin{eqnarray}\label{2-Pxy-1}
\begin{array}{l}
\langle \mathcal{P}^{xy}_{k+1}x_{k+1},y_{k+1}\rangle=\big\langle \mathcal{P}^{xy}_{k+1}\big[(\mathcal{A}_kx_k+\mathcal{B}_ku_k)\\
\quad +(\mathcal{C}_kx_k+\mathcal{D}_ku_k)w_k\big],\big[\mathcal{F}_ky_k+ \mathcal{G}_ky_kv_k\big] \big\rangle\\
=\mathbb{E}\big[(\mathcal{F}_ky_k)^T\mathcal{P}^{xy}_{k+1}(\mathcal{A}_kx_k+\mathcal{B}_ku_k)\\
\quad +(\mathcal{G}_ky_kv_k)^T\mathcal{P}^{xy}_{k+1}(\mathcal{A}_kx_k+\mathcal{B}_ku_k)\big]\\
\quad +\mathbb{E}\big[(\mathcal{F}_ky_k)^T\mathcal{P}^{xy}_{k+1}(\mathcal{C}_kx_k+\mathcal{C}_ku_k)w_k\\
\quad +(\mathcal{G}_ky_kv_k)^T\mathcal{P}^{xy}_{k+1}(\mathcal{C}_kx_k+\mathcal{D}_ku_k)w_k\big]\\
=\mathbb{E}\big[(\mathcal{F}_ky_k)^T\mathcal{P}^{xy}_{k+1}(\mathcal{A}_kx_k+\mathcal{B}_ku_k)\big]\\
\quad +\rho\mathbb{E}\big[(\mathcal{G}_ky_k)^T\mathcal{P}^{xy}_{k+1}(\mathcal{C}_kx_k+\mathcal{D}_ku_k)\big]\\
=\big\langle \mathcal{P}^{xy}_{k+1}(\mathcal{A}_kx_k+\mathcal{B}_ku_k), \mathcal{F}_ky_k \rangle\\
\quad +\rho\langle \mathcal{P}^{xy}_{k+1}(\mathcal{C}_kx_k+\mathcal{D}_ku_k), \mathcal{G}_ky_k \big\rangle,
\end{array}
\end{eqnarray}
\begin{eqnarray}\label{2-Py-1}
\begin{array}{l}
\langle \mathcal{P}^y_{k+1}y_{k+1},y_{k+1}\rangle\\
=\big\langle \mathcal{P}^y_{k+1}(\mathcal{F}_ky_k+ \mathcal{G}_ky_kv_k), (\mathcal{F}_ky_k+ \mathcal{G}_ky_kv_k) \big\rangle\\
=\mathbb{E}\big[(\mathcal{F}_ky_k)^T\mathcal{P}^y_{k+1}\mathcal{F}_ky_k+(\mathcal{G}_ky_kv_k)^T\mathcal{P}^y_{k+1}\mathcal{F}_ky_k\big]\\
\quad +\mathbb{E}\big[(\mathcal{F}_ky_k)^T\mathcal{P}^y_{k+1}\mathcal{G}_ky_kv_k+(\mathcal{G}_ky_kv_k)^T\mathcal{P}^y_{k+1}\mathcal{G}_ky_kv_k\big]\\
=\mathbb{E}\big[(\mathcal{F}_ky_k)^T\mathcal{P}^y_{k+1}\mathcal{F}_ky_k\big]+\mathbb{E}\big[(\mathcal{G}_ky_k)^T\mathcal{P}^y_{k+1}\mathcal{G}_ky_k\big]\\
=\big\langle \mathcal{P}^y_{k+1}\mathcal{F}_ky_k, \mathcal{F}_ky_k \rangle +\langle \mathcal{P}^y_{k+1}\mathcal{G}_ky_k, \mathcal{G}_ky_k \big\rangle.
\end{array}
\end{eqnarray}
We now consider the problem by backward recursion with only $k\in\{N-1,N\}$. By (\ref{2-Px-1})--(\ref{2-Py-1}), we have
\begin{eqnarray}\label{2-P-1}
\begin{array}{l}
\langle \mathcal{P}^x_{N}x_{N},x_{N}\rangle =\big\langle \mathcal{P}^x_{N}(\mathcal{A}_{N-1}x_k+\mathcal{B}_{N-1}u_{N-1}),\\
\qquad (\mathcal{A}_{N-1}x_{N-1}+\mathcal{B}_{N-1}u_{N-1}) \big\rangle +\big\langle \mathcal{P}^x_{N}(\mathcal{C}_{N-1}x_{N-1}\\
\qquad +\mathcal{D}_{N-1}u_{N-1}),(\mathcal{C}_{N-1}x_{N-1}+\mathcal{D}_{N-1}u_{N-1}) \big\rangle,\\
\langle \mathcal{P}^{xy}_{N}x_{N},y_{N}\rangle\\
=\big\langle \mathcal{P}^{xy}_{N}(\mathcal{A}_{N-1}x_{N-1}+\mathcal{B}_{N-1}u_{N-1}), \mathcal{F}_{N-1}y_{N-1} \big\rangle\\
\qquad +\rho\big\langle \mathcal{P}^{xy}_{N}(\mathcal{C}_{N-1}x_{N-1}+\mathcal{D}_{N-1}u_{N-1}), \mathcal{G}_{N-1}y_{N-1} \big\rangle,\\
\langle \mathcal{P}^y_{N}y_{N},y_{N}\rangle =\big\langle \mathcal{P}^y_{N}\mathcal{F}_{N-1}y_{N-1}, \mathcal{F}_{N-1}y_{N-1} \big\rangle\\
\qquad +\big\langle \mathcal{P}^y_{N}\mathcal{G}_{N-1}y_{N-1}, \mathcal{G}_{N-1}y_{N-1} \big\rangle.
\end{array}
\end{eqnarray}
Let $\mathcal{P}^x_{N},-\mathcal{P}^{xy}_{N},~\mathcal{P}^y_{N}=\mathcal{Q}_{N}$.  By taking (\ref{2-P-1}) into $J_{N-1}^N$, we have
\begin{eqnarray}\label{2-cost-2}
\begin{array}{l}
J_{N-1}^N(x_{N-1}, y_{N-1}, u)\\
= \langle \mathcal{Q}_{N-1}x_{N-1}, x_{N-1}\rangle -2\langle \mathcal{Q}_{N-1}x_{N-1},y_{N-1} \rangle\\
\quad + \langle \mathcal{Q}_{N-1}y_{N-1},y_{N-1} \rangle +\langle \mathcal{R}_{N-1}u_{N-1}, u_{N-1}\rangle\\
\quad +\langle \mathcal{Q}_Nx_N, x_N\rangle -2\langle \mathcal{Q}_Nx_N,y_N \rangle + \langle \mathcal{Q}_Ny_N,y_N \rangle\\
=\langle \mathcal{Q}_{N-1}x_{N-1}, x_{N-1}\rangle -2\langle \mathcal{Q}_{N-1}x_{N-1},y_{N-1} \rangle\\
\quad + \langle \mathcal{Q}_{N-1}y_{N-1},y_{N-1} \rangle +\langle \mathcal{R}_{N-1}u_{N-1}, u_{N-1}\rangle\\
\quad +\big\langle \mathcal{P}^x_{N}(\mathcal{A}_{N-1}x_{N-1}+\mathcal{B}_{N-1}u_{N-1}), \\
\quad (\mathcal{A}_{N-1}x_{N-1}+\mathcal{B}_{N-1}u_{N-1}) \big\rangle +\big\langle \mathcal{P}^x_{N}(\mathcal{C}_{N-1}x_{N-1}\\
\quad +\mathcal{D}_{N-1}u_{N-1}), (\mathcal{C}_{N-1}x_{N-1} +\mathcal{D}_{N-1}u_{N-1}) \big\rangle\\
\quad +2\big\langle \mathcal{P}^{xy}_{N}(\mathcal{A}_{N-1}x_{N-1}+\mathcal{B}_{N-1}u_{N-1}), \mathcal{F}_{N-1}y_{N-1} \big\rangle\\
\quad +2\rho\big\langle \mathcal{P}^{xy}_{N}(\mathcal{C}_{N-1}x_{N-1}+\mathcal{D}_{N-1}u_{N-1}), \mathcal{G}_{N-1}y_{N-1} \big\rangle\\
\quad +\big\langle \mathcal{P}^y_{N}\mathcal{F}_{N-1}y_{N-1}, \mathcal{F}_{N-1}y_{N-1} \big\rangle\\
\quad +\big\langle \mathcal{P}^y_{N}\mathcal{G}_{N-1}y_{N-1}, \mathcal{G}_{N-1}y_{N-1} \big\rangle\\
=2\langle \Theta_{1,N-1}u_{N-1},x_{N-1} \rangle +\langle \Theta_{2,N-1}u_{N-1},u_{N-1} \rangle\\
\quad +2\langle \Theta_{3,N-1}u_{N-1},y_{N-1} \rangle +\langle \Lambda_{1,N-1}x_{N-1},x_{N-1} \rangle\\
\quad +2\langle \Lambda_{2,N-1}x_{N-1},y_{N-1} \rangle +\langle \Lambda_{3,N-1}y_{N-1},y_{N-1} \rangle,
\end{array}
\end{eqnarray}
where
\begin{eqnarray}\label{2-TL-1}
\left\{\begin{array}{l}
\Theta_{1,N-1}=\mathcal{A}_{N-1}^T\mathcal{P}^x_N\mathcal{B}_{N-1}+\mathcal{C}_{N-1}^T\mathcal{P}^x_N\mathcal{D}_{N-1},\\
\Theta_{2,N-1}=\mathcal{R}_{N-1}+\mathcal{B}_{N-1}^T\mathcal{P}^x_N\mathcal{B}_{N-1}+\mathcal{D}_{N-1}^T\mathcal{P}^x_N\mathcal{D}_{N-1},\\
\Theta_{3,N-1}=\mathcal{F}_{N-1}^T\mathcal{P}^{xy}_N\mathcal{B}_{N-1}+\rho\mathcal{G}_{N-1}^T\mathcal{P}^{xy}_N\mathcal{D}_{N-1},\\
\Lambda_{1,N-1}=\mathcal{Q}_{N-1}+\mathcal{A}_{N-1}^T\mathcal{P}^x_N\mathcal{A}_{N-1}+\mathcal{C}_{N-1}^T\mathcal{P}^x_N\mathcal{C}_{N-1},\\
\Lambda_{2,N-1}=-\mathcal{Q}_{N-1}+\mathcal{F}_{N-1}^T\mathcal{P}^{xy}_N\mathcal{A}_{N-1}+\rho\mathcal{G}_{N-1}^T\mathcal{P}^{xy}_N\mathcal{C}_{N-1},\\
\Lambda_{3,N-1}=\mathcal{Q}_{N-1}+\mathcal{F}_{N-1}^T\mathcal{P}^y_N\mathcal{F}_{N-1}+\mathcal{G}_{N-1}^T\mathcal{P}^y_N\mathcal{G}_{N-1}.\\
\end{array}\right.
\end{eqnarray}
If $\Theta_{2,N-1}$ is positive definite and self-adjoint, (\ref{2-cost-2}) can then be rewritten as
\begin{eqnarray}\label{2-cost-5}
\begin{array}{l}
J^N_{N-1}(x_{N-1},y_{N-1},u)\\
=\big\langle (\Lambda_{1,N-1}-\Theta_{1,N-1}\Theta_{2,N-1}^{-1}\Theta_{1,N-1}^*)x_{N-1},x_{N-1} \big\rangle\\
\quad +2\big\langle (\Lambda_{2,N-1}-\Theta_{3,N-1}\Theta_{2,N-1}^{-1}\Theta_{1,N-1}^*)x_{N-1},y_{N-1} \big\rangle\\
\quad +\big\langle (\Lambda_{3,N-1}-\Theta_{3,N-1}\Theta_{2,N-1}^{-1}\Theta_{3,N-1}^*)y_{N-1},y_{N-1} \big\rangle\\
\quad +\big\langle \Theta_{2,N-1}(u_{N-1}+\Theta_{2,N-1}^{-1}\Theta_{1,N-1}^*x_{N-1}+\Theta_{2,N-1}^{-1}\Theta_{3,N-1}^*\\
\quad \cdot y_{N-1}),(u_{N-1}+\Theta_{2,N-1}^{-1}\Theta_{1,N-1}^*x_{N-1}+\Theta_{2,N-1}^{-1}\Theta_{3,N-1}^*y_{N-1}) \big\rangle.\\
\end{array}
\end{eqnarray}
Let
\begin{eqnarray}\label{2-P-5}
\left\{\begin{array}{l}
\mathcal{P}^x_{N-1}=\Lambda_{1,N-1}-\Theta_{1,N-1}\Theta_{2,N-1}^{-1}\Theta_{1,N-1}^*,\\
\mathcal{P}^{xy}_{N-1}=\Lambda_{2,N-1}-\Theta_{3,N-1}\Theta_{2,N-1}^{-1}\Theta_{1,N-1}^*,\\
\mathcal{P}^y_{N-1}=\Lambda_{3,N-1}-\Theta_{3,N-1}\Theta_{2,N-1}^{-1}\Theta_{3,N-1}^*,
\end{array}\right.
\end{eqnarray}
and the cost functional (\ref{2-cost-5}) may achieve its minimum if we select 
\begin{eqnarray}\label{2-control-2}
\begin{array}{l}
u_{N-1}^*=-\Theta_{2,N-1}^{-1}\Theta_{1,N-1}^*x_{N-1}-\Theta_{2,N-1}^{-1}\Theta_{3,N-1}^*y_{N-1},
\end{array}
\end{eqnarray}
where the minimum is
\begin{eqnarray}\label{2-cost-3}
\begin{array}{l}
J^N_{N-1}(x_{N-1},y_{N-1},u^*)=\langle \mathcal{P}^x_{N-1}x_{N-1},x_{N-1} \rangle \\
\quad +2\langle \mathcal{P}^{xy}_{N-1}x_{N-1},y_{N-1} \rangle +\langle \mathcal{P}^y_{N-1}x_{N-1},x_{N-1} \rangle .
\end{array}
\end{eqnarray}

For this, we reach to the following lemma, which gives a compact form of $\mathcal{P}^x_{N-1}$, $\mathcal{P}^{xy}_{N-1}$ and $\mathcal{P}^y_{N-1}$.
\begin{lemma}\label{1-lemma-1}
If
\begin{eqnarray}\label{2-condition-1}
\begin{array}{l}
Q_{N-1}, Q_{N-1}+\bar{Q}_{N-1}\geq 0,~Q_N,Q_N+\bar{Q}_N\geq 0,\\
R_{N-1}, R_{N-1}+\bar{R}_{N-1}>0,
\end{array}
\end{eqnarray}
then $\mathcal{P}^x_{N-1}$, $\mathcal{P}^{xy}_{N-1}$, $\mathcal{P}^y_{N-1}$ defined in (\ref{2-P-1}) have the following form
\begin{eqnarray}\label{2-P-2}
\left\{\begin{array}{l}
\mathcal{P}^x_{N-1}=(I-\mathbb{E}^*)S^x_{N-1}(I-\mathbb{E})+\mathbb{E}^*T^x_{N-1}\mathbb{E},\\
\mathcal{P}^{xy}_{N-1}=(I-\mathbb{E}^*)S^{xy}_{N-1}(I-\mathbb{E})+\mathbb{E}^*T^{xy}_{N-1}\mathbb{E},\\
\mathcal{P}^y_{N-1}=(I-\mathbb{E}^*)S^y_{N-1}(I-\mathbb{E})+\mathbb{E}^*T^y_{N-1}\mathbb{E}.
\end{array}\right.
\end{eqnarray}
where
\begin{eqnarray*}
\left\{\begin{array}{l}
S_{N-1}^x=Q_{N-1}+A_{N-1}^T Q_{N} A_{N-1}+C_{N-1}^T Q_{N} C_{N-1}\\
\qquad -{H}_{N-1}^{(1)}({W}_{N-1}^{(1)})^{-1}({H}_{N-1}^{(1)})^T,\\
T_{N-1}^x=Q_{N-1}+\bar{Q}_{N-1}+(A_{N-1}+\bar{A}_{N-1})^T (Q_N+\bar{Q}_N)\\
\qquad \cdot(A_{N-1}+\bar{A}_{N-1})+(C_{N-1}+\bar{C}_{N-1})^T Q_{N}\\
\qquad \cdot(C_{N-1}+\bar{C}_{N-1})-{H}_{N-1}^{(2)}({W}_{N-1}^{(2)})^{-1}({H}_{N-1}^{(2)})^T,\\
S_{N-1}^{xy}=-Q_{N-1}+F_{N-1}^T Q_{N} A_{N-1}+\rho G_{N-1}^T Q_{N} C_{N-1}\\
\qquad -{H}_{N-1}^{(3)}({W}_{N-1}^{(1)})^{-1}({H}_{N-1}^{(1)})^T,\\
T_{N-1}^{xy}=-Q_{N-1}-\bar{Q}_{N-1}+(F_{N-1}+\bar{F}_{N-1})^T (Q_N+\bar{Q}_N)\\
\qquad \cdot(A_{N-1}+\bar{A}_{N-1})+\rho(G_{N-1}+\bar{G}_{N-1})^T Q_{N}\\
\qquad \cdot(C_{N-1}+\bar{C}_{N-1})-{H}_{N-1}^{(4)}({W}_{N-1}^{(2)})^{-1}({H}_{N-1}^{(2)})^T, \\
S_{N-1}^y=Q_{N-1}+F_{N-1}^T Q_{N} F_{N-1}+G_{N-1}^T Q_{N} G_{N-1}\\
\qquad -{H}_{N-1}^{(3)}({W}_{N-1}^{(1)})^{-1}({H}_{N-1}^{(3)})^T,\\
T_{N-1}^y=Q_{N-1}+\bar{Q}_{N-1}+(F_{N-1}+\bar{F}_{N-1})^T (Q_N+\bar{Q}_N)\\
\qquad \cdot(F_{N-1}+\bar{F}_{N-1})+(G_{N-1}+\bar{G}_{N-1})^T Q_{N}\\
\qquad \cdot(G_{N-1}+\bar{G}_{N-1})-{H}_{N-1}^{(4)}({W}_{N-1}^{(2)})^{-1}({H}_{N-1}^{(4)})^T,\\
S_N^x=Q_N,~T_N^x=Q_N+\bar{Q}_N,~S_N^{xy}=-Q_N,\\
T_N^{xy}=-Q_N-\bar{Q}_N,~S_N^y=Q_N,~T_N^y=Q_N+\bar{Q}_N.
\end{array}\right.
\end{eqnarray*}
with
\begin{eqnarray*}
\left\{\begin{array}{l}
W_{N-1}^{(1)}=R_{N-1}+B_{N-1}^TQ_{N}B_{N-1}+D_{N-1}^TQ_{N}D_{N-1},\\
W_{N-1}^{(2)}=R_{N-1}+\bar{R}_{N-1}+(B_{N-1}+\bar{B}_{N-1})^T(Q_N+\bar{Q}_N)\\
\cdot(B_{N-1}+\bar{B}_{N-1})+(D_{N-1}+\bar{D}_{N-1})^TQ_N(D_{N-1}+\bar{D}_{N-1}),\\
H_{N-1}^{(1)}=A_{N-1}^TQ_{N}B_{N-1}+C_{N-1}^TQ_{N}D_{N-1},\\
H_{N-1}^{(2)}=(A_{N-1}+\bar{A}_{N-1})^T(Q_N+\bar{Q}_N)(B_{N-1}+\bar{B}_{N-1})\\
\quad +(C_{N-1}+\bar{C}_{N-1})^TQ_{N}(D_{N-1}+\bar{D}_{N-1}),\\
H_{N-1}^{(3)}=F_{N-1}^TQ_{N}B_{N-1}+\rho G_{N-1}^TQ_{N}D_{N-1},\\
H_{N-1}^{(4)}=(F_{N-1}+\bar{F}_{N-1})^T (Q_N+\bar{Q}_N) (B_{N-1}+\bar{B}_{N-1})\\
\quad +\rho(G_{N-1}+\bar{G}_{N-1})^TQ_{N}(D_{N-1}+\bar{D}_{N-1}),
\end{array}\right.
\end{eqnarray*}
\end{lemma}
\emph{Proof}. See the Appendix.

We now express the optimal control using Riccati difference equations by the following theorem.
\begin{theorem}\label{2-theorem-1}
Let
\begin{eqnarray}\label{2-condition-5}
Q_k, Q_k+\bar{Q}_k\geq 0,,~k\in\bar{\mathbb{N}},~R_k, R_k+\bar{R}_k>0,
\end{eqnarray}
and introduce Riccati equations
\begin{eqnarray}\label{2-ST-2}
\left\{\begin{array}{l}
S_k^x=Q_k+A_k^TS_{k+1}^xA_k+C_k^TS_{k+1}^xC_k\\
\quad -{H}_k^{(1)}({W}_k^{(1)})^{-1}({H}_k^{(1)})^T,\\
T_k^x=Q_k+\bar{Q}_k+(A_k+\bar{A}_k)^TT_{k+1}^x(A_k+\bar{A}_k)+(C_k\\
\quad +\bar{C}_k)^TS_{k+1}^x(C_k+\bar{C}_k)-{H}_k^{(2)}({W}_k^{(2)})^{-1}({H}_k^{(2)})^T,\\
S_k^{xy}=-Q_k+F_k^TS_{k+1}^{xy} A_k+\rho G_k^TS_{k+1}^{xy}C_k\\
\quad -{H}_k^{(3)}({W}_k^{(1)})^{-1}({H}_k^{(1)})^T,\\
T_k^{xy}=-Q_k-\bar{Q}_k+(F_k+\bar{F}_k)^TT_{k+1}^{xy}(A_k+\bar{A}_k)+\rho(G_k\\
\quad +\bar{G}_k)^TS_{k+1}^{xy}(C_k+\bar{C}_k)-{H}_k^{(4)}({W}_k^{(2)})^{-1}({H}_k^{(2)})^T,\\
S_k^y=Q_k+F_k^TS_{k+1}^yF_k+G_k^TS_{k+1}^yG_k\\
\quad -{H}_k^{(3)}({W}_k^{(1)})^{-1}({H}_k^{(3)})^T,\\
T_k^y=Q_k+\bar{Q}_k+(F_k+\bar{F}_k)^TT_{k+1}^y(F_k+\bar{F}_k)+(G_k\\
\quad +\bar{G}_k)^TS_{k+1}^y(G_k+\bar{G}_k)-{H}_k^{(4)}({W}_k^{(2)})^{-1}({H}_k^{(4)})^T,\\
S_N^x=Q_N,~T_N^x=Q_N+\bar{Q}_N,~S_N^{xy}=-Q_N,\\
T_N^{xy}=-Q_N-\bar{Q}_N,~S_N^y=Q_N,~T_N^y=Q_N+\bar{Q}_N.
\end{array}\right.
\end{eqnarray}
with
\begin{eqnarray}\label{2-WH-1}
\left\{
\begin{array}{l}
W_k^{(1)}=R_k+B_k^TS^x_{k+1}B_k+D_k^TS^x_{k+1}D_k,\\
W_k^{(2)}=R_k+\bar{R}_k+(B_k+\bar{B}_k)^TT^x_{k+1}(B_k+\bar{B}_k)\\
\qquad +(D_k+\bar{D}_k)^TS^x_{k+1}(D_k+\bar{D}_k),\\
H_k^{(1)}=A_k^TS^x_{k+1}B_k+C_k^TS^x_{k+1}D_k,\\
H_k^{(2)}=(A_k+\bar{A}_k)^TT^x_{k+1}(B_k+\bar{B}_k)\\
\qquad +(C_k+\bar{C}_k)^TS^x_{k+1}(D_k+\bar{D}_k),\\
H_k^{(3)}= F_k^TS^{xy}_{k+1}B_k+\rho G_k^TS^{xy}_{k+1}D_k,\\
H_k^{(4)}=(F_k+\bar{F}_k)^TT^{xy}_{k+1}(B_k+\bar{B}_k)\\
\qquad +\rho(G_k+\bar{G}_k)^TS^{xy}_{k+1}(D_k+\bar{D}_k).
\end{array}
\right.
\end{eqnarray}
The unique optimal control for Problem (MF-LQ) is
\begin{eqnarray}\label{2-control-1}
\begin{array}{l}
u_k^*=-({W}_k^{(1)})^{-1}({H}_k^{(1)})^T(x_k-\mathbb{E}x_k)\\
\quad -({W}_k^{(2)})^{-1}({H}_k^{(2)})^T\mathbb{E}x_k-({W}_k^{(1)})^{-1}({H}_k^{(3)})^T(y_k-\mathbb{E}y_k)\\
\quad -({W}_k^{(2)})^{-1}({H}_k^{(4)})^T\mathbb{E}y_k,~k\in\mathbb{N}.
\end{array}
\end{eqnarray}
\end{theorem}
\emph{Proof}. See the Appendix.

\subsection{Optimal Control via Completing the Square\label{4.3}}
In the previous subsections, we present two methods to solve the Problem (MF-LQ) and obtain the optimal controls. In this subsection, we shall show that controls obtained in last two subsections coincide, for which the extension of classical method of completing the square of LQ problem will be used. Firstly, (\ref{1-Px-1})--(\ref{1-barPy-1}) can be transformed to
\begin{eqnarray}\label{2-Px-2}
\begin{array}{l}
P^x_k = Q_k+A_k^TP^x_{k+1}A_k+C_k^TP^x_{k+1}C_k\\
\quad +(L^{xo}_k)^T\big(R_k+B_k^TP^x_{k+1}B_k+D_k^TP^x_{k+1}D_k\big)L^{xo}_k\\
\quad +(L^{xo}_k)^T\big(B_k^TP^x_{k+1}A_k+D_k^TP^x_{k+1}C_k\big)\\
\quad +\big(A_k^TP^x_{k+1}B_k+C_k^TP^x_{k+1}D_k\big) L^{xo}_k\\
= Q_k+A_k^TP^x_{k+1}A_k+C_k^TP^x_{k+1}C_k-\bar{H}_k^{(1)}(\bar{W}_k^{(1)})^{-1}(\bar{H}_k^{(1)})^T,
\end{array}
\end{eqnarray}
\begin{eqnarray}\label{2-Px-3}
\begin{array}{l}
P^x_k+\bar{P}^x_k = Q_k+\bar{Q}_k+(C_k+\bar{C_k})^TP^x_{k+1}(C_k+\bar{C}_k)\\
\quad +(A_k+\bar{A}_k)^T(P^x_{k+1}+\bar{P}^x_{k+1})(A_k+\bar{A}_k)+(L^{xo}_k+\bar{L}^{xo}_k)^T\\
\quad \cdot\big[R_k+\bar{R}_k+(B_k+\bar{B}_k)^T(P^x_{k+1}+\bar{P}^x_{k+1})(B_k+\bar{B}_k)\\
\quad +(D_k+\bar{D}_k)^TP^x_{k+1}(D_k+\bar{D}_k)\big](L^{xo}_k+\bar{L}^{xo}_k)\\
\quad +(L^{xo}_k+\bar{L}^{xo}_k)^T\big[(D_k+\bar{D}_k)^TP^x_{k+1}(C_k+\bar{C}_k)\\
\quad +(B_k+\bar{B}_k)^T(P^x_{k+1}+\bar{P}^x_{k+1})(A_k+\bar{A}_k)\big]\\
\quad +\big[(A_k+\bar{A})^T(P^x_{k+1}+\bar{P}^x_{k+1})(B_k+\bar{B}_k)\\
\quad +(C_k+\bar{C}_k)^TP^x_{k+1}(D_k+\bar{D}_k)\big](L^{xo}_k+\bar{L}^{xo}_k)\\
=Q_k+\bar{Q}_k+(A_k+\bar{A}_k)^T(P^x_{k+1}+\bar{P}^x_{k+1})(A_k+\bar{A_k})\\
\quad +(C_k+\bar{C_k})^TP^x_{k+1}(C_k+\bar{C}_k)-\bar{H}_k^{(2)}(\bar{W}_k^{(2)})^{-1}(\bar{H}_k^{(2)})^T,
\end{array}
\end{eqnarray}
\begin{eqnarray}\label{2-Pxy-2}
\begin{array}{l}
P^{xy}_k = -Q_k+F_k^TP^{xy}_{k+1} A_k+\rho G_k^TP^{xy}_{k+1}C_k\\
\quad +(L^{yo}_k)^T\big(R_k+B_k^TP^x_{k+1}B_k+D_k^TP^x_{k+1}D_k\big)L^{xo}_k\\
\quad +(L^{yo}_k)^T\big(B_k^TP^{x}_{k+1}A_k+ D_k^TP^{x}_{k+1}C_k\big)\\
\quad +\big(F_k^TP^{xy}_{k+1}B_k+\rho G_k^TP^{xy}_{k+1}D_k\big)L^{xo}_k\\
= -Q_k+F_k^TP^{xy}_{k+1}A_k+\rho G_k^TP^{xy}_{k+1}C_k-\bar{H}_k^{(3)}(\bar{W}_k^{(1)})^{-1}(\bar{H}_k^{(1)})^T,
\end{array}
\end{eqnarray}
\begin{eqnarray}\label{2-Pxy-3}
\begin{array}{l}
P^{xy}_k+\bar{P}^{xy}_k = -Q_k-\bar{Q}_k+\rho(G_k+\bar{G_k})^TP^{xy}_{k+1}(C_k+\bar{C}_k)\\
\quad +(F_k+\bar{F}_k)^T(P^{xy}_{k+1}+\bar{P}^{xy}_{k+1})(A_k+\bar{A}_k)+(L^{yo}_k+\bar{L}^{yo}_k)^T\\
\quad \cdot\big[R_k+\bar{R}_k+(B_k+\bar{B}_k)^T(P^x_{k+1}+\bar{P}^x_{k+1})(B_k+\bar{B}_k)\\
\quad +(D_k+\bar{D}_k)P^x_{k+1}(D_k+\bar{D}_k)\big](L^{xo}_k+\bar{L}^{xo}_k)\\
\quad +(L^{yo}_k+\bar{L}^{yo}_k)^T\big[(D_k+\bar{D}_k)^TP^{x}_{k+1}(C_k+\bar{C}_k)\\
\quad +(B_k+\bar{B}_k)^T(P^{x}_{k+1}+\bar{P}^{x}_{k+1})(A_k+\bar{A}_k)\big]\\
\quad +\big[(F_k+\bar{F})^T(P^{xy}_{k+1}+\bar{P}^{xy}_{k+1})(B_k+\bar{B}_k)\\
\quad +\rho(G_k+\bar{G}_k)^TP^{xy}_{k+1}(D_k+\bar{D}_k)\big](L^{xo}_k+\bar{L}^{xo}_k)\\
= -Q_k-\bar{Q}_k+(F_k+\bar{F}_k)^T(P^{xy}_{k+1}+\bar{P}^{xy}_{k+1})(A_k+\bar{A_k})\\\quad +\rho(G_k+\bar{G}_k)^TP^{xy}_{k+1}(C_k+\bar{C}_k)-\bar{H}_k^{(4)}(\bar{W}_k^{(2)})^{-1}(\bar{H}_k^{(2)})^T,
\end{array}
\end{eqnarray}
\begin{eqnarray}\label{2-Py-2}
\begin{array}{l}
P^y_k = Q_k+F_k^TP^{y}_{k+1}F_k+G_k^TP^y_{k+1}G_k\\
\quad +(L^{yo}_k)^T\big(R_k+B_k^TP^x_{k+1}B_k+D_k^TP^x_{k+1}D_k\big)L^{yo}_k\\
 +(L^{yo}_k)^T\big[B_k^T(P^{xy}_{k+1})^T F_k+\rho D_k^T(P^{xy}_{k+1})^T G_k\big]\\
\quad +\big(F_k^T P^{xy}_{k+1} B_k+\rho G_k^TP^{xy}_{k+1}D_k\big)L^{yo}_k\\
= Q_k+F_k^TP^{y}_{k+1}F_k+G_k^TP^y_{k+1}G_k-\bar{H}_k^{(3)}(\bar{W}_k^{(1)})^{-1}(\bar{H}_k^{(3)})^T,
\end{array}
\end{eqnarray}
\begin{eqnarray}\label{2-Py-3}
\begin{array}{l}
P^y_k+\bar{P}^y_k = Q_k+\bar{Q}_k+(G_k+\bar{G_k})^TP^y_{k+1}(G_k+\bar{G}_k)\\
\quad +(F_k+\bar{F}_k)^T(P^y_{k+1}+\bar{P}^y_{k+1})(F_k+\bar{F}_k)+(L^{yo}_k+\bar{L}^{yo}_k)^T\\
\quad \cdot\big[R_k+\bar{R}_k+(B_k+\bar{B}_k)^T(P^x_{k+1}+\bar{P}^x_{k+1})(B_k+\bar{B}_k)\\
\quad +(D_k+\bar{D}_k)^TP^x_{k+1}(D_k+\bar{D}_k)\big](L^{yo}_k+\bar{L}^{yo}_k)\\
\quad +(L^{yo}_k+\bar{L}^{yo}_k)^T\big[\rho(D_k+\bar{D}_k)^T(P^{xy}_{k+1})^T(G_k+\bar{G}_k)\\
\quad +(B_k+\bar{B}_k)^T(P^{xy}_{k+1}+\bar{P}^{xy}_{k+1})^T(F_k+\bar{F}_k)\big]\\
\quad +\big[(F_k+\bar{F})^T(P^{xy}_{k+1}+\bar{P}^{xy}_{k+1})(B_k+\bar{B}_k)\\
\quad +\rho(G_k+\bar{G}_k)^TP^{xy}_{k+1}(D_k+\bar{D}_k)\big](L^{yo}_k+\bar{L}^{yo}_k)\\
= Q_k+\bar{Q}_k+(F_k+\bar{F}_k)^T(P^y_{k+1}+\bar{P}^y_{k+1})(F_k+\bar{F_k})\\
\quad +(G_k+\bar{G_k})^TP^y_{k+1}(G_k+\bar{G}_k)-\bar{H}_k^{(4)}(\bar{W}_k^{(2)})^{-1}(\bar{H}_k^{(4)})^T.
\end{array}
\end{eqnarray}
Notice that the boundary conditions, (\ref{2-Px-2})--(\ref{2-Py-3}) are those of (\ref{2-ST-2}) with
\begin{eqnarray*}
\left\{\begin{array}{l}
\bar{W}^{(1)}_k=W^{(1)}_k, ~\bar{W}^{(2)}_k=W^{(2)}_k, ~\bar{H}^{(1)}_k=H^{(1)}_k,\\
\bar{H}^{(2)}_k=H^{(2)}_k, ~\bar{H}^{(3)}_k=H^{(3)}_k, ~\bar{H}^{(4)}_k=H^{(4)}_k,~~k\in\mathbb{N},\\
P^x_k=S^x_k, ~P^x_k+\bar{P}^x_k=T^x_k, ~P^{xy}_k=S^{xy}_k,\\
P^{xy}_k+\bar{P}^{xy}_k=T^{xy}_k,~P^y_k=S^y_k, ~P^y_k+\bar{P}^y_k=T^y_k,~k\in \bar{\mathbb{N}}.
\end{array}
\right.
\end{eqnarray*}
Thus, the results (\ref{1-Px-1})--(\ref{1-barPy-1}) of Theorem \ref{1-theorem-1} and Riccati equations (\ref{2-ST-2}) of Theorem \ref{2-theorem-1} coincide. Next, we consider taking representations $\left\{\mathbb{E}x_k, ~\mathbb{E}y_k\right\}, ~\left\{x_k-\mathbb{E}x_k, ~y_k-\mathbb{E}y_k\right\}$ into $J(\zeta^x,\zeta^y,u)$ directly to deal with Problem (MF-LQ).  Then the system state (\ref{1-system-1}) can be separated into two:
\begin{eqnarray*}
\left\{\begin{array}{l}
\mathbb{E}x_{k+1}= (A_k+\bar{A}_k)\mathbb{E}x_k+(B_k+\bar{B}_k)\mathbb{E}u_k,\\
\mathbb{E}y_{k+1}= (F_k+\bar{F}_k)\mathbb{E}y_k,\\
\mathbb{E}x_{0}= \mathbb{E}\zeta^x,~~\mathbb{E}y_{0}=\mathbb{E}\zeta^y,
\end{array}\right.
\end{eqnarray*}
and
\begin{eqnarray*}
\left\{\begin{array}{l}
x_{k+1}-\mathbb{E}x_{k+1}=\big[A_k(x_k-\mathbb{E}x_k)+B_k(u_k-\mathbb{E}u_k)\big]\\
\qquad +\big[C_k(x_k-\mathbb{E}x_k)+(C_k+\bar{C}_k)\mathbb{E}x_k\\
\qquad +D_k(u_k-\mathbb{E}u_k)+(D_k+\bar{D}_k)\mathbb{E}u_k\big]w_k,\\
y_{k+1}-\mathbb{E}y_{k+1}=F_k(y_k-\mathbb{E}y_k)\\
\qquad +\big[G_k(y_k-\mathbb{E}y_k)+(G_k+\bar{G}_k)\mathbb{E}y_k\big]v_k,\\
x_0-\mathbb{E}x_{0}=\zeta^x-\mathbb{E}\zeta^x,~y_0-\mathbb{E}y_{0}=\zeta^x-\mathbb{E}\zeta^y.
\end{array}
\right.
\end{eqnarray*}
The optimal control (\ref{1-control-2}) in Theorem \ref{1-theorem-1} is equivalent to
\begin{eqnarray*}
\begin{array}{rl}
u_k=&L^{x}_k(x_k-\mathbb{E}x_k)+(L^{x}_k+\bar{L}^{x}_k)\mathbb{E}x_k\\
&+L^y_k(y_k-\mathbb{E}y_k)+(L^{y}_k+\bar{L}^{y}_k)\mathbb{E}y_k,~k\in\mathbb{N}.
\end{array}
\end{eqnarray*}
To proceed, six sequences of symmetric matrices $\left\{S^x_k, k \in  \bar{\mathbb{N}}\right\}$, $\left\{T^x_k, k \in  \bar{\mathbb{N}}\right\}$, $\left\{S^{xy}_k, k \in  \bar{\mathbb{N}}\right\}$, $\left\{T^{xy}_k, k \in  \bar{\mathbb{N}}\right\}$, $\left\{S^y_k, k \in  \bar{\mathbb{N}}\right\}$ and $\left\{T^y_k, k \in  \bar{\mathbb{N}}\right\}$, which are (\ref{2-ST-2}), are introduced to complete the square by the classical method. Thus, the cost functional can be expressed as
\begin{eqnarray*}
\begin{array}{l}
J(\zeta^x,\zeta^y,u)\\
=\ds\sum^{N-1}_{k=0}\mathbb{E}\big[ (x_k-\mathbb{E}x_k)^TQ_k(x_k-\mathbb{E}x_k)+(\mathbb{E}x_k)^T(Q_k+\bar{Q}_k)\mathbb{E}x_k\\
\quad -2(y_k-\mathbb{E}y_k)^TQ_k(x_k-\mathbb{E}x_k)-2(\mathbb{E}y_k)^T(Q_k+\bar{Q}_k)\mathbb{E}x_k\\
\quad +(y_k-\mathbb{E}y_k)^TQ_k(y_k-\mathbb{E}y_k)+(\mathbb{E}y_k)^T(Q_k+\bar{Q}_k)\mathbb{E}y_k\\
\quad +(u_k-\mathbb{E}u_k)^TR_k(u_k-\mathbb{E}u_k)+(\mathbb{E}u_k)^T(R_k+\bar{R}_k)\mathbb{E}u_k\big]\\
\quad +\mathbb{E}\big[(x_N-\mathbb{E}x_N)^TQ_N(x_N-\mathbb{E}x_N)+(\mathbb{E}x_N)^T(Q_N+\bar{Q}_N)\mathbb{E}x_N\\
\quad -2(y_N-\mathbb{E}y_N)^TQ_N(x_N-\mathbb{E}x_N)-2(\mathbb{E}y_N)^T(Q_N+\bar{Q}_N)\mathbb{E}x_N\\
\quad +(y_N-\mathbb{E}y_N)^TQ_N(y_N-\mathbb{E}y_N)+(\mathbb{E}y_N)^T(Q_N+\bar{Q}_N)\mathbb{E}y_N\big]\\
=\ds\sum^{N-1}_{k=0}\Delta_k +\mathbb{E}\big[(\zeta^x-\mathbb{E}\zeta^x)^TS^x_0(\zeta^x-\mathbb{E}\zeta^x)+(\mathbb{E}\zeta^x)^TT^x_0\mathbb{E}\zeta^x \\
\quad +2(\zeta^y-\mathbb{E}\zeta^y)^TS^{xy}_0(\zeta^x-\mathbb{E}\zeta^x)+2(\mathbb{E}\zeta^y)^TT^{xy}_0\mathbb{E}\zeta^x\\
\quad +(\zeta^y-\mathbb{E}\zeta^y)^TS^y_0(\zeta^y-\mathbb{E}\zeta^y)+(\mathbb{E}\zeta^y)^TT^y_0\mathbb{E}\zeta^y \big],
\end{array}
\end{eqnarray*}
where
\begin{eqnarray*}
\begin{array}{l}
\Delta_k =\mathbb{E}\big[ (x_k-\mathbb{E}x_k)^TQ_k(x_k-\mathbb{E}x_k)+(x_{k+1}-\mathbb{E}x_{k+1})^T\\
\quad \cdot S^x_{k+1}(x_{k+1}-\mathbb{E}x_{k+1})-(x_k-\mathbb{E}x_k)^TS^x_k(x_k-\mathbb{E}x_k)\\
\quad -2(y_k-\mathbb{E}y_k)^TQ_k(x_k-\mathbb{E}x_k)+2(y_{k+1}-\mathbb{E}y_{k+1})^T\\
\quad \cdot S^{xy}_{k+1}(x_{k+1}-\mathbb{E}x_{k+1})-2(y_k-\mathbb{E}y_k)^TS^{xy}_k(x_k-\mathbb{E}x_k)\\
\quad +(y_k-\mathbb{E}y_k)^TQ_k(y_k-\mathbb{E}y_k)+(y_{k+1}-\mathbb{E}y_{k+1})^TS^y_{k+1}\\
\quad \cdot(y_{k+1}-\mathbb{E}y_{k+1})-(y_k-\mathbb{E}y_k)^TS^y_k(y_k-\mathbb{E}y_k)\\
\quad +(\mathbb{E}x_k)^T(Q_k+\bar{Q}_k)\mathbb{E}x_k+(\mathbb{E}x_{k+1})^TT^x_{k+1}\mathbb{E}x_{k+1}\\
\quad -2(\mathbb{E}y_k)^T(Q_k+\bar{Q}_k)\mathbb{E}x_k+2(\mathbb{E}y_{k+1})^TT^{xy}_{k+1}\mathbb{E}x_{k+1}\\
\end{array}
\end{eqnarray*}
can be reduced to
\begin{eqnarray*}
\begin{array}{l}
\Delta_k
=\mathbb{E}\Big[\Big(u_k-\mathbb{E}u_k+({W}_k^{(1)})^{-1}({H}_k^{(1)})^T(x_k-\mathbb{E}x_k)+({W}_k^{(1)})^{-1}\\
\qquad \cdot({H}_k^{(3)})^T(y-\mathbb{E}y_k)\Big)^T W^{(1)}_k \Big(u_k-\mathbb{E}u_k+({W}_k^{(1)})^{-1}({H}_k^{(1)})^T\\
\qquad \cdot(x_k-\mathbb{E}x_k)+({W}_k^{(1)})^{-1}({H}_k^{(3)})^T(y_k-\mathbb{E}y_k)\Big)+\Big(\mathbb{E}u_k\\
\qquad +({W}_k^{(2)})^{-1}({H}_k^{(2)})^T\mathbb{E}x_k+{W}_k^{(2)})^{-1}({H}_k^{(4)})^T\mathbb{E}y_k\Big)^TW^{(2)}_k\\
\qquad \cdot \Big(\mathbb{E}u_k+({W}_k^{(2)})^{-1}({H}_k^{(2)})^T\mathbb{E}x_k+({W}_k^{(2)})^{-1}({H}_k^{(4)})^T\mathbb{E}y_k \Big)\Big].
\end{array}
\end{eqnarray*}
Consequently,
\begin{eqnarray*}
\begin{array}{l}
J(\zeta^x,\zeta^y,u)=\ds\sum^{N-1}_{k=0}\mathbb{E}\Big[\Big(u_k-\mathbb{E}u_k+({W}_k^{(1)})^{-1}({H}_k^{(1)})^T\\
\quad \cdot(x_k-\mathbb{E}x_k)+({W}_k^{(1)})^{-1}({H}_k^{(3)})^T(y_k-\mathbb{E}y_k)\Big)^TW^{(1)}_k\\
\quad \cdot \Big(u_k-\mathbb{E}u_k+({W}_k^{(1)})^{-1}({H}_k^{(1)})^T(x_k-\mathbb{E}x_k)+({W}_k^{(1)})^{-1}\\
\quad \cdot({H}_k^{(3)})^T(y_k-\mathbb{E}y_k)\Big)+\Big(\mathbb{E}u_k+({W}_k^{(2)})^{-1}({H}_k^{(2)})^T\mathbb{E}x_k\\
\quad +({W}_k^{(2)})^{-1}({H}_k^{(4)})^T\mathbb{E}y_k\Big)^T W^{(2)}_k\Big(\mathbb{E}u_k+({W}_k^{(2)})^{-1}\\
\quad \cdot({H}_k^{(2)})^T\mathbb{E}x_k+({W}_k^{(2)})^{-1}({H}_k^{(4)})^T\mathbb{E}y_k\Big)\Big]\\
\quad +\mathbb{E}\big[(\zeta^x-\mathbb{E}\zeta^x)^TS^x_0(\zeta^x-\mathbb{E}x_0)+(\mathbb{E}\zeta^x)^TT^x_0\mathbb{E}\zeta^x\\
\quad +2(\zeta^y-\mathbb{E}\zeta^y)^TS^{xy}_0(\zeta^x-\mathbb{E}\zeta^x)+2(\mathbb{E}\zeta^y)^TT^{xy}_0\mathbb{E}\zeta^x\\
\quad +(\zeta^y-\mathbb{E}\zeta^y)^TS^y_0(\zeta^y-\mathbb{E}\zeta^y)+(\mathbb{E}\zeta^y)^TT^y_0\mathbb{E}\zeta^y\big].
\end{array}
\end{eqnarray*}
Let
\begin{eqnarray*}
\begin{array}{l}
\bar{u}_k=-({W}_k^{(1)})^{-1}({H}_k^{(1)})^T(x_k-\mathbb{E}x_k)-({W}_k^{(2)})^{-1}\\
\quad \cdot({H}_k^{(2)})^T\mathbb{E}x_k-({W}_k^{(1)})^{-1}({H}_k^{(3)})^T(y_k-\mathbb{E}y_k)\\
\quad -({W}_k^{(2)})^{-1}({H}_k^{(4)})^T\mathbb{E}y_k,
\end{array}
\end{eqnarray*}
for each $k\in \mathbb{N}$, which is equivalent to
\begin{eqnarray*}
\left\{\begin{array}{l}
\bar{u}_k-\mathbb{E}\bar{u}_k = -({W}_k^{(1)})^{-1}({H}_k^{(1)})^T(x_k-\mathbb{E}x_k)\\
\qquad -({W}_k^{(1)})^{-1}({H}_k^{(3)})^T(y_k-\mathbb{E}y_k),\\
\mathbb{E}\bar{u}_k=-({W}_k^{(2)})^{-1}({H}_k^{(2)})^T\mathbb{E}x_k-({W}_k^{(2)})^{-1}({H}_k^{(4)})^T\mathbb{E}y_k.
\end{array}\right.
\end{eqnarray*}
We have
\begin{eqnarray*}
\begin{array}{l}
J(\zeta^x,\zeta^y,\bar{u})
=\mathbb{E}\big[(\zeta^x-\mathbb{E}\zeta^x)^TS^x_0(\zeta^x-\mathbb{E}\zeta^x)+(\mathbb{E}\zeta^x)^TT^x_0\mathbb{E}\zeta^x \\
\quad +2(\zeta^y-\mathbb{E}\zeta^y)^TS^{xy}_0(\zeta^x-\mathbb{E}\zeta^x)+2(\mathbb{E}\zeta^y)^TT^{xy}_0\mathbb{E}\zeta^x\\
\quad +(\zeta^y-\mathbb{E}\zeta^y)^TS^{y}_0(\zeta^y-\mathbb{E}\zeta^y)+(\mathbb{E}\zeta^y)^TT^y_0\mathbb{E}\zeta^y]\\
\leq J(\zeta^x,\zeta^y,u)
\end{array}
\end{eqnarray*}
for any $u=(u_0,...,u_{N-1})$ with $u_k \in \mathcal{U}_k$. The above optimal strategy coincides with the result (\ref{2-control-1}) of Theorem \ref{2-theorem-1}.

\textbf{Remark 4.2} In fact, we can extend Theorem \ref{2-theorem-1} to the case with multidimensional noises. Specifically, consider the following system equations for $k\in\mathbb{N}$:
\begin{eqnarray}\label{3-remark-system-1}
\left\{\begin{array}{l}
x_{k+1}=(A_kx_k+\bar{A}_k\mathbb{E}x_k+B_ku_k+\bar{B}_k\mathbb{E}u_k)\\
\quad +\sum^m_{i=1}(C^i_{k}x_k+\bar{C}^i_{k}\mathbb{E}x_k+D^i_{k}u_k+\bar{D}^i_{k}\mathbb{E}u_k)w^i_{k}, \\
y_{k+1}=(F_ky_k+\bar{F}_k\mathbb{E}y_k)+ \sum^m_{i=1}(G^i_{k}y_k+\bar{G}^i_{k}\mathbb{E}y_k)v^i_{k}, \\
x_0=\zeta^x, y_0=\zeta^y.
\end{array}\right.
\end{eqnarray}
Here, $\left\{w_k=(w^1_{k},...,w^m_{k})^T,k=0,1,...,N-1\right\}$ and $\left\{v_k=(v^1_{k},...,v^m_{k})^T,k=0,1,...,N-1\right\}$ are vector-valued martingale difference sequences defined on a probability space $(\Omega,\mathcal{F},P)$ in the sense that
\begin{eqnarray}\label{3-remark-moment-1}\begin{array}{l}
\mathbb{E}\big[w_{k+1}|\mathcal{F}_k\big]=\mathbb{E}\big[v_{k+1}|\mathcal{F}_k\big]=0,~
\mathbb{E}\big[w_{k+1}w_{k+1}^T|\mathcal{F}_k\big]=\alpha_{k+1},\\
\mathbb{E}\big[v_{k+1}v_{k+1}^T|\mathcal{F}_k\big]=\beta_{k+1},~
\mathbb{E}\big[w_{k+1}v_{k+1}^T|\mathcal{F}_k\big]=\gamma_{k+1},
\end{array}
\end{eqnarray}
where $\alpha_{k+1}=(\alpha^{ij}_{k+1})_{m\times m}$, $\beta_{k+1}=(\beta^{ij}_{k+1})_{m\times m}$, $\gamma_{k+1}=(\gamma^{ij}_{k+1})_{m\times m}$ are deterministic positive semi-definite matrices. Recall that $\mathcal{F}_k$ is the $\sigma$-algebra generated by $\{\zeta^x, w_l, l=0, 1,\cdots,k\}$, $\{\zeta^y,v_l, l=0, 1,\cdots,k\}$. By taking expectations in both sides of (\ref{3-remark-system-1}), we get

\begin{eqnarray*}
\left\{\begin{array}{ll}
\mathbb{E}x_{k+1}=(A_k+\bar{A}_k)\mathbb{E}x_k+(B_k+\bar{B}_k)\mathbb{E}u_k,\\
\mathbb{E}y_{k+1}=(F_k+\bar{F}_k)\mathbb{E}y_k,\\
\mathbb{E}x_0=\mathbb{E}\zeta^x, \mathbb{E}y_0=\mathbb{E}\zeta^y.
\end{array}\right.
\end{eqnarray*}
Therefore,
\begin{eqnarray}\label{3-remark-system-2}
\left\{\begin{array}{l}
x_{k+1}-\mathbb{E}x_{k+1}=A_k(x_k-\mathbb{E}x_k)+B_k(u_k-\mathbb{E}u_k)\\
\quad +\sum^n_{i=1}\big[C^i_{k}(x_k-\mathbb{E}x_k)+(C^i_k+\bar{C}^i_k)\mathbb{E}x_k\\
\quad +D^i_{k}(u_k-\mathbb{E}u_k)+(D^i_k+\bar{D}^i_k)\mathbb{E}x_k\big]w^i_{k},\\
y_{k+1}-\mathbb{E}y_{k+1}=F_k(y_k-\mathbb{E}y_k)+ \sum^n_{i=1}G^i_{k}(y_k-\mathbb{E}y_k)v^i_{k},\\
x_0-\mathbb{E}x_0=\zeta^x-\mathbb{E}\zeta^x,~y_0-\mathbb{E}y_0=\zeta^y-\mathbb{E}\zeta^y.
\end{array}\right.
\end{eqnarray}
Similar to Theorem \ref{2-theorem-1}, we have the following results.
\begin{theorem}\label{2-theorem-2}
The Riccati equations corresponding to (\ref{3-remark-system-1}) and cost functional (\ref{1-cost-1}) are
\begin{eqnarray}\label{2-remark-ST-1}
\left\{\begin{array}{l}
S_k^x=Q_k+A_k^TS_{k+1}^xA_k+\sum^n_{i,j=1}\alpha^{ij}_{k} (C^i_{k})^TS_{k+1}^xC^j_{k}\\
\qquad -{H}_k^{(1)}({W}_k^{(1)})^{-1}({H}_k^{(1)})^T,\\
T_k^x=Q_k+\bar{Q}_k+(A_k+\bar{A}_k)^TT_{k+1}^x(A_k+\bar{A}_k)\\
\qquad +\sum^n_{i,j=1}\alpha^{ij}_{k} (C^i_{k}+\bar{C}^i_{k})^TS_{k+1}^x(C^j_{k}+\bar{C}^j_{k})\\
\qquad -{H}_k^{(2)}({W}_k^{(2)})^{-1}({H}_k^{(2)})^T,\\
S_k^{xy}=-Q_k+F_k^TS_{k+1}^{xy} A_k+\sum^n_{i,j=1}\gamma^{ij}_{k} (G^i_{k})^TS_{k+1}^{xy}C^j_{k}\\
\qquad -{H}_k^{(3)}({W}_k^{(1)})^{-1}({H}_k^{(1)})^T,\\
T_k^{xy}=-Q_k-\bar{Q}_k+(F_k+\bar{F}_k)^TT_{k+1}^{xy}(A_k+\bar{A}_k)\\
\qquad +\sum^n_{i,j=1}\gamma^{ij}_{k} (G^i_{k}+\bar{G}^i_{k})^TS_{k+1}^{xy}(C^j_{k}+\bar{C}^j_{k})\\
\qquad -{H}_k^{(4)}({W}_k^{(2)})^{-1}({H}_k^{(2)})^T,\\
S_k^y=Q_k+F_k^TS_{k+1}^yF_k+\sum^n_{i,j=1}\beta^{ij}_{k}(G^i_{k})^TS_{k+1}^yG^j_{k}\\
\qquad -{H}_k^{(3)}({W}_k^{(1)})^{-1}({H}_k^{(3)})^T,\\
T_k^y=Q_k+\bar{Q}_k+(F_k+\bar{F}_k)^TT_{k+1}^y(F_k+\bar{F}_k)\\
\qquad +\sum^n_{i,j=1}\beta^{ij}_{k} (G^i_{k}+\bar{G}^i_{k})^TS_{k+1}^y(G^j_{k}+\bar{G}^j_{k})\\
\qquad -{H}_k^{(4)}({W}_k^{(2)})^{-1}({H}_k^{(4)})^T,\\
S_N^x=Q_N,~T_N^x=Q_N+\bar{Q}_N,~S_N^{xy}=-Q_N,\\
T_N^{xy}=-Q_N-\bar{Q}_N,~S_N^y=Q_N,~T_N^y=Q_N+\bar{Q}_N,
\end{array}\right.
\end{eqnarray}
with
\begin{eqnarray}\label{2-remark-WH-1}
\left\{\begin{array}{l}
W_k^{(1)}=R_k+B_k^TS^x_{k+1}B_k+\sum^n_{i,j=1}\alpha^{ij}_{k}(D^i_{k})^TS^x_{k+1}D^j_{k},\\
W_k^{(2)}=R_k+\bar{R}_k+(B_k+\bar{B}_k)^TT^x_{k+1}(B_k+\bar{B}_k)\\
\qquad +\sum^n_{i,j=1}\alpha^{ij}_{k} (D^i_{k}+\bar{D}^i_{k})^TS^x_{k+1}(D^j_{k}+\bar{D}^j_{k}),\\
H_k^{(1)}=A_k^TS^x_{k+1}B_k+\sum^n_{i,j=1}\alpha^{ij}_{k}(C^i_{k})^TS^x_{k+1}D^j_{k},\\
H_k^{(2)}=(A_k+\bar{A}_k)^TT^x_{k+1}(B_k+\bar{B}_k)\\
\qquad +\sum^n_{i,j=1}\alpha^{ij}_{k} (C^i_{k}+\bar{C}^i_{k})^TS^x_{k+1}(D^j_{k}+\bar{D}^j_{k}),\\
H_k^{(3)}= F_k^TS^{xy}_{k+1}B_k+\sum^n_{i,j=1}\gamma^{ij}_{k}(G^i_{k})^TS^{xy}_{k+1}D^j_{k},\\
H_k^{(4)}=(F_k+\bar{F}_k)^TT^{xy}_{k+1}(B_k+\bar{B}_k)\\
\qquad +\sum^n_{i,j=1}\gamma^{ij}_{k} (G^i_{k}+\bar{G}^i_{k})^TS^{xy}_{k+1}(D^j_{k}+\bar{D}^j_{k}).
\end{array}\right.
\end{eqnarray}
Under conditions
\begin{eqnarray}\label{2-remark-condition-2}
Q_k, Q_k+\bar{Q}_k\geq 0,~k\in\bar{\mathbb{N}},~R_k, R_k+\bar{R}_k>0,~k\in\mathbb{N},
\end{eqnarray}
the unique optimal control for Problem (MF-LQ) is
\begin{eqnarray}\label{3-remark-control-2}
\begin{array}{l}
u_k^*=-({W}_k^{(1)})^{-1}({H}_k^{(1)})^T(x_k-\mathbb{E}x_k)-({W}_k^{(2)})^{-1}\\
\quad \cdot({H}_k^{(2)})^T\mathbb{E}x_k-({W}_k^{(1)})^{-1}({H}_k^{(3)})^T(y_k-\mathbb{E}y_k)\\
\quad -({W}_k^{(2)})^{-1}({H}_k^{(4)})^T\mathbb{E}y_k,~~k\in\mathbb{N}.
\end{array}
\end{eqnarray}
\end{theorem}
\emph{Proof}. The proof of above is a straightforward extension of that of Theorem \ref{2-theorem-1} and hence omitted here.
\hfill $\square$

\section{An Example\label{5}}
Basing upon the general theory in previous sections, in this section, we consider an example extended from financial application in asset-liability management with numerical results.
\subsection{Example Setting\label{6}}

Suppose a investment market and a loan market consisting of $m$ risky investment acceptable assets, one risk-free asset and one loan product within a time horizon $N$. Let $B_k=(B^1_k,...,B^m_k)$ be the row vector of random excess returns of the $m$ risky assets, $a_k$ and $f_k$ are given return of the risk-free asset and repayment of loan at time period $k$. We assume that vectors $B_k,~k=0,1,...,N-1$ are statistically independent and the only information known about the random excess return vector $B_k$ is its first two moments: its mean $\mathbb{E}(B_k)$ and covariance $Cov(B_k)$.

Let $x_k$ and $y_k$ be the total asset and liability at the beginning of the $k$-th period, respectively. Let $u^i_k,i=1,2,...,m,$ be the amount invested in the $i$-th risky asset at period $k$. The system combined with assets and liabilities at the beginning of the $(k+1)$-th period is given by
\begin{eqnarray}\label{3-system-1}
\left\{\begin{array}{l}
x_{k+1}=a_kx_k+B_ku_k,\\
y_{k+1}=f_ky_k,\\
x_0=\zeta^x,~y_0=\zeta^y,
\end{array}\right.
\end{eqnarray}
We shall transform (\ref{3-system-1}) into the form of (\ref{3-remark-system-1}), by which the general theory in above sections works. Define
\begin{eqnarray*}
\left\{\begin{array}{l}
D^i_k=(0,...,0,1,0,...,0),\mbox{where~1~is~the~$i$th~entry},\\
w^i_k=B^i_k-\mathbb{E}(B^i_k),~w_k=(w^1_k,w^2_k,...,w^m_k)^T,\\
i=1,...,m,~~k=1,...,N-1.
\end{array}\right.
\end{eqnarray*}
These lead to
\begin{eqnarray}\label{3-system-2}
\left\{\begin{array}{l}
x_{k+1}=a_kx_k+\mathbb{E}B_ku_k+\sum^m_{i=1}D^i_ku_kw^i_k,\\
y_{k+1}=f_ky_k.\\
x_0=\zeta^x,~y_0=\zeta^y.
\end{array}\right.
\end{eqnarray}
Clearly, $x_k, y_k\in \mathbb{R}$, $k \in \mathbb{N}$. $a_k, f_k\in \mathbb{R}$ and $\mathbb{E}B_k, D^i_k\in \mathbb{R}^{1\times m}$ are deterministic. By taking expectation of the system state, we have
\begin{eqnarray*}
\left\{\begin{array}{l}
\mathbb{E}x_{k+1}=a_k\mathbb{E}x_k+\mathbb{E}B_k\mathbb{E}u_k,\\
\mathbb{E}y_{k+1}=f_k\mathbb{E}y_k,\\
\mathbb{E}x_0=\mathbb{E}\zeta^x,~\mathbb{E}y_0=\mathbb{E}\zeta^y.
\end{array}\right.
\end{eqnarray*}
Hence,
\begin{eqnarray}\label{3-system-3}
\left\{\begin{array}{l}
x_{k+1}-\mathbb{E}x_{k+1}=a_k(x_k-\mathbb{E}x_k)+\mathbb{E}B_k(u_k-\mathbb{E}u_k)\\
\quad +\sum^m_{i=1}D^i_ku_kw^i_k,\\
y_{k+1}-\mathbb{E}y_{k+1}=f_k(y_k-\mathbb{E}y_k).\\
x_0=\zeta^x,~y_0=\zeta^y.\\
\end{array}\right.
\end{eqnarray}
Define $\mathcal{F}_k'$ by the information set at the beginning of period $k$ which is generated by $\left\{\zeta^x,w_l,l=0,1,...,k\right\}$. Recall that $w_k$ is a martingale difference sequence defined on a probability space $(\Omega,\mathcal{F},P)$, where $\mathbb{E}\big[w_{k+1}w_{k+1}^T|\mathcal{F}_k'\big]=\alpha_{k+1}=Cov(B_{k+1})$. The cost functional (an extension of variance function) associated with (\ref{3-system-1}) is
\begin{eqnarray}\label{3-cost-1}
\begin{array}{l}
J(\zeta^x,\zeta^y,u) = \ds\sum_{k=0}^{N-1}\mathbb{E}\big(u_k^TR_ku_k\big)+\mathbb{E}\big(q_N(x_N-y_N)^2\big)\\
\qquad\qquad +\bar{q}_N\big(\mathbb{E}(x_N-y_N)\big)^2,
\end{array}
\end{eqnarray}
where $q_N, \bar{q}_N, R_k, \bar{R}_k, k\in\mathbb{N}$ are deterministic symmetric matrices with appropriate dimensions. In this paper, we consider the case where short-selling of stock is allowed, i.e., $u^i_k,i=1,...,k$, could take values in $\mathbb{R}$. Hence, the admissible policy set of $u=(u_0,u_1,\cdots,u_{N-1})$ in this section
\begin{eqnarray*}
\begin{array}{c}
\mathcal{U}_{ad}\equiv\left\{u~\left\vert~u_k \in \mathbb{R}^m, \mbox{\ is } \mathcal{F}_k' \mbox{-measurable }, \mathbb{E}|u_k|^2<\infty\right.\right\}.
\end{array}
\end{eqnarray*}
Problem (MF-LQ) extended from asset-liability management is represented as follows:

\noindent
\textbf{Problem (MF-AL)}.
\emph{For any given square-integrable initial values $\zeta^x$ and $\zeta^y$, find  $u^o\in \mathcal{U}_{ad}$ such that}
\begin{eqnarray} \label{3-problem-1}
J(\zeta^x,\zeta^y,u^o)=\inf_{u\in \mathcal{U}_{ad}}J(\zeta^x,\zeta^y,u).
\end{eqnarray}
\emph{We then call $u^o$ an optimal control for Problem (MF-AL).}

To proceed, we recall the following lemma \cite{Dunne1993}.
\begin{lemma}
Let $M \in \mathbb{R}^{n\times n}$, $c \in \mathbb{R}^n$. If  $c \in Range(M)$, then
\begin{eqnarray*}
\begin{array}{c}
\ds (M\pm cc^T)^{\dagger}=M-\frac {M^{\dagger}cc^TM^{\dagger}}{c^TM^{\dagger}c}. 
\end{array}
\end{eqnarray*}
\end{lemma}
The optimal strategy can be obtained by the following theorem.
\begin{theorem}\label{3-theorem-1}
Suppose that
\begin{eqnarray}\label{2-remark-condition-2}
R_k>0,~k\in\mathbb{N},~~q_N, q_N+\bar{q}_N\geq 0.
\end{eqnarray}
The unique optimal strategy for Problem (MF-AL) is given by
\begin{eqnarray}\label{3-control-1}
\begin{array}{l}
u_k^o=-({W}_k^{(1)})^{-1}({H}_k^{(1)})^T(x_k-\mathbb{E}x_k)-({W}_k^{(2)})^{-1}\\
\quad \cdot({H}_k^{(2)})^T\mathbb{E}x_k-({W}_k^{(1)})^{-1}({H}_k^{(3)})^T(y_k-\mathbb{E}y_k)\\
\quad -({W}_k^{(2)})^{-1}({H}_k^{(4)})^T\mathbb{E}y_k,~k\in\mathbb{N}.
\end{array}
\end{eqnarray}
where
\begin{eqnarray}\label{3-WH-1}
\left\{
\begin{array}{l}
W_k^{(1)}=R_k+S^x_{k+1}\mathbb{E}(B_k^TB_k),\\
W_k^{(2)}=R_k+T^x_{k+1}\mathbb{E}B_k^T\mathbb{E}B_k+S^x_{k+1}Cov(r_k),\\
H_k^{(1)}=a_kS^x_{k+1}\mathbb{E}B_k,\\
H_k^{(2)}=a_kT^x_{k+1}\mathbb{E}B_k,\\
H_k^{(3)}= f_kS^{xy}_{k+1}\mathbb{E}B_k,\\
H_k^{(4)}= f_kT^{xy}_{k+1}\mathbb{E}B_k,\\
\end{array}
\right.
\end{eqnarray}
and
\begin{eqnarray}\label{3-ST-1}
\left\{
\begin{array}{l}
S_k^x=a_k^2S_{k+1}^x\big[1-S_{k+1}^x\mathbb{E}B_k\big(R_k+S_{k+1}^x\mathbb{E}(B_k^TB_k)\big)^{-1}\mathbb{E}B_k^T\big],\\
T_k^x=a_k^2T_{k+1}^x\big[1-T_{k+1}^x\mathbb{E}B_k\big(R_k+T^x_{k+1}\mathbb{E}B_k^T\mathbb{E}B_k\\
\qquad +S_{k+1}^xCov(B_k)\big)^{-1}\mathbb{E}B_k^T\big],\\
S_k^{xy}=a_kf_kS_{k+1}^{xy}\big[1-S_{k+1}^x\mathbb{E}B_k\big(R_k+S_{k+1}^x\mathbb{E}(B_k^TB_k)\big)^{-1}\mathbb{E}B_k^T\big],\\
T_k^{xy}=a_kf_kT_{k+1}^{xy}\big[1-T_{k+1}^x\mathbb{E}B_k\big(R_k+T^x_{k+1}\mathbb{E}B_k^T\mathbb{E}B_k\\
\qquad +S_{k+1}^xCov(B_k)\big)^{-1}\mathbb{E}B_k^T\big],\\
S_k^y=f_k^2\big[S_{k+1}^y-(S_{k+1}^{xy})^2\mathbb{E}B_k\big(R_k+S_{k+1}^x\mathbb{E}(B_k^TB_k)\big)^{-1}\mathbb{E}B_k^T\big],\\
T_k^y=f_k^2\big[T_{k+1}^y-(T_{k+1}^{xy})^2\mathbb{E}B_k\big(R_k+T^x_{k+1}\mathbb{E}B_k^T\mathbb{E}B_k\\
\qquad +S_{k+1}^xCov(B_k)\big)^{-1}\mathbb{E}B_k^T\big],\\
S_N^x=q_N,~T_N^x=q_N+\bar{q}_N,~S_N^{xy}=-q_N,\\
T_N^{xy}=-q_N-\bar{q}_N,~S_N^y=q_N,~T_N^y=q_N+\bar{q}_N.
\end{array}
\right.
\end{eqnarray}
Under the optimal strategy (\ref{3-control-1}), the optimal solution of cost functional is
\begin{eqnarray}\label{3-efficient-frontier}
\begin{array}{l}
J(\zeta^x, \zeta^y, u^o)\\
=\mathbb{E}\big[S^x_0(\zeta^x-\mathbb{E}\zeta^x)^2+T^x_0(\mathbb{E}\zeta^x)^2+2S^{xy}_0(\zeta^x-\mathbb{E}\zeta^x)\\
\cdot(\zeta^y-\mathbb{E}\zeta^y)+T^{xy}_0\mathbb{E}\zeta^x\mathbb{E}\zeta^y +S^{y}_0(\zeta^y-\mathbb{E}\zeta^y)^2+T^y_0(\mathbb{E}\zeta^y)^2\big],
\end{array}
\end{eqnarray}
and its related expectation of system state in $N$-th period is
\begin{eqnarray}\label{3-expected-equity}
\begin{array}{l}
\mathbb{E}(x_N-y_N)\\
= \prod\limits_{k=1}^{N-1} N_k\mathbb{E}\zeta^x+\Big(\sum^{N-1}_{k=0}\prod\limits_{j=k+1}^{N-1}N_j M_k\prod\limits_{j=1}^{k-1}f_k-\prod\limits_{k=1}^{N-1} f_k\Big)\mathbb{E}\zeta^y
\end{array}
\end{eqnarray}
with
\begin{eqnarray*}
\left\{
\begin{array}{l}
N_k=a_k\big(1-T^x_{k+1}\mathbb{E}B_k(W^{(2)}_k)^{-1}\mathbb{E}B_k^T\big),\\
M_k=-f_kT^{xy}_{k+1}\mathbb{E}B_k(W^{(2)}_k)^{-1}\mathbb{E}B_k^T.
\end{array}\right.
\end{eqnarray*}
\end{theorem}
\noindent\textbf{Proof}. Due to the Theorem \ref{2-theorem-2}, we have the following Riccati equations:
\begin{eqnarray}\label{3-ST-2}
\left\{\begin{array}{l}
S_k^x=a_k^2S_{k+1}^x-{H}_k^{(1)}({W}_k^{(1)})^{-1}({H}_k^{(1)})^T,\\
T_k^x=a_k^2T_{k+1}^x-{H}_k^{(2)}({W}_k^{(2)})^{-1}({H}_k^{(2)})^T,\\
S_k^{xy}=a_kf_kS_{k+1}^{xy}-{H}_k^{(3)}({W}_k^{(1)})^{-1}({H}_k^{(1)})^T,\\
T_k^{xy}=a_kf_kT_{k+1}^{xy}-{H}_k^{(4)}({W}_k^{(2)})^{-1}({H}_k^{(2)})^T,\\
S_k^y=f_k^2S_{k+1}^y-{H}_k^{(3)}({W}_k^{(1)})^{-1}({H}_k^{(3)})^T,\\
T_k^y=f_k^2T_{k+1}^y-{H}_k^{(4)}({W}_k^{(2)})^{-1}({H}_k^{(4)})^T,\\
S_N^x=q_N,~~T_N^x=q_N+\bar{q}_N,~~S_N^{xy}=-q_N,\\
T_N^{xy}=-q_N-\bar{q}_N,~~S_N^y=q_N,~~T_N^y=q_N+\bar{q}_N.
\end{array}\right.
\end{eqnarray}
Here,
\begin{eqnarray*}
\left\{
\begin{array}{l}
W_k^{(1)}=R_k+\mathbb{E}B_k^TS^x_{k+1}\mathbb{E}B_k+\ds\ds\ds\ds\sum^n_{i=1}\ds\ds\ds\ds\sum^m_{j=1}\alpha_k^{ij}(D^i_k)^TS^x_{k+1}D^j_k\\
\qq =R_k+S^x_{k+1}\mathbb{E}(B_k^TB_k),\\
W_k^{(2)}=R_k+\mathbb{E}B_k^T T^x_{k+1}\mathbb{E}B_k+\ds\ds\ds\ds\sum^m_{i=1}\ds\ds\ds\ds\sum^m_{j=1}\alpha_k^{ij}(D^i_k)^TS^x_{k+1}D^j_k\\
\qq =R_k+T^x_{k+1}\mathbb{E}B_k^T\mathbb{E}B_k+S^x_{k+1}Cov(B_k),\\
H_k^{(1)}=a_kS^x_{k+1}\mathbb{E}B_k,\\
H_k^{(2)}=a_kT^x_{k+1}\mathbb{E}B_k,\\
H_k^{(3)}=f_kS^{xy}_{k+1}\mathbb{E}B_k,\\
H_k^{(4)}=f_kT^{xy}_{k+1}\mathbb{E}B_k.
\end{array}
\right.
\end{eqnarray*}
Then, we should derive $S^x_k,~T^x_k\geq 0$ for $~k\in \bar{\mathbb{N}}$. Under condition (\ref{2-remark-condition-2}), we have  $R_k > 0$ and $S^x_N,~T^x_N \geq 0$. Therefore, based on Lemma 5.1,
\begin{eqnarray*}
\begin{array}{l}
S_{N-1}^x=a_{N-1}^2S_{N}^x\big[1-S_{N}^x\mathbb{E}B_{N-1}\big(R_{N-1}\\
\qquad +S_{N}^x\mathbb{E}(B_{N-1}^TB_{N-1})\big)^{-1}\mathbb{E}B_{N-1}^T\big]\\
=a_{N-1}^2S_{N}^x\big[1-S^x_{N}\mathbb{E}B_{N-1}\big(R_{N-1}+S^x_{N}Cov(B_{N-1})\\
\qquad +S^x_{N}\mathbb{E}B_{N-1}^T\mathbb{E}B_{N-1})\big)^{-1}\mathbb{E}B_{N-1}^T\big]\\
\ds =\frac {1}{1+S^x_{N}\mathbb{E}B_{N-1}\big(R_{N-1}+S^x_{N}Cov(B_{N-1})\big)^{-1}\mathbb{E}B_{N-1}^T}\geq 0,
\end{array}
\end{eqnarray*}
and
\begin{eqnarray*}
\begin{array}{l}
T_{N-1}^x=a_{N-1}^2T_{N}^x\big[1-T_{N}^x\mathbb{E}B_{N-1}\big(R_{N-1}\\
\qquad +T^x_{N}\mathbb{E}B_{N-1}^T\mathbb{E}B_{N-1}+S_{N}^xCov(B_{N-1})\big)^{-1}\mathbb{E}B_{N-1}^T\big]\\
=a_{N-1}^2T_{N}^x\big[1-T_{N}^x\mathbb{E}B_{N-1}\big(R_{N-1}+S_{N}^xCov(B_{N-1})\\
\qquad +T^x_{N}\mathbb{E}B_{N-1}^T\mathbb{E}B_{N-1}\big)^{-1}\mathbb{E}B_{N-1}^T\big]\\
\ds =\frac {1}{1+T^x_{N}\mathbb{E}B_{N-1}\big(R_{N-1}+S_{N}^xCov(B_{N-1})\big)^{-1}\mathbb{E}B_{N-1}^T} \geq 0.
\end{array}
\end{eqnarray*}
By induction, $S^x_{k},~T^x_{k} \geq 0$ are satisfied. Finally, Based on the content of subsection 4.3, we have 
\begin{eqnarray*}
\begin{array}{l}
J(\zeta^x,\zeta^y,u)\\
=\ds\sum^{N-1}_{k=0}\mathbb{E}\Big[\Big(u_k-\mathbb{E}u_k+({W}_k^{(1)})^{-1}({H}_k^{(1)})^T(x_k-\mathbb{E}x_k)\\
+({W}_k^{(1)})^{-1}({H}_k^{(3)})^T(y_k-\mathbb{E}y_k)\Big)^TW^{(1)}_k\Big(u_k-\mathbb{E}u_k+({W}_k^{(1)})^{-1}\\
\cdot({H}_k^{(1)})^T(x_k-\mathbb{E}x_k)+({W}_k^{(1)})^{-1}({H}_k^{(3)})^T(y_k-\mathbb{E}y_k)\Big)\\
+\Big(\mathbb{E}u_k+({W}_k^{(2)})^{-1}({H}_k^{(2)})^T\mathbb{E}x_k+({W}_k^{(2)})^{-1}({H}_k^{(4)})^T\mathbb{E}y_k\Big)^T\\
\cdot W^{(2)}_k\Big(\mathbb{E}u_k+({W}_k^{(2)})^{-1}({H}_k^{(2)})^T\mathbb{E}x_k+({W}_k^{(2)})^{-1}({H}_k^{(4)})^T\mathbb{E}y_k\Big)\Big]\\
+\mathbb{E}\big[S^x_0(\zeta^x-\mathbb{E}\zeta^x)^2+T^x_0(\mathbb{E}\zeta^x)^2+2S^{xy}_0(\zeta^x-\mathbb{E}\zeta^x)\\
\cdot(\zeta^y-\mathbb{E}\zeta^y)+2T^{xy}_0\mathbb{E}\zeta^x\mathbb{E}\zeta^y +S^y_0(\zeta^y-\mathbb{E}\zeta^y)^2+T^y_0(\mathbb{E}\zeta^y)^2\big].
\end{array}
\end{eqnarray*}
Therefore, under the optimal strategy (\ref{3-control-1}), $J(\zeta^x, \zeta^y, u)$ reach to the minimum. By simple calculation, we obtain (\ref{3-expected-equity}) from Remark 4.1. This completes the proof.
\hfill $\square$\\

\subsection{Numerical Results\label{7}}
Consider a 3-period numerical example. Coefficients are given as follows:\\
$a_k=0.5$, $f_k=0.6$, $\mathbb{E}B_k=\big(0.2,~0.3,~0.4\big)$, $R_k=I$, $\bar{R}_k=0$, $q_3=1$, $\bar{q}_3=-1$, 
$Cov(r_k)=\begin{pmatrix} 
1 & ~0.2& ~0.3\\[-2mm]
0.2&~1 & ~0.6\\[-2mm]
0.3 & ~0.6 &~1
\end{pmatrix}.
$\\
By simple calculation, we have\\
$\mathbb{E}B_k^T\mathbb{E}B_k=\begin{pmatrix} 
0.040 & ~~0.060 & ~~0.080\\[-2mm]
0.060 & ~~0.090 & ~~0.120\\[-2mm]
0.080 & ~~0.120 & ~~0.160
\end{pmatrix}
$\\
$\mathbb{E}(B_k^TB_k)=\begin{pmatrix} 
1.040 & ~~0.260 & ~~0.380\\[-2mm]
0.260 & ~~1.090 & ~~0.720\\[-2mm]
0.380 & ~~0.720 & ~~1.160
\end{pmatrix}.
$\\

Based on Theorem \ref{3-theorem-1}, the Riccati solutions for $S_k$ and $T_k$ for $k=0,1,2,3$ are given by
\begin{eqnarray*}
\begin{array}{rcl}
S^x_3=1,~~~~~~ &~~S^{xy}_3=-1,~~~~~ &~~S^y_3=1,\\[-1mm]
S^x_2=0.2260, &~~S^{xy}_2=-0.2712, &~~S^y_2=0.3254,\\[-1mm]
S^x_1=0.0540, &~~S^{xy}_1=-0.0777, &~~S^y_1=0.1119,\\[-1mm]
S^x_0=0.0133, &~~S^{xy}_0=-0.0230, &~~S^y_0=0.0397
\end{array}
\end{eqnarray*}
and $T^x_k=T^{xy}_k=T^y_k=0$, which lead $N_k=a_k$ and $M_k=0$. We also obtain the optimal control, that is $u^o_k=O^{x}_k(x_k-\mathbb{E}x_k)+O^{y}_k(y_k-\mathbb{E}y_k),~k=0,1,2$, where
\begin{eqnarray*}
\begin{array}{l}
O^{x}_2=\big(-0.0300~-0.0429~-0.0730\big),\\ [-1mm]
O^{y}_2=\big(~~0.0359~~~0.0515~~~0.0876\big),\\[-1mm]
O^{x}_1=\big(-0.0150~-0.0223~-0.0319\big),\\[-1mm]
O^{y}_1=\big(~~0.0216~~~0.0321~~~0.0460\big),\\[-1mm]
O^{x}_0=\big(-0.0048~-0.0072~-0.0098\big),\\[-1mm]
O^{y}_0=\big(~~0.0069~~~0.0104~~~0.0141\big).
\end{array}
\end{eqnarray*}

\section{Conclusion}

In this paper, we formulate the Problem (MF-LQ) and give necessary and sufficient conditions for the solvability of problem.
Two approaches, dynamical optimization by matrix minimum principle and operator linear-quadratic method, are investigated to derived the optimal control, where six Riccati equations is obtained concomitantly.
The solution of these two approaches are derived to be coincided through method of completing the square.
Also, after concerning with the solution of this problem under multidimensional noise assumption, we give an financial application with numerical results.
For future research, we may study a model with relaxation of conditions such as indefinite mean-field stochastic linear-quadratic optimal control problem and may also expend the problem from finite horizon to infinite, where the stability of system should be considered first.

\section*{Appendix}

\subsection*{Proof of Theorem \ref{1-theorem-1}.}

\emph{Proof}. From Proposition \ref{1-proposition-1}, we know that Problem (MF-LQ) admits a unique minimizer under condition (\ref{1-condition-2}). Thus, the optimal control uniquely exists. We now introduce the Lagrangian function associated with
Problem (\ref{1-problem-2}),
\begin{eqnarray}\label{1-lagrangian-1}
\begin{array}{l}
\mathfrak{L} = \ds\sum_{k=0}^{N-1} \mathfrak{L}_k+Tr\big[Q_N\big(X_N-XY_N-(XY_N)^T+Y_N\big)\big]\\
\quad +Tr\big[ \bar{Q}_N \big(\bar{X}_N-\bar{XY}_N-(\bar{XY}_N)^T+\bar{Y}_N\big)\big]\\
= \ds\sum_{k=0}^{N-1} \mathfrak{L}_k+Tr\left[\left(\begin{array}{c}
Q_N\\[-2mm]
\bar{Q}_N\\[-2mm]
-2Q_N\\[-2mm]
-2\bar{Q}_N\\[-2mm]
Q_N\\[-2mm]
\bar{Q}_N
\end{array}\right)^T\centerdot
\left(\begin{array}{c}X_N\\[-2mm] 
\bar{X}_N\\[-2mm]
XY_N\\[-2mm]
\bar{XY}_N\\[-2mm]
Y_N\\[-2mm]
\bar{Y}_N
\end{array}\right)\right],
\end{array}
\end{eqnarray}
where
\begin{eqnarray}\label{1-lagrangian-2}
\begin{array}{l}
\mathfrak{L}_k=Tr\big[Q_k\big(X_k-XY_k-(XY_k)^T+Y_k\big)\big]\\
\quad +Tr\big[ \bar{Q}_k \big(\bar{X}_k-\bar{XY}_k-(\bar{XY}_k)^T+\bar{Y}_k\big)\big]\\
\quad +Tr\big[(L^x_k)^TR_kL^x_k X_k +\big((L^x_k)^TR_k\bar{L}^x_k+(\bar{L}^x_k)^TR_kL^x_k\\
\quad +(\bar{L}^x_k)^TR_k\bar{L}^x_k\big)\bar{X}_k \big]+Tr\big[(L^x_k)^TR_kL^y_k (XY_k)^T\\
\quad +\big((L^x_k)^TR_k\bar{L}^y_k+(\bar{L}^x_k)^TR_kL^y_k+(\bar{L}^x_k)^TR_k\bar{L}^y_k\big)(\bar{XY}_k)^T \big]\\
\quad +Tr\big[(L^y_k)^TR_kL^x_k XY_k +\big((L^y_k)^TR_k\bar{L}^x_k+(\bar{L}^y_k)^TR_kL^x_k\\
\quad +(\bar{L}^y_k)^TR_k\bar{L}^x_k\big)\bar{XY}_k\big]+Tr\big[(L^y_k)^TR_kL^y_k Y_k\\
\quad +\big((L^y_k)^TR_k\bar{L}^y_k+(\bar{L}^y_k)^TR_kL^y_k+(\bar{L}^y_k)^TR_k\bar{L}^y_k\big)\bar{Y}_k \big]\\
\quad +Tr\big[(L^x_k+\bar{L}^x_k)^T\bar{R}_k(L^x_k+\bar{L}^x_k\big)\bar{X}_k+ (L^x_k+\bar{L}^x_k)^T\bar{R}_k\\
\quad \cdot(L^y_k+\bar{L}^y_k\big)(\bar{XY}_k)^T \big]+Tr\big[(L^y_k+\bar{L}^y_k)^T\bar{R}_k(L^x_k+\bar{L}^x_k\big)\\
\quad \cdot\bar{XY}_k+ (L^y_k+\bar{L}^y_k)^T\bar{R}_k(L^y_k+\bar{L}^y_k\big)\bar{Y}_k \big]\\
+Tr\left[{\left(\begin{array}{c}
P^x_{k+1}\\[-2mm]
\bar{P}^x_{k+1}\\[-2mm]
~2P^{xy}_{k+1}\\[-2mm]
2\bar{P}^{xy}_{k+1}\\[-2mm]
P^{y}_{k+1}\\[-2mm]
\bar{P}^{y}_{k+1}
\end{array}\right)}^T\centerdot
\left(\begin{array}{c}
\mathcal{X}_k(L^x_k, \bar{L}^x_k, L^y_k, \bar{L}^y_k)-X_{k+1}\\[-2mm] 
\bar{\mathcal{X}}_k(L^x_k, \bar{L}^x_k, L^y_k, \bar{L}^y_k)-\bar{X}_{k+1}\\[-2mm]
\mathcal{XY}_k(L^x_k, \bar{L}^x_k, L^y_k, \bar{L}^y_k)-XY_{k+1}\\[-2mm]
\bar{\mathcal{XY}}_k(L^x_k, \bar{L}^x_k, L^y_k, \bar{L}^y_k)-\bar{XY}_{k+1}\\[-2mm]
\mathcal{Y}_k-Y_{k+1}\\[-2mm]
\bar{\mathcal{Y}}_k-\bar{Y}_{k+1}\\[-2mm]
\end{array}\right)\right],
\end{array}
\end{eqnarray}
and $P^x_{k+1}~\bar{P}^x_{k+1}~P^{xy}_{k+1}~\bar{P}^{xy}_{k+1}~P^y_{k+1}~\bar{P}^y_{k+1},~k\in\mathbb{N}$ are the Lagrangian
multipliers.
Denote $\mathbb{{P}}_{k+1}=\big(P^x_{k+1}~\bar{P}^x_{k+1}~2P^{xy}_{k+1}~2\bar{P}^{xy}_{k+1}~P^{y}_{k+1}~\bar{P}^{y}_{k+1}\big)$ and 
$\mathbb{X}_k=\big(X_k~\bar{X}_k~XY_k~\bar{XY}_k~Y_k~\bar{Y}_k\big)$. Clearly, by the matrix minimum principle \cite{Athans1967}, the optimal feedback gains $L^o_k$ and Lagrangian multipliers
$\mathbb{P}_{k}$ satisfy the following first-order
necessary conditions
\begin{eqnarray*}
\left\{
\begin{array}{l}
\ds\frac{\partial\mathfrak{L}_k}{\partial L^o_k}=0,~~\mathbb{P}_k=\frac{\partial\mathfrak{L}_k}{\partial \mathbb{X}_k},~~k\in\mathbb{N},\\
\mathbb{P}_N=\big(Q_N~~\bar{Q}_N~-2Q_N ~-2\bar{Q}~~Q_N ~~\bar{Q}_N\big),
\end{array}\right.
\end{eqnarray*}
i.e.,
\begin{eqnarray}\label{1-P-1}
\left\{\begin{array}{l}
\ds\frac{\partial\mathfrak{L}_k}{\partial L^{xo}_k}=0,~\frac{\partial\mathfrak{L}_k}{\partial \bar{L}^{xo}_k}=0,~\frac{\partial\mathfrak{L}_k}{\partial L^{yo}_k}=0,~\frac{\partial\mathfrak{L}_k}{\partial \bar{L}^{yo}_k}=0,\\
\ds {P}^x_k=\frac{\partial\mathfrak{L}_k}{\partial {X}_k},~\bar{P}^x_k=\frac{\partial\mathfrak{L}_k}{\partial \bar{{X}}_k},~P^{xy}_k=\frac{1}{2}\frac{\partial\mathfrak{L}_k}{\partial {XY}_k},\\
\ds\bar{P}^{xy}_k=\frac{1}{2}\frac{\partial\mathfrak{L}_k}{\partial \bar{{XY}}_k},~{P}^y_k=\frac{\partial\mathfrak{L}_k}{\partial {Y}_k},~\bar{P}^y_k=\frac{\partial\mathfrak{L}_k}{\partial \bar{{Y}}_k},~k\in\mathbb{N}\\
P^x_N=Q_N,~\bar{P}^x_N=\bar{Q}_N,~P^{xy}_N=-Q_N,\\
\bar{P}^{xy}_N=-\bar{Q}_N,~P^y_N=Q_N,~\bar{P}^y_N=\bar{Q}_N.\\
\end{array}\right.
\end{eqnarray}
Now, we shall calculate several gradient matrices. Firstly, we have
\begin{eqnarray}\label{1-Lx-1}
\begin{array}{l}
\ds\frac{\partial \mathfrak{L}_k}{\partial L^x_k}=2X_k\big[(L^x_k)^T(R_k+B_k^TP^x_{k+1}B_k+D_k^TP^x_{k+1}D_k)\\
\quad +A_k^TP^x_{k+1}B_k+C_k^TP^x_{k+1}D_k\big]+2\bar{X}_k\big\{(\bar{L}^x_k)^T\big[R_k+\bar{R}_k\\
\quad +2 XY_k\big[(L^y_k)^T(R_k+  B_k^TP^x_{k+1}B_k+D_k^TP^x_{k+1}D_k)\\
\quad +F_k^T P^{xy}_{k+1} B_k+ \rho G_k^T P^{xy}_{k+1} D_k \big]+2 \bar{XY}_k\big\{(\bar{L}^y_k)^T\big[R_k+\bar{R}_k\\
\quad +(B_k+\bar{B}_k)^T(P^x_{k+1}+\bar{P}^x_{k+1})(B_k+\bar{B}_k)\\
\quad +(D_k+\bar{D}_k)^T P^x_{k+1}(D_k+\bar{D}_k)\big]+(L^y_k)^T[\bar{R}_k\\
\quad +(B_k+\bar{B}_k)^T(P^x_{k+1}+\bar{P}^x_{k+1})(B_k+\bar{B}_k)-B_k^TP^x_{k+1}B_k\\
\quad +(D_k+\bar{D}_k)^T P^x_{k+1}(D_k+\bar{D}_k)-D_k^TP^x_{k+1}D_k\big]\\
\quad +(F_k+\bar{F}_k)^T (P^{xy}_{k+1}+\bar{P}^{xy}_{k+1})(B_k+\bar{B}_k)-F_k^TP^{xy}_{k+1}B_k\\
\quad +\rho(G_k+\bar{G}_k)^T P^{xy}_{k+1} (D_k+\bar{D}_k)-\rho{G}_k^TP^{xy}_{k+1}D_k\big\}\\
=2(X_k-\bar{X}_k)\big((L^x_k)^T\bar{W}_k^{(1)}+\bar{H}_k^{(1)}\big)\\
\quad +2\bar{X}_k\big((L^x_k+\bar{L}^x_k)^T\bar{W}_k^{(2)}+\bar{H}_k^{(2)}\big)\\
\quad +2(XY_k-\bar{XY}_k)\big((L^y_k)^T\bar{W}_k^{(1)}+\bar{H}_k^{(3)}\big)\\
\quad +2\bar{XY}_k\big((L^y_k+\bar{L}^y_k)^T\bar{W}_k^{(2)}+\bar{H}_k^{(4)}\big),
\end{array}
\end{eqnarray}
\begin{eqnarray}\label{1-barLx-1}
\begin{array}{l}
\ds\frac{\partial \mathfrak{L}_k}{\partial \bar{L}^x_k}= 2\bar{X}_k\big\{(L^x_k+\bar{L}^x_k)^T\big[R_k+\bar{R}_k\\
\quad+(B_k+\bar{B}_k)^T(P^x_{k+1}+\bar{P}^x_{k+1})(B_k+\bar{B}_k)\\ 
\quad +(D_k+\bar{D}_k)^TP^x_{k+1}(D_k+\bar{D}_k)\big]\\
\quad +(A_k+\bar{A}_k)^T(P^x_{k+1}+\bar{P}^x_{k+1})(B_k+\bar{B}_k)\\
\quad +(C_k+\bar{C}_k)^TP^x_{k+1}(D_k+\bar{D}_k)\big\}\\
\quad +2\bar{XY}_k\big\{(L^y_k+\bar{L}^y_k)^T\big[R_k+\bar{R}_k\\
\quad +(B_k+\bar{B}_k)^T(P^x_{k+1}+\bar{P}^x_{k+1})(B_k+\bar{B}_k)\\
\quad +(D_k+\bar{D}_k)^TP^x_{k+1}(D_k+\bar{D}_k)\big]\\
\quad +(F_k+\bar{F}_k)^T(P^{xy}_{k+1}+\bar{P}^{xy}_{k+1})(B_k+\bar{B}_k)\\
\quad +\rho(G_k+\bar{G}_k)^T P^{xy}_{k+1} (D_k+\bar{D}_k)\big\}\\
=2\bar{X}_k\big((L^x_k+\bar{L}^x_k)^T\bar{W}_k^{(2)}+\bar{H}_k^{(2)}\big)\\
\quad +2\bar{XY}_k\big((L^y_k+\bar{L}^y_k)^T\bar{W}_k^{(2)}+\bar{H}_k^{(4)}\Big),
\end{array}
\end{eqnarray}
\begin{eqnarray}\label{1-Ly-1}
\begin{array}{l}
\ds\frac{\partial \mathfrak{L}_k}{\partial L^y_k}=2 Y_k\big[(L^y_k)^T(R_k+B_k^TP^x_{k+1}B_k+D_k^TP^x_{k+1}D_k)\\
\quad +F_k^T P^{xy}_{k+1} B_k+ \rho G_k^T P^{xy}_{k+1} D_k \big]+2\bar{Y}_k\big\{(\bar{L}^y_k)^T\big[R_k+\bar{R}_k\\
\quad +(B_k+\bar{B}_k)^T(P^x_{k+1}+\bar{P}^x_{k+1})(B_k+\bar{B}_k)\\
\quad +(D_k+\bar{D}_k)^TP^x_{k+1}(D_k+\bar{D}_k)\big]+(L^y_k)^T\big[\bar{R}_k\\
\quad +(B_k+\bar{B}_k)^T(P^x_{k+1}+\bar{P}^x_{k+1})(B_k+\bar{B}_k)-B_k^T P^x_{k+1}B_k\\
\quad +(D_k+\bar{D}_k)^TP^x_{k+1}(D_k+\bar{D}_k)^T-D_k^TP^x_{k+1}{D}_k\big]\\
\quad +(F_k+\bar{F}_k)^T(P^{xy}_{k+1}+\bar{P}^{xy}_{k+1})(B_k+\bar{B}_k)-F_k^T P^{xy}_{k+1} B_k\\
\quad +\rho(G_k+\bar{G}_k)^T P^{xy}_{k+1} (D_k+\bar{D}_k)- \rho G_k^T P^{xy}_{k+1}{D}_k \big\}\\
\quad +2(XY_k)^T\big[(L^x_k)^T(R_k+B_k^T P^x_{k+1}B_k+D_k^T P^x_{k+1}D_k)\\
\quad +A^T P^x_{k+1}B_k + C_k^T P^x_{k+1}D_k\big]+2(\bar{XY}_k)^T\big\{(\bar{L}^x_k)^T\big[R_k+\bar{R}_k\\
\quad +(B_k+\bar{B}_k)^T(P^x_{k+1}+\bar{P}^x_{k+1})(B_k+\bar{B}_k)\\
\quad +(D_k+\bar{D}_k)^TP^x_{k+1}(D_k+\bar{D}_k)\big]+(L^x_k)^T\big[\bar{R}_k\\
\quad +(B_k+\bar{B}_k)^T(P^x_{k+1}+\bar{P}^x_{k+1})(B_k+\bar{B}_k)-B_k^TP^x_{k+1}B_k\\
\quad +(D_k+\bar{D}_k)^TP^x_{k+1}(D_k+\bar{D}_k)-{D}_k^TP^x_{k+1}D_k\big]\\
\quad +(A_k+\bar{A}_k)^T(P^x_{k+1}+\bar{P}_{k+1})(B_k+\bar{B}_k)-A^T_kP^x_{k+1}B_k\\
\quad +(C_k+\bar{C}_k)^TP^x_{k+1}(D_k+\bar{D}_k)k-C^T_kP^x_{k+1}D_k\big\}\\
= 2(Y_k-\bar{Y}_k)\big((L^y_k)^T\bar{W}_k^{(1)}+\bar{H}_k^{(3)}\big)\\
\quad +2\bar{Y}_k\big(\big(L^y_k+\bar{L}^y_k\big)^T\bar{W}_k^{(2)}+\bar{H}_k^{(4)}\big)\\
\quad +\ds2(XY_k-\bar{XY}_k)^T\big((L^x_k)^T\bar{W}_k^{(1)}+\bar{H}_k^{(1)}\big)\\
\quad +2(\bar{XY}_k)^T\big((L^x_k+\bar{L}^x_k)^T\bar{W}_k^{(2)}+\bar{H}_k^{(2)}\big),
\end{array}
\end{eqnarray}
\begin{eqnarray}\label{1-barLy-1}
\begin{array}{l}
\ds\frac{\partial \mathfrak{L}_k}{\partial \bar{L}^y_k}=2\bar{Y}_k\{(L^y_k+\bar{L}^y_k)^T[R_k+\bar{R}_k\\
\quad +(B_k+\bar{B}_k)^T(P^x_{k+1}+\bar{P}^x_{k+1})(B_k+\bar{B}_k)\\
\quad +(D_k+\bar{D}_k)^TP^x_{k+1}(D_k+\bar{D}_k)\big]\\
\quad +(F_k+\bar{F}_k)^T (P^{xy}_{k+1}+\bar{P}^{xy}_{k+1})(B_k+\bar{B}_k)\\
\quad +\rho(G_k+\bar{G}_k)^T P^{xy}_{k+1} (D_k+\bar{D}_k)\big\}\\
\quad +2(\bar{XY}_k)^T\{(\bar{L}^x_k+L^x_k)^T\big[R_k+\bar{R}_k\\
\quad +(B_k+\bar{B}_k)^T(P^x_{k+1}+\bar{P}^x_{k+1})(B_k+\bar{B}_k)\\
\quad +(D_k+\bar{D}_k)^TP^x_{k+1}(D_k+\bar{D}_k)\big]\\
\quad +(A_k+\bar{A}_k)^T(P^x_{k+1}+\bar{P}^x_{k+1})(B_k+\bar{B}_k)\\
\quad +(C_k+\bar{C}_k)^TP^x_{k+1}(D_k+\bar{D}_k)\big\}\\
= \ds 2\bar{Y}_k\big((L^y_k+\bar{L}^y_k)^T\bar{W}_k^{(2)}+\bar{H}_k^{(4)}\big)\\
\quad + (\bar{XY}_k)^T\big((L^x_k+\bar{L}^x_k)^T\bar{W}_k^{(2)}+\bar{H}_k^{(2)}\big).
\end{array}
\end{eqnarray}
Here, $\bar{W}_k^{(i)}, i=1,2, \bar{H}_k^{(j)}, j=1,2,3,4$, are defined in (\ref{1-WH-1}).
%
%
Combining (\ref{1-Lx-1})--(\ref{1-barLy-1}), $L^{xo}_k, \bar{L}^{xo}_k, L^{yo}_k$ and $\bar{L}^{yo}_k$ must satisfy
\begin{eqnarray}\label{1-mid-1}
\left\{\begin{array}{l}
(X_k-\bar{X}_k)\big((L^x_k)^T\bar{W}_k^{(1)}+\bar{H}_k^{(1)}\big)\\
\quad +\bar{X}_k\big((L^x_k+\bar{L}^x_k)^T\bar{W}_k^{(2)}+\bar{H}_k^{(2)}\big)\\
\quad +(XY_k-\bar{XY}_k)\big((L^y_k)^T\bar{W}_k^{(1)}+\bar{H}_k^{(3)}\big)\\
\quad +\bar{XY}_k\big((L^y_k+\bar{L}^y_k)^T\bar{W}_k^{(2)}+\bar{H}_k^{(4)}\big)=0,\\
\bar{X}_k\big((L^x_k+\bar{L}^x_k)^T\bar{W}_k^{(2)}+\bar{H}_k^{(2)}\big)\\
\quad +\bar{XY}_k\big((L^y_k+\bar{L}^y_k)^T\bar{W}_k^{(2)}+\bar{H}_k^{(4)}\big)=0,\\
(Y_k-\bar{Y}_k)\big((L^y_k)^T\bar{W}_k^{(1)}+\bar{H}_k^{(3)}\big)\\
\quad +\bar{Y}_k\big((L^y_k+\bar{L}^y_k)^T\bar{W}_k^{(2)}+\bar{H}_k^{(4)}\big)\\
\quad +(XY_k-\bar{XY}_k)^T\big((L^x_k)^T\bar{W}_k^{(1)}+\bar{H}_k^{(1)}\big)\\
\quad +(\bar{XY}_k)^T\big((L^x_k+\bar{L}^x_k)^T\bar{W}_k^{(2)}+\bar{H}_k^{(2)}\big)=0,\\
\bar{Y}_k\big((L^y_k+\bar{L}^y_k)^T\bar{W}_k^{(2)}+\bar{H}_k^{(4)}\big)\\
\quad +(\bar{XY}_k)^T\big((L^x_k+\bar{L}^x_k)^T\bar{W}_k^{(2)}+\bar{H}_k^{(2)}\big)=0.
\end{array}\right.
\end{eqnarray}
Note that (\ref{1-mid-1}) holds for any initial values $X_0-\bar{X}_0$, $\bar{X}_0, XY_0-\bar{XY}_0$, $\bar{XY}_0$, $Y_0-\bar{Y}_0$ and $\bar{Y}_0$. Hence, (\ref{1-mid-1}) reduces to
\begin{eqnarray}
\left\{\begin{array}{l}
(L^x_k)^T\bar{W}_k^{(1)}+\bar{H}_k^{(1)}=0, \\
(L^x_k+\bar{L}^x_k)^T\bar{W}_k^{(2)}+\bar{H}_k^{(2)}=0,\\
(L^y_k)^T\bar{W}_k^{(1)}+\bar{H}_k^{(3)}=0, \\
(L^y_k+\bar{L}^y_k)^T\bar{W}_k^{(2)}+\bar{H}_k^{(4)}=0,\\
\end{array}\right.
\end{eqnarray}
which is obtained by letting coefficients be zero in (\ref{1-mid-1}).
Clearly, we obtain the optimal feedback gains within the class of controls (\ref{1-control-2})
\begin{eqnarray*}\label{1-L-1}
\left\{\begin{array}{l}
L^{xo}_k=-(\bar{W}_k^{(1)})^{-1}(\bar{H}_k^{(1)})^T,\\
\bar{L}^{xo}_k=-(\bar{W}_k^{(2)})^{-1}(\bar{H}_k^{(2)})^T+(\bar{W}_k^{(1)})^{-1}(\bar{H}_k^{(1)})^T,\\
L^{yo}_k=-(\bar{W}_k^{(1)})^{-1}(\bar{H}_k^{(3)})^T,\\
\bar{L}^{yo}_k=-(\bar{W}_k^{(2)})^{-1}(\bar{H}_k^{(4)})^T+(\bar{W}_k^{(1)})^{-1}(\bar{H}_k^{(3)})^T.
\end{array}\right.
\end{eqnarray*}
We now derive the equations that $P^x_k$, $\bar{P}^x_k$, $P^{xy}_k$, $\bar{P}^{xy}_k$,  $P^y_k, \bar{P}^y_k$ satisfy. By (\ref{1-P-1}), we have
\begin{eqnarray*}
\begin{array}{l}
P^x_k= \ds\frac{\partial\mathfrak{L}_k}{\partial {X}_k}\Big{|}_{L^x_k=L_k^{xo}}\\
=Q_k+(L_k^{xo})^TR_kL^{xo}_k+(A_k+B_kL^{xo}_k)^TP^x_{k+1}(A_k\\
\quad +B_kL^{xo}_k)+(C_k+D_kL^{xo}_k)^TP^x_{k+1}(C_k+D_kL^{xo}_k),
\end{array}%
\end{eqnarray*}
which are (\ref{1-Px-1})--(\ref{1-barPy-1}). 
The final step is to assure that for any $k \in \bar{\mathbb{N}}$, $P^x_k,~P^x_k+\bar{P}^x_k\geq 0$. We prove this by backward induction. Clearly, $P^x_N,~P^x_N+\bar{P}^x_N\geq 0$ by definition. For $k=N-1$, $P^x_{N-1}$ is positive semi-definite and
\begin{eqnarray*}
\begin{array}{l}
P^x_{N-1}+\bar{P}^x_{N-1}=Q_{N-1}+\bar{Q}_{N-1}+(L_{N-1}^{xo}+\bar{L}_{N-1}^{xo})^TR_{N-1}\\
\quad \cdot(L_{N-1}^{xo}+\bar{L}_{N-1}^{xo})+\big[A_{N-1}+\bar{A}_{N-1}+(B_{N-1}+\bar{B}_{N-1})\\
\quad \cdot(L^{xo}_{N-1}+\bar{L}^{xo}_{N-1})\big]^T(P^x_N+\bar{P}^x_N)\big[A_{N-1}+\bar{A}_{N-1}+(B_{N-1}\\
\quad +\bar{B}_{N-1})(L^{xo}_{N-1}+\bar{L}^{xo}_{N-1})\big]+\big[C_{N-1}+\bar{C}_{N-1}+(D_{N-1}\\
\quad +\bar{D}_{N-1})(L^{xo}_{N-1}+\bar{L}^{xo}_{N-1})\big] ^T P^x_N \big[C_{N-1}+\bar{C}_{N-1}+(D_{N-1}\\
\quad +\bar{D}_{N-1})(L^{xo}_{N-1}+\bar{L}^{xo}_{N-1})\big]\geq 0.
\end{array}
\end{eqnarray*}
After simple calculation, we have
\begin{eqnarray*}
%
\end{eqnarray*}
We can easily induce that $P^x_k,~P^x_k+\bar{P}^x_k\geq 0$. Note that
\begin{eqnarray*}
%
\end{eqnarray*}
This completes the proof. 
\hfill $\square$

\subsection*{Proof of Lemma \ref{1-lemma-1}.}

\emph{Proof}. We now reformulate equations (\ref{2-Operator-1}) and (\ref{2-Operator-2}) in terms of $\mathbb{E}$, $I-\mathbb{E}$ as follows
\begin{eqnarray*}
\left\{%
\right.
\end{eqnarray}
To proceed, we need $\Theta_{2,N-1}^{-1}$. Under condition (\ref{2-condition-1}), $\mathcal{P}^x_N\geq0 $ and  
\begin{eqnarray}\label{2-R-1}
%
\end{eqnarray}
Here $\vert \cdot \vert_m$ denotes the norm in $\mathbb{R}^m$; $\lambda^{(N-1)}_1, \lambda^{(N-1)}_2$ are the smallest eigenvalues of matrices $R_{N-1}$ and $R_{N-1}+\bar{R}_{N-1}$, respectively, and $\lambda^{(N-1)}=\mbox{min}\left\{ \lambda^{(N-1)}_1, \lambda^{(N-1)}_2 \right\}$; $\Vert \cdot \Vert_m$ is the norm induced by inner product in $\mathcal{U}_{N-1}$. Hence, $\Theta_{2,N-1}$ must be positive definite and self-adjoint. So far, (\ref{2-P-5})--(\ref{2-cost-3}) are established. Furthermore, the technique of operator pseudo-inverse is used to compute $\Theta_{2,N-1}^{-1}$\cite{Elliott2013,Beutler1965}.

Clearly,
$(I-\mathbb{E})(I-\mathbb{E})^{\dagger}= 
\bigg(\, \begin{matrix} 
I &0 \\[-3mm]
0 &0 
\end{matrix}\, \bigg)
$
, $\mathbb{E}\mathbb{E}^{\dagger}=
\bigg(\, \begin{matrix} 
0 & 0 \\[-3mm]
0 & I 
\end{matrix}\, \bigg)$
and
\begin{eqnarray*}
\left\{\begin{array}{l}
(I-\mathbb{E})(I-\mathbb{E})^{\dagger}+\mathbb{E}\mathbb{E}^{\dagger}=I,\\
\mathbb{E}(I-\mathbb{E})^{\dagger}=0, (I-\mathbb{E})\mathbb{E}^{\dagger}=0.\\
\end{array}\right.
\end{eqnarray*}
From (\ref{2-Th-2}) and (\ref{2-La-2}), we have
\begin{eqnarray}\label{2-T-2}
\begin{array}{l}
\Theta_{2,N-1}^{-1}=(\mathcal{R}_{N-1}+ \mathcal{B}_{N-1}^T\mathcal{P}^x_N\mathcal{B}_{N-1}+\mathcal{D}_{N-1}^T\mathcal{P}^x_N\mathcal{D}_{N-1})^{-1}\\
=(I-\mathbb{E})^{\dagger} (W_{N-1}^{(1)})^{-1} (I-\mathbb{E}^*)^{\dagger}+\mathbb{E}^{\dagger} (W_{N-1}^{(2)})^{-1} (\mathbb{E}^*)^{\dagger}.
\end{array}
\end{eqnarray}
In fact,
\begin{eqnarray*}
\begin{array}{l}
\big[(I-\mathbb{E}^*) W_{N-1}^{(1)} (I-\mathbb{E})+\mathbb{E}^* W_{N-1}^{(2)} \mathbb{E}\big]\\
\cdot\big[(I-\mathbb{E})^{\dagger} (W_{N-1}^{(1)})^{-1} (I-\mathbb{E}^*)^{\dagger}+\mathbb{E}^{\dagger} (W_{N-1}^{(2)})^{-1} (\mathbb{E}^*)^{\dagger}\big]\\
=(I-\mathbb{E}^*)(I-\mathbb{E})(I-\mathbb{E})^{\dagger}(I-\mathbb{E}^*)^{\dagger}+\mathbb{E}^*\mathbb{E}\mathbb{E}^{\dagger}(\mathbb{E}^*)^{\dagger}\\
=I.
\end{array}
\end{eqnarray*}
By simple calculation,
\begin{eqnarray}\label{2-Px-4}
\begin{array}{l}
\mathcal{P}^x_{N-1}= \Lambda_{1,N-1}-\Theta_{1,N-1}\Theta_{2,N-1}^{-1}\Theta_{1,N-1}^*\\
=(I-\mathbb{E}^*)\big[ Q_{N-1}+A_{N-1}^T Q_N A_{N-1}+C_{N-1}^T Q_N C_{N-1}\\
\quad - {H}_{N-1}^{(1)}({W}_{N-1}^{(1)})^{-1}({H}_{N-1}^{(1)})^T \big](I-\mathbb{E})+\mathbb{E}^*\big[ Q_{N-1}\\
\quad +\bar{Q}_{N-1}+(A_{N-1}+\bar{A}_{N-1})^T (Q_N+\bar{Q}_N) (A_{N-1}\\
\quad +\bar{A}_{N-1})+(C_{N-1}+\bar{C}_{N-1})^T Q_N (C_{N-1}+\bar{C}_{N-1})\\
\quad -{H}_{N-1}^{(2)}({W}_{N-1}^{(2)})^{-1}({H}_{N-1}^{(2)})^T \big]\mathbb{E},
\end{array}
\end{eqnarray}
\begin{eqnarray}\label{2-Pxy-4}
\begin{array}{l}
\mathcal{P}^{xy}_{N-1}=\Lambda_{2,N-1}-\Theta_{3,N-1}\Theta_{2,N-1}^{-1}\Theta_{1,N-1}^*\\
= (I-\mathbb{E}^*)\big[-Q_{N-1}+F_{N-1}^T Q_{N} A_{N-1}+\rho G_{N-1}^T Q_{N} C_{N-1}\\
\quad -{H}_{N-1}^{(3)}({W}_{N-1}^{(1)})^{-1}({H}_{N-1}^{(1)})^T\big](I-\mathbb{E})+\mathbb{E}^*\big[-Q_{N-1}\\
\quad -\bar{Q}_{N-1}+(F_{N-1}+\bar{F}_{N-1})^T (Q_N+\bar{Q}_N) (A_{N-1}\\
\quad +\bar{A}_{N-1})+\rho(G_{N-1}+\bar{G}_{N-1})^T Q_{N} (C_{N-1}+\bar{C}_{N-1})\\
\quad -{H}_{N-1}^{(4)}({W}_{N-1}^{(2)})^{-1}({H}_{N-1}^{(2)})^T \big]\mathbb{E},
\end{array}
\end{eqnarray}
\begin{eqnarray}\label{2-Py-4}
\begin{array}{l}
\mathcal{P}^y_{N-1}=\Lambda_{3,N-1}-\Theta_{3,N-1}\Theta_{2,N-1}^{-1}\Theta_{3,N-1}^*\\
= (I-\mathbb{E}^*)\big[ Q_{N-1}+F_{N-1}^T Q_{N} F_{N-1}+G_{N-1}^T Q_{N} G_{N-1}\\
\quad -{H}_{N-1}^{(3)}({W}_{N-1}^{(1)})^{-1}({H}_{N-1}^{(3)})^T \big](I-\mathbb{E})+\mathbb{E}^*\big[ Q_{N-1}\\
\quad +\bar{Q}_{N-1}+(F_{N-1}+\bar{F}_{N-1})^T (Q_N+\bar{Q}_N) (F_{N-1}\\
\quad +\bar{F}_{N-1})+(G_{N-1}+\bar{G}_{N-1})^T Q_{N} (G_{N-1}+\bar{G}_{N-1})\\
\quad -{H}_{N-1}^{(4)}({W}_{N-1}^{(2)})^{-1}({H}_{N-1}^{(4)})^T\big]\mathbb{E},
\end{array}
\end{eqnarray}
which are (\ref{2-P-2}). Hence, we can easily derive that $\mathcal{P}^x_{N-1}$ is positive definite. This completes the proof.
\hfill $\square$

\subsection*{Proof of Theorem \ref{2-theorem-1}.}

\emph{Proof}. Suppose that $\mathcal{P}^x_{k+1}\geq0$. By combining (\ref{2-Px-1})--(\ref{2-Py-1}), we get
\begin{eqnarray*}
\begin{array}{l}
\langle \mathcal{P}^x_{k+1}x_{k+1},x_{k+1}\rangle =\langle \mathcal{P}^x_{k+1}(\mathcal{A}_kx_k+\mathcal{B}_ku_k), (\mathcal{A}_kx_k+\mathcal{B}_ku_k) \rangle\\
\qquad +\langle \mathcal{P}^x_{k+1}(\mathcal{C}_kx_k+\mathcal{D}_ku_k), (\mathcal{C}_kx_k+\mathcal{D}_ku_k) \rangle,\\
\langle \mathcal{P}^{xy}_{k+1}x_{k+1},y_{k+1}\rangle =\langle \mathcal{P}^{xy}_{k+1}(\mathcal{A}_kx_k+\mathcal{B}_ku_k), \mathcal{F}_ky_k \rangle\\
\qquad +\rho\langle \mathcal{P}^{xy}_{k+1}(\mathcal{C}_kx_k+\mathcal{D}_ku_k), \mathcal{G}_ky_k \rangle,\\
\langle \mathcal{P}^y_{k+1}y_{k+1},y_{k+1}\rangle =\langle \mathcal{P}^y_{k+1}\mathcal{F}_ky_k, \mathcal{F}_ky_k \rangle +\langle \mathcal{P}^y_{k+1}\mathcal{G}_ky_k, \mathcal{G}_ky_k, \rangle.\\
\end{array}
\end{eqnarray*}
And we may isolate the following term from (\ref{2-cost-1}) in terms of $k$:
\begin{eqnarray*}
\begin{array}{l}
\langle \mathcal{Q}_{k}x_{k}, x_{k}\rangle -2\langle \mathcal{Q}_{k}x_{k},y_{k} \rangle + \langle \mathcal{Q}_{k}y_{k},y_{k} \rangle +\langle \mathcal{R}_{k}u_{k}, u_{k}\rangle\\ 
+\langle \mathcal{P}^x_{k+1}x_{k+1}, x_{k+1}\rangle-2\langle \mathcal{P}^{xy}_{k+1}x_{k+1},y_{k+1} \rangle + \langle \mathcal{P}^y_{k+1}y_{k+1},y_{k+1} \rangle\\
=2\langle \Theta_{1,k}u_{k},x_{k} \rangle +\langle \Theta_{2,k}u_{k},u_{k} \rangle +2\langle \Theta_{3,k}u_{k},y_{k} \rangle\\
\quad +\langle \Lambda_{1,k}x_{k},x_{k} \rangle +2\langle \Lambda_{2,k}x_{k},y_{k} \rangle +\langle \Lambda_{3,k}y_{k},y_{k} \rangle\\
=\big\langle (\Lambda_{1,k}-\Theta_{1,k}\Theta_{2,k}^{-1}\Theta_{1,k}^*)x_{k},x_{k} \big\rangle +2\big\langle (\Lambda_{2,k}-\Theta_{3,k}\Theta_{2,k}^{-1}\\
\quad \cdot\Theta_{1,k}^*)x_{k},y_{k} \big\rangle +\big\langle (\Lambda_{3,k}-\Theta_{3,k}\Theta_{2,k}^{-1}\Theta_{3,k}^*)y_{k},y_{k} \big\rangle\\
\quad +\big\langle \Theta_{2,k}(u_{k}+\Theta_{2,k}^{-1}\Theta_{1,k}^*x_{k}+\Theta_{2,k}^{-1}\Theta_{3,k}^*y_{k}),\\
\quad (u_{k}+\Theta_{2,k}^{-1}\Theta_{1,k}^*x_{k}+\Theta_{2,k}^{-1}\Theta_{3,k}^*y_{k}) \big\rangle,
\end{array}
\end{eqnarray*}
where
\begin{eqnarray}\label{2-TL-3}
\left\{
\begin{array}{rl}
\Theta_{1,k}=&\mathcal{A}_{k}^T\mathcal{P}^x_{k+1}\mathcal{B}_{k}+\mathcal{C}_{k}^T\mathcal{P}^x_{k+1}\mathcal{D}_{k}\\
=&(I-\mathbb{E}^*) {H}_{k}^{(1)} (I-\mathbb{E})+\mathbb{E}^* {H}_{k}^{(2)} \mathbb{E},\\
\Theta_{2,k}=&\mathcal{R}_{k}+\mathcal{B}_{k}^T\mathcal{P}^x_{k+1}\mathcal{B}_{k}+\mathcal{D}_{k}^T\mathcal{P}^x_{k+1}\mathcal{D}_{k}\\
=&(I-\mathbb{E}^*) {W}_{k}^{(1)} (I-\mathbb{E})+\mathbb{E}^* {W}_{k}^{(2)} \mathbb{E},\\
\Theta_{3,k}=&\mathcal{F}_{k}^T\mathcal{P}^{xy}_{k+1}\mathcal{B}_{k}+\rho\mathcal{G}_{k}^T\mathcal{P}^{xy}_{k+1}\mathcal{D}_{k}\\
=&(I-\mathbb{E}^*) {H}_{k}^{(3)} (I-\mathbb{E})+\mathbb{E}^* {H}_{k}^{(4)} \mathbb{E},
\end{array}
\right.
\end{eqnarray}
\begin{eqnarray}\label{2-TL-3}
\left\{\begin{array}{l}
\Lambda_{1,k}=\mathcal{Q}_{k}+\mathcal{A}_{k}^T\mathcal{P}^x_{k+1}\mathcal{A}_{k}+\mathcal{C}_{k}^T\mathcal{P}^x_{k+1}\mathcal{C}_{k}\\
\quad =(I-\mathbb{E}^*) \big(Q_{k}+A_{k}^T S^x_{k+1} A_{k}+C_{k}^T S^x_{k+1} C_{k}\big) (I-\mathbb{E})\\
\qquad +\mathbb{E}^* \big[Q_{k}+\bar{Q}_{k}+(A_{k}+\bar{A}_k)^T T^x_{k+1} (A_{k}+\bar{A}_{k})\\
\qquad +(C_{k}+\bar{C}_{k})^T S^x_{k+1} (C_{k}+\bar{C}_{k})\big] \mathbb{E},\\
\Lambda_{2,k}=-\mathcal{Q}_{k}+\mathcal{F}_{k}^T\mathcal{P}^{xy}_{k+1}\mathcal{A}_{k}+\rho\mathcal{G}_{k}^T\mathcal{P}^{xy}_{k+1}\mathcal{C}_{k}\\
\quad =(I-\mathbb{E}^*) \big(-Q_{k}+F_{k}^T S^{xy}_{k+1} A_{k}+\rho G_{k}^T S^{xy}_{k+1} C_{k}\big) (I-\mathbb{E})\\
\qquad +\mathbb{E}^* \big[-Q_{k}-\bar{Q}_{k}+(F_{k}+\bar{F}_{k})^T T^{xy}_{k+1} (A_{k}+\bar{A}_{k})\\
\qquad +\rho(G_{k}+\bar{G}_{k})^T S^{xy}_{k+1} (C_{k}+\bar{C}_{k})\big] \mathbb{E},\\
\Lambda_{3,k}=\mathcal{Q}_{k}+\mathcal{F}_{k}^T\mathcal{P}^y_{k+1}\mathcal{F}_{k}+\mathcal{G}_{k}^T\mathcal{P}^y_{k+1}\mathcal{G}_{k}\\
\quad =(I-\mathbb{E}^*) \big(Q_{k}+F_{k}^T S^y_{k+1} F_{k}+G_{k}^T S^y_{k+1} G_{k}\big) (I-\mathbb{E})\\
\qquad +\mathbb{E}^* \big[Q_{k}+\bar{Q}_{k}+(F_{k}+\bar{F}_{k})^T T^y_{k+1}(F_{k}+\bar{F}_{k})\\
\qquad +(G_{k}+\bar{G}_{k})^T S^y_{k+1} (G_{k}+\bar{G}_{k})\big] \mathbb{E}.
\end{array}\right.
\end{eqnarray}
Similar to (\ref{2-R-1}), we have a positive-definite and self-adjoint $\Theta_{2,k}$ under condition (\ref{2-condition-5}) and hence
\begin{eqnarray*}
\begin{array}{rl}
\Theta_{2,k}^{-1} = (I-\mathbb{E})^{\dagger} (W_{k}^{(1)})^{-1} (I-\mathbb{E}^*)^{\dagger}+\mathbb{E}^{\dagger} (W_{k}^{(2)})^{-1}(\mathbb{E}^*)^{\dagger}.
\end{array}
\end{eqnarray*}
Let
\begin{eqnarray}\label{2-P-3}
\left\{\begin{array}{rl}
\mathcal{P}^x_{k}=&\Lambda_{1,k}-\Theta_{1,k}\Theta_{2,k}^{-1}\Theta_{1,k}^*\\
=&(I-\mathbb{E}^*)S^x_{k}(I-\mathbb{E})+\mathbb{E}^*T^x_{k}\mathbb{E},\\
\mathcal{P}^{xy}_{k}=&\Lambda_{2,k}-\Theta_{3,k}\Theta_{2,k}^{-1}\Theta_{1,k}^*\\
=&(I-\mathbb{E}^*)S^{xy}_{k}(I-\mathbb{E})+\mathbb{E}^*T^{xy}_{k}\mathbb{E},\\
\mathcal{P}^y_{k}=&\Lambda_{3,k}-\Theta_{3,k}\Theta_{2,k}^{-1}\Theta_{3,k}^*\\
=&(I-\mathbb{E}^*)S^y_{k}(I-\mathbb{E})+\mathbb{E}^*T^y_{k}\mathbb{E}.
\end{array}\right.
\end{eqnarray}
Then, in a backward recursion,
\begin{eqnarray}\label{2-cost-6}
\begin{array}{l}
J_k^N(x_k, y_k, u^o|_{\{k, k+1,\cdots,N-1\}})\\
=\langle \mathcal{P}^x_{k}x_{k},x_{k} \rangle +2\langle \mathcal{P}^{xy}_{k}x_{k},y_{k} \rangle +\langle \mathcal{P}^y_{k}x_{k},x_{k} \rangle\\
\quad +\ds\sum_{i=k}^{N-1}\big\langle \Theta_{2,i}(u_{i}+\Theta_{2,i}^{-1}\Theta_{1,i}^*x_{i}+\Theta_{2,i}^{-1}\Theta_{3,i}^*y_{i}),\\
\qquad (u_{i}+\Theta_{2,i}^{-1}\Theta_{1,i}^*x_{i}+\Theta_{2,i}^{-1}\Theta_{3,i}^*y_{i}) \big\rangle.
\end{array}
\end{eqnarray}
We can prove $\mathcal{P}^x_{k}\geq0$ by induction. Consequently,
\begin{eqnarray*}
\begin{array}{l}
J(\zeta^x,\zeta^y,u^o)=\langle \mathcal{P}^x_{0}x_{0},x_{0} \rangle +2\langle \mathcal{P}^{xy}_{0}x_{0},y_{0} \rangle +\langle \mathcal{P}^y_{0}x_{0},x_{0} \rangle\\
\quad +\ds\sum_{k=0}^{N-1}\big\langle \Theta_{2,k}(u_{k}+\Theta_{2,k}^{-1}\Theta_{1,k}^*x_{k}+\Theta_{2,k}^{-1}\Theta_{3,k}^*y_{k}),\\
\qquad (u_{k}+\Theta_{2,k}^{-1}\Theta_{1,k}^*x_{k}+\Theta_{2,k}^{-1}\Theta_{3,k}^*y_{k}) \big\rangle,
\end{array}
\end{eqnarray*}
and the optimal control
\begin{eqnarray*}
\begin{array}{rl}
u_k^* =& -\Theta_{2,k}^{-1}\Theta_{1,k}^*x_{k}-\Theta_{2,k}^{-1}\Theta_{3,k}^*y_{k}\\
=&-(I-\mathbb{E})^{\dagger}({W}_k^{(1)})^{-1}({H}_k^{(1)})^T(x_k-\mathbb{E}x_k)\\
& -\mathbb{E}^{\dagger}({W}_k^{(2)})^{-1}({H}_k^{(2)})^T\mathbb{E}x_k\\
& -(I-\mathbb{E})^{\dagger}({W}_k^{(1)})^{-1}({H}_k^{(3)})^T(y_k-\mathbb{E}y_k)\\
& -\mathbb{E}^{\dagger}({W}_k^{(2)})^{-1}({H}_k^{(4)})^T\mathbb{E}y_k, ~~k\in\mathbb{N},
\end{array}
\end{eqnarray*}
which is (\ref{1-control-2}) by computing $u_k^*=(I+\mathbb{E})u_k^*+\mathbb{E}u_k^*$. This completes the proof.
\hfill $\square$

\end{document}